\documentclass[hidelinks,onefignum,onetabnum]{siamart220329}

\usepackage{lipsum}
\usepackage{amsfonts}
\usepackage{graphicx}
\usepackage{epstopdf}
\usepackage{algorithmic}
\ifpdf
  \DeclareGraphicsExtensions{.eps,.pdf,.png,.jpg}
\else
  \DeclareGraphicsExtensions{.eps}
\fi

\newsiamremark{remark}{Remark}
\newsiamremark{hypothesis}{Hypothesis}
\crefname{hypothesis}{Hypothesis}{Hypotheses}
\newsiamthm{claim}{Claim}

\headers{Globally divergence-free entropy stable DG method for MHD}{Yuchang Liu, Wei Guo, Yan Jiang, and Mengping Zhang}

\title{A globally divergence-free entropy stable nodal DG method for conservative ideal MHD equations\thanks{Submitted to the editors DATE.}}

\author{Yuchang Liu\thanks{School of Mathematical Sciences,
         University of Science and Technology of China,
         Hefei, Anhui 230026, P.R. China.
  (\email{lissandra@mail.ustc.edu.cn}).}
\and Wei Guo\thanks{Department of Mathematics and Statistics, Texas Tech University, Lubbock, TX, 70409, USA. 
  (\email{weimath.guo@ttu.edu}). Research supported by NSF grant NSF-DMS-2111383, Air Force Office of Scientific Research FA9550-18-1-0257. }
\and Yan Jiang\thanks{ School of Mathematical Sciences,
         University of Science and Technology of China, Hefei, 
         Anhui 230026, P.R. China.  
         (\email{jiangy@ustc.edu.cn}). 
         Research supported by NSFC grant 12271499. }
\and Mengping Zhang\thanks{ School of Mathematical Sciences, University of Science and Technology of China, Hefei,
         Anhui 230026, P.R. China.  
         (\email{mpzhang@ustc.edu.cn}). }  
}

\usepackage{amsopn}

\ifpdf
\hypersetup{
  pdftitle={A globally divergence-free entropy stable nodal DG method for conservative ideal MHD equations},
  pdfauthor={Yuchang Liu, Wei Guo, Yan Jiang, and Mengping Zhang}
}
\fi

\allowdisplaybreaks
\usepackage{subfigure}

\begin{document}

\maketitle

\begin{abstract}
We propose an arbitrarily high-order globally divergence-free entropy stable nodal discontinuous Galerkin (DG) method to directly solve the conservative form of the ideal MHD equations using appropriate quadrature rules. The method ensures a globally divergence-free magnetic field by updating it at interfaces with a constraint-preserving formulation \cite{balsara2021globally} and employing a novel least-squares reconstruction technique. Leveraging this property, the semi-discrete nodal DG scheme is proven to be entropy stable. To handle the problems with strong shocks, we introduce a novel limiting strategy that suppresses unphysical oscillations while preserving the globally divergence-free property. Numerical experiments verify the accuracy and efficacy of our method.
\end{abstract}



\section{Introduction}
\label{sec1}

In this paper, we focus on simulating the compressible ideal magnetohydrodynamic (MHD) equations, which describe the dynamics of perfectly conducting fluids and have a wide range of applications in space weather forecasting, astrophysics, plasma physics, among others. The ideal MHD equations are derived by combing the Euler equations from fluid dynamics and pre-Maxwell's equations from electromagnetism and are formulated as a set of nonlinear hyperbolic conservation laws, along with a divergence-free constraint on the magnetic field. Solving the MHD equations accurately and reliably presents significant challenges. As with other nonlinear conservation laws, the desired entropy solution of the MHD equations is known to develop discontinuities and other complex solution features in finite time and satisfy the entropy inequality, aligning with the second law of thermodynamics. In addition, if the magnetic field is initially divergence-free, it will remain so for all time, reflecting the absence of magnetic monopoles in nature. Failure to preserve these physical constraints on the discrete level may lead to unphysical artifacts and instability of MHD simulations \cite{balsara2004comparison}. 
In this work, we aim for a globally divergence-free shock capture MHD solver that further fulfills a discrete entropy inequality, termed entropy stability.

In recent years, various effective high-order methods have been proposed to solve the MHD equations, including finite difference methods \cite{abd2003finite, christlieb2014finite, minoshima2019high}, finite volume methods \cite{felker2018fourth}, spectral methods \cite{schneider2011pseudo, luo2014simulation}, finite element methods \cite{hu2017stable}, and discontinuous Galerkin (DG) methods \cite{li2005locally,liu2018entropy,balsara2021globally}. In this work, we utilize the DG discretization \cite{reed1973triangular, cockburn1989tvb, cockburn1989tvb3, cockburn1990runge2, cockburn1991runge, cockburn1998runge}, leveraging its well-known advantages, such as high order accuracy, good conservation properties, the ability to handle hanging nodes, and high parallel efficiency. Furthermore, its flexibility enables the design of effective techniques to enforce the desired constraints, specifically the divergence-free property and entropy stability. 

To achieve the divergence-free property, several effective approaches have been proposed in the literature, including the 8-wave formulation \cite{powell1999solution}, the projection methods \cite{toth2000b}, hyperbolic divergence-cleaning technique \cite{dedner2002hyperbolic}, and the constrained transport (CT) method \cite{evans1988simulation, balsara2004second, christlieb2014finite, yu2020free}.  In the context of DG discretization, a locally divergence-free DG method was developed in \cite{li2005locally}, which employs a piecewise divergence-free polynomial vector space  to approximate the magnetic field. In \cite{li2011central}, the globally divergence-free central DG methods involving two copies of solutions were formulated for solving the two-dimensional ideal MHD equation. In such methods, the magnetic field is approximated within the exactly divergence-free $H(div)$-conforming Brezzi–Douglas–Marini  finite element space \cite{brezzi1985two}. In particular, this method defines the magnetic field components on the interface as univariate polynomials, utilizing them at each time step to reconstruct the magnetic field, thereby ensuring the magnetic field is divergence-free internally and continuous in the normal direction. In \cite{fu2018globally}, Fu \textit{et al.} extended it to DG method by using multidimensional Riemann solver at vertices. The magnetic field reconstruction requirements of the above two methods do not explicit contain the divergence-free constraint, and the divergence-free property is guaranteed by a theorem, in which the scheme on the boundary and internal must be closely related. However, for shock problems, they fail in designing a limiter for the magnetic field at interfaces since the limiter will destroy the theorem. Balsara \textit{et al.} made an improvement in \cite{balsara2021globally}. Their reconstruction requirements explicit contain the divergence-free constraint, hence the reconstruction is compatible with the limiter.  

Meanwhile, the development of entropy stable numerical methods for hyperbolic conservation laws has attracted significant research interest over the past two decades, as they are crucial for capturing the physically relevant solution that obeys the second law of thermodynamics.
For instance, the famous high order entropy stable finite volume ENO scheme, termed the TeCNO scheme, was developed in \cite{fjordholm2012arbitrarily}, leveraging fundamental techniques such as entropy conservative and entropy stable fluxes \cite{tadmor1987numerical, tadmor2003entropy} and the sign property of ENO reconstruction \cite{cockburn1998essentially}. For DG discretization, two primary approaches have been developed to ensure entropy stability. The first approach utilizes the summation-by-parts (SBP) methodology \cite{gassner2013skew}. In \cite{chen2017entropy}, Chen and Shu constructed the entropy stable nodal DG method by introducing SBP operators on Gauss-Lobatto points in conjunction. In \cite{chan2018discretely, chan2019skew}, it was extended to modal DG formulation by designing the hybridized SBP operators that combine the volume and surface quadrature nodes. The second approach aims to explicitly control the production of entropy within each element by adding an artificial term, see \cite{abgrall2018general, gaburro2023high, liu2024non}. 

In \cite{godunov1972symmetric}, Godunov pointed out that the ideal MHD equations with the divergence-free constraint are not symmetrizable, unless an additional source term is included, which vanishes if the magnetic field is divergence-free. A symmetrizable hyperbolic system is guaranteed to possess an entropy function. Hence, entropy analysis and the divergence-free condition are closely related at both the PDE and discrete levels. In \cite{chandrashekar2016entropy}, Chandrashekar et al. proposed an entropy stable finite volume entropy scheme for this modified symmetrizable form of MHD by introducing a novel entropy conservative flux.  In \cite{liu2018entropy}, Liu \emph{et al.} extended this approach to construct an entropy stable nodal DG method by carefully handling the added source term. However, since both methods do not satisfy the divergence-free constraint at the discrete level, they are inherently nonconservative by discretizing the modified MHD equations. This can lead to numerical instability and nonphysical features in the approximated solutions \cite{toth2000b, chandrashekar2016entropy}.  In this work, we follow the entropy stable nodal DG framework \cite{liu2018entropy} and integrate the globally divergence-free technique \cite{balsara2021globally}. This allows us to directly simulate the conservative form of the ideal MHD equations, as opposed to the existing approaches, ensuring both entropy stability and globally divergence-free conditions. For problems with strong shocks, we develop a scaling limiter to suppress unphysical oscillations, inspired by the key idea of the oscillation-free DG method and subsequent works \cite{lu2021oscillation, peng2024oedg, wei2024jump}. This limiter preserves both entropy stability and the divergence-free property.

The contributions of this work are summarized as follows:
\begin{itemize}

\item[1.] The magnetic field is globally divergence-free.

\item[2.] The semi-discrete scheme in each cell is entropy stable.

\item[3.] The framework conserves locally the mass, momentum and magnetic field, while total energy conservation is not maintained. Such a conservation error depends on the on the discrepancy between the magnetic field before and after the reconstruction step. 

\end{itemize}

The rest of this paper is organized as follows: Section \ref{sec2} introduces the ideal MHD equations and the associated entropy condition. Section \ref{sec3} presents the proposed globally divergence-free entropy stable (ES-GDF) nodal DG scheme for the conservative form of MHD, along with a theoretical analysis. Section \ref{sec4} describes a consistent limiting strategy that preserves both entropy stability and the globally divergence-free property, as well as the fully discrete formulation. In Section \ref{sec5}, various numerical examples are provided to demonstrate the efficiency and effectiveness of the proposed scheme.  Section \ref{sec6} concludes with a discussion of the findings.

\section{Ideal MHD equations}
\label{sec2}
\setcounter{equation}{0}
\setcounter{figure}{0}

\subsection{Governing equations}

In this paper, we consider the two-dimensional ideal MHD equations in Cartesian coordinates, formulated as a system of conservation laws
\begin{equation}
\label{eq:MHD}
\frac{\partial \mathbf U}{\partial t} + \frac{\partial \mathbf F(\mathbf U)}{\partial x} + \frac{\partial \mathbf G(\mathbf U)}{\partial y} = \mathbf{0},
\end{equation}
where $
\mathbf{U}=\left( \rho ,u_x,u_y,u_z,\mathcal{E} ,B_x,B_y,B_z \right) ^T,
$ and
$$
\mathbf{F}\left( \mathbf{U} \right) =\left[ \begin{array}{l}
	\rho u_x\\
	\rho u_{x}^{2}+p^{\star}-B_{x}^{2}\\
	\rho u_xu_y-B_xB_y\\
	\rho u_xu_z-B_xB_z\\
	u_x\left( \mathcal{E} +p^{\star} \right) -B_x\left( \mathbf{u}\cdot \mathbf{B} \right)\\
	0\\
	u_xB_y-u_yB_x\\
	u_xB_z-u_zB_x\\
\end{array} \right],\quad \mathbf{G}\left( \mathbf{U} \right) =\left[ \begin{array}{l}
	\rho u_y\\
	\rho u_xu_y-B_xB_y\\
	\rho u_{y}^{2}+p^{\star}-B_{y}^{2}\\
	\rho u_yu_z-B_yB_z\\
	u_y\left( \mathcal{E} +p^{\star} \right) -B_y\left( \mathbf{u}\cdot \mathbf{B} \right)\\
	u_yB_x-u_xB_y\\
	0\\
	u_yB_z-u_zB_y\\
\end{array} \right], 
$$
with 
\begin{equation}\label{eq:EOS}\mathcal E=\frac{p}{\gamma - 1} + \frac{1}{2}\rho\left\|\mathbf u\right\|^2+\frac{1}{2}\left\|\mathbf B\right\|^2.\end{equation}
Here, $\rho$ represents the mass density, $\rho\mathbf{u}$ is the momentum, $\mathcal{E}$ denotes the total energy, $p$ is the hydrodynamic pressure, $\mathbf{B}$ represents the magnetic field, $\gamma = 5/3$ is the ideal gas constant, $p^\star = p + \left\|\mathbf{B}\right\|^2/2$ is the total pressure, and $\|\cdot\|$ denotes the Euclidean vector norm. Taking divergence of the magnetic field equation yields
$
\frac{\partial \left( \nabla \cdot \mathbf{B} \right)}{\partial t}=0,
$
 indicating
 $\nabla\cdot\mathbf B(\mathbf x,t) = \nabla\cdot\mathbf B(\mathbf x,0).$
Consequently, if the divergence of the magnetic field is initially zero, it will remain zero for all time, i.e. 
\begin{equation}\label{eq:divfree0}
\nabla\cdot\mathbf B = 0.
\end{equation}
This is called the \textit{divergence-free property}. 
Preservation of the discrete divergence-free property in MHD simulations is critical to prevent unphysical artifacts and numerical instabilities.

\subsection{Entropy analysis for the ideal MHD equations}

It is well-known that weak solutions to conservation laws are not unique. Imposing an entropy condition ensures uniqueness for scalar problems \cite{godlewski1991hyperbolic}, although for systems, the uniqueness remains an open problem. 
Nevertheless, the entropy condition is a valuable design principle that enhances numerical stability and ensures compliance with the second law of thermodynamics.

\begin{definition}\label{def:entropy}
A convex scalar function $\mathcal U(\mathbf U)$ is called an entropy function of system \eqref{eq:MHD}, if there exist scalar entropy fluxes $\mathcal F(\mathbf U)$ and $\mathcal G(\mathbf U)$ such that 
\begin{equation}\label{eq:entropy}
\mathcal F'\left( \mathbf{U} \right) =\mathcal U'\left( \mathbf{U} \right) \mathbf{F}'\left( \mathbf{U} \right) ,\quad \mathcal G'\left( \mathbf{U} \right) =\mathcal U'\left( \mathbf{U} \right) \mathbf{G}'\left( \mathbf{U} \right) .
\end{equation}
The functions $(\mathcal U,\mathcal F,\mathcal G)$ are called an entropy pair. 
Here, $\mathcal U'(\mathbf U)$, $\mathcal F'(\mathbf U)$, $\mathcal G'(\mathbf U)$ are viewed as row vectors, and $\mathbf{F}'\left( \mathbf{U} \right)$, $\mathbf{G}'\left( \mathbf{U} \right)$ are the Jacobian matrixes. 
\end{definition}

If a system of conservation law admits an entropy pair, a weak solution that further satisfies the following inequality in the weak sense is called the \textit{entropy solution}, 
\begin{equation}\label{eq:entropy_condition}
\frac{\partial \mathcal U\left( \mathbf{U} \right)}{\partial t}+\frac{\partial \mathcal F\left( \mathbf{U} \right)}{\partial x}+\frac{\partial \mathcal G\left( \mathbf{U} \right)}{\partial y}\le 0,
\end{equation}
which is known as the \textit{entropy condition}. 
In particular, this inequality becomes an equality for smooth solutions.

 Define the entropy variable $\mathbf V=\mathcal U'(\mathbf U)^T$. If $\mathcal{U}$ is strictly convex, the mapping $\mathbf U\to\mathbf V$ is one-to-one. It is known that the existence of an entropy pair is closely related to the symmetrization of the underlying system of conservation laws. In particular, the system \eqref{eq:MHD} possesses a strictly convex entropy function $\mathcal U$  if and only if $\partial\mathbf U/\partial\mathbf V$ is a symmetric, positive definite matrix and $\partial \mathbf F(\mathbf U(\mathbf V))/\partial\mathbf V, \partial\mathbf G(\mathbf U(\mathbf V))/\partial\mathbf V$ are symmetric matrices \cite{godlewski2013numerical}.
 Therefore, if the transformation $\mathbf U\to\mathbf V$ symmetrizes the system, then there exist twice differentiable scalar functions $\varphi(\mathbf V),\psi_F(\mathbf V),\psi_G(\mathbf V)$ with $\varphi(\mathbf V)$ being strictly convex such that
$$
\mathbf{U}\left( \mathbf{V} \right) ^T=\frac{\partial \varphi}{\partial \mathbf{V}},\quad \mathbf{F}\left( \mathbf{V} \right) ^T=\frac{\partial \psi _F}{\partial \mathbf{V}},\quad \mathbf{G}\left( \mathbf{V} \right) ^T=\frac{\partial \psi _G}{\partial \mathbf{V}},
$$
where $\frac{\partial\varphi}{\partial\mathbf V},\frac{\partial\psi_F}{\partial\mathbf V},\frac{\partial\psi_G}{\partial\mathbf V}$ are viewed as row vectors. It is straightforward to verify that
$$\begin{aligned}
\varphi \left( \mathbf{V} \right) &=\mathbf{V}^T\mathbf{U}(\mathbf V)-\mathcal U(\mathbf{U}(\mathbf V)),\quad\\ \psi _F\left( \mathbf{V} \right) &=\mathbf{V}^T\mathbf{F}(\mathbf{U}(\mathbf V))-\mathcal F(\mathbf{U}(\mathbf V)),\quad\\ \psi _G\left( \mathbf{V} \right) &=\mathbf{V}^T\mathbf{G}(\mathbf{U}(\mathbf V))-\mathcal G(\mathbf{U}(\mathbf V)).\end{aligned}
$$

The entropy pair for ideal MHD equations  is introduced in light of the thermodynamic entropy $s=\ln(p\rho^{-\gamma})$. In particular,  the following equation for $\rho s$ can be derived from the MHD equation:
$$
\frac{\partial \rho s}{\partial t}+\frac{\partial \rho u_xs}{\partial x}+\frac{\partial \rho u_ys}{\partial y}+\left( \gamma -1 \right) \frac{\rho \left( \mathbf{u}\cdot \mathbf{B} \right)}{p}\left( \nabla \cdot \mathbf{B} \right) =0.
$$
We can verify that the following pair satisfies the condition  \eqref{eq:entropy_condition}
\begin{equation}
\label{eq:entropy_pair}
\mathcal U=-\frac{\rho s}{\gamma -1},\quad \mathcal F=-\frac{\rho u_xs}{\gamma -1},\quad \mathcal G=-\frac{\rho u_ys}{\gamma -1},
\end{equation}
and consequently, $\mathbf{V}$ is given by
$$
\mathbf{V}=\left(\frac{\partial \mathcal U}{\partial \mathbf{U}}\right)^T=\left[ \frac{\gamma -s}{\gamma -1}-\beta \left\| \mathbf{u} \right\| ^2,2\beta \mathbf{u},-2\beta ,2\beta\mathbf{B} \right] ^T,\quad \beta=\frac{\rho }{2p}.
$$
However, the change of variables $\mathbf U\to\mathbf V$ fails to symmetrize the ideal MHD equations \eqref{eq:MHD}, and hence $(\mathcal U,\mathcal F,\mathcal G)$ does not constitute an entropy pair.  To overcome the difficulty, Godunov \cite{godunov1972symmetric} introduced a modified form 
\begin{equation}\label{eq:modMHD}
\frac{\partial \mathbf{U}}{\partial t}+\frac{\partial \mathbf{F}\left( \mathbf{U} \right)}{\partial x}+\frac{\partial \mathbf{G}\left( \mathbf{U} \right)}{\partial y}+\phi'\left( \mathbf{V} \right)^T \nabla \cdot \mathbf{B}=0
\end{equation}
with $\phi(\mathbf V) = 2\beta(\mathbf u\cdot\mathbf B)$, which can be symmetrized  with the above transformation $\mathbf U\to\mathbf V$.  Since $\nabla\cdot\mathbf B=0$, the above modification is consistent with the ideal MHD equation \eqref{eq:MHD}. In particular, $\phi(\mathbf V)$ is a homogeneous function of degree one, that is
$ \mathbf V\cdot\phi'(\mathbf V)^T= \phi(\mathbf V), $
with Jacobian matrix
$\phi'(\mathbf V)= [0,\mathbf B,\mathbf u\cdot\mathbf B, \mathbf u]. $

\section{Globally divergence-free entropy stable nodal DG scheme for the conservative form of MHD}
\label{sec3}
\setcounter{equation}{0}

The modified MHD equation \eqref{eq:modMHD} is equivalent to the original conservative  form  \eqref{eq:MHD} only if the magnetic field is  divergence-free. Many existing entropy stable numerical schemes discretize the modified form \eqref{eq:modMHD}, and the additional term is treated as a non-conservative source term.  Thus, without strictly preserving the diveregence-free property, such schemes can ensure the mass conservation but may fail to conserve other quantities, including the momentum, total energy, and magnetic field. Meanwhile, only conservative schemes can guarantee the correct jump conditions and capture a weak solution upon convergence \cite{toth2000b}. Hence, non-conservative schemes are susceptible to producing non-physical results, as demonstrated in the literature \cite{toth2000b, chandrashekar2016entropy, liu2018entropy, wu2018provably}. To partly address the difficulty, we propose to develop a globally divergence-free scheme, enabling us to discrete the conservative form  \eqref{eq:MHD} while still ensuring the entropy condition.

\subsection{Nodal DG method for conservation laws}

The proposed method follows the methodology in \cite{chen2017entropy}, which established a general framework for designing  entropy stable nodal DG schemes for hyperbolic conservation laws. 
The main ingredients include SBP operators derived from the Gauss-Lobatto quadrature, the entropy conservative flux within elements, and the entropy stable flux at element interfaces. This approach has been successfully applied to the ideal MHD in the non-conservative form \cite{liu2018entropy}.

Assume that the 2D spatial domain $\Omega=[a,b]\times[c,d]$ is divided into a $N_x\times N_y$ uniform
rectangular mesh with cells $\mathcal K=\{K_{ij}=[x_{i-1/2},x_{i+1/2}]\times [y_{j-1/2},y_{j+1/2}]\}$, and the step sizes in $x$ and $y$ directions are denoted by $\Delta x$ and $\Delta y$, respectively.

Denote the $k+2$ Gauss–Lobatto quadrature points
$$-1=X_0<X_1<\cdots<X_{k+1}=1,$$
and the associated quadrature weights $\{\omega_l\}_{l=0}^{k+1}$ over the reference element $I=[-1,1]$. 
The difference matrix $D$ is defined as $D_{jl}={L_l}'(X_j)$, where $L_l$ is the $l$-th Lanrange basis polynomial satisfying $L_l(X_j)=\delta_{lj}$, and the mass matrix $M$ and the stiffness matrix $S$ are defined as
 $ M=\mathrm{diag}\{\omega_0,\omega_1,\cdots,\omega_{k+1}\},\  S=MD. $
We recall the properties of these matrices \cite{chen2017entropy}:

\begin{lemma} 
\label{lem:SBP}
Set the boundary matrix
$$B=\mathrm{diag}\{-1,0,0,\cdots,0,0,1\}=:\mathrm{diag}\{\tau_0,\cdots,\tau_{k+1}\},$$
then $S+S^T=B$.
\end{lemma}

\begin{lemma}
For each $0\le j\le k+1$, 
$$\displaystyle
\sum_{l=0}^{k+1}{D_{jl}}=\sum_{l=0}^{k+1}{S_{jl}=0},\quad \sum_{l=0}^{k+1}{S_{lj}}=\tau _j.
$$
\end{lemma}

We further introduce some useful notations. Define the piecewise tensor-product polynomial space  $$
V_{h}=\left\{ w\left( x,y \right): w\left( x,y \right)|_{K_{ij}} \in Q^{k+1}(K_{ij}),\,\forall K_{ij}\in \mathcal{K} \right\},
$$ 
where $Q^{k+1}(K_{ij})$ denotes the set of tensor-product polynomials of degree at most $k+1$ over cell $K_{ij}$. Denote the DG solution $\mathbf U_h\in [V_h]^8$ that approximates  $\mathbf U$. Let
$$
x_i\left( X \right) = x_{i} +\frac{\Delta x}{2}X,\quad  
y_j\left( Y \right) = y_{j} +\frac{\Delta y}{2}Y,$$
and denote the nodal values at the local quadrature points by
$$	\mathbf{U}_{i_1,j_1}^{i,j}=\mathbf{U}_h\left( x_i\left( X_{i_1} \right) ,y_j\left( Y_{j_1} \right) \right). 
$$
Then the solution is represented by
$$\mathbf U_h|_{K_{ij}} = \sum_{i_1,j_1=0}^{k+1} \mathbf U^{i,j}_{i_1,j_1} L_{i_1}\left(\frac{x-x_i}{\Delta x/2}\right) L_{j_1}\left(\frac{y-y_j}{\Delta y/2}\right).$$ 
Denote the nodal values $\mathbf{F}_{i_1,j_1}^{i,j} = \mathbf{F}(\mathbf{U}_{i_1,j_1}^{i,j})$ and $\mathbf{G}_{i_1,j_1}^{i,j} = \mathbf{G}(\mathbf{U}_{i_1,j_1}^{i,j})$,
and define
$$
    \mathbf{F}_{i_1,j_1}^{*,i,j}=\left\{ \begin{array}{ll}
    \mathbf{\hat{F}}\left( \mathbf U^{i-1,j}_{k+1,j_1},\mathbf U^{i,j}_{0,j_1} \right)=:\hat{\mathbf F}^{i-1/2,j}_{j_1} ,& i_1=0,\\	
    \mathbf{0},& 1\leq i_1\leq k,\\
    \mathbf{\hat{F}}\left( \mathbf U^{i,j}_{k+1,j_1},\mathbf U^{i+1,j}_{0,j_1} \right)=:\hat{\mathbf F}^{i+1/2,j}_{j_1} ,& i_1=k+1,\\	
\end{array}\right.
$$
$$
    \mathbf{G}_{i_1,j_1}^{*,i,j}=\left\{ \begin{array}{ll}
	\mathbf{\hat{G}}\left( \mathbf U^{i,j-1}_{i_1,k+1},\mathbf U^{i,j}_{i_1,0} \right)=:\hat{\mathbf G}^{i,j-1/2}_{i_1} ,& j_1=0,\\	
        \mathbf{0},& 1\leq j_1\leq k,\\
	\mathbf{\hat{G}}\left( \mathbf U^{i,j}_{i_1,k+1},\mathbf U^{i,j+1}_{i_1,0} \right)=:\hat{\mathbf G}^{i,j+1/2}_{i_1} ,& j_1=k+1.\\	
	\end{array}\right.
	$$
Here, $\hat{\mathbf F}(\cdot, \cdot),\hat{\mathbf G}(\cdot, \cdot)$ are the numerical fluxes at interfaces. 

Following \cite{chen2017entropy}, the proposed nodal DG scheme for the ideal MHD equations in the conservative form \eqref{eq:MHD} is given by 
\begin{equation}\label{eq:nodalES}\begin{aligned}
	\frac{\mathrm{d}\mathbf{U}_{i_1,j_1}^{i,j}}{\mathrm{d}t}=&-\frac{2}{\Delta x}\sum_{l=0}^{k+1}{2D_{i_1,l}}\mathbf{F}^S\left( \mathbf{U}_{i_1,j_1}^{i,j},\mathbf{U}_{l,j_1}^{i,j} \right) +\frac{2}{\Delta x}\frac{\tau _{i_1}}{\omega _{i_1}}\left( \mathbf{F}_{i_1,j_1}^{i,j}-\mathbf{F}_{i_1,j_1}^{*,i,j} \right) 
	\\
	&-\frac{2}{\Delta y}\sum_{l=0}^{k+1}{2D_{j_1,l}}\mathbf{G}^S\left( \mathbf{U}_{i_1,j_1}^{i,j},\mathbf{U}_{i_1,l}^{i,j} \right) +\frac{2}{\Delta y}\frac{\tau _{j_1}}{\omega _{j_1}}\left( \mathbf{G}_{i_1,j_1}^{i,j}-\mathbf{G}_{i_1,j_1}^{*,i,j} \right), \end{aligned}
	\end{equation}
for $1\le i\le N_x$, $1\le j\le N_y$, $0\le i_1,j_1\le k+1$. Note that the scheme is represented as a matrix vector formulation based on nodal values. 
Here, $\mathbf F^S(\mathbf U_L,\mathbf U_R)$ and $\hat{\mathbf F}(\mathbf U_L,\mathbf U_R)$ are the entropy conservative flux and the entropy stable flux, respectively, which are defined as follows.

\begin{definition}
A consistent, symmetric two-point numerical flux $\mathbf F^S(\mathbf U_L,\mathbf U_R)$ is called entropy conservative, if for the given entropy
function $\mathcal U$,
\begin{equation}\label{eq:ESflux}
\left( \mathbf{V}_R-\mathbf{V}_L \right) ^T \mathbf{F}^S\left( \mathbf{U}_L,\mathbf{U}_R \right) =\left( \psi _R-\psi _L \right) -\left( \phi _R-\phi _L \right) \frac{B_{x,R}+B_{x,L}}{2}.
\end{equation}
In particular, $\mathbf{F}^S\left( \mathbf{U},\mathbf{U} \right) = \mathbf{F} \left( \mathbf{U}\right)$.
\end{definition}

\begin{definition}
A consistent two-point numerical flux $\hat{\mathbf F}(\mathbf U_L,\mathbf U_R)$ is called entropy stable, if for the given entropy
function $\mathcal U$,
\begin{equation}\label{eq:ESflux2}
   \left( \mathbf{V}_R-\mathbf{V}_L \right) ^T\hat{\mathbf{F}}\left( \mathbf{U}_L,\mathbf{U}_R \right) \le \left( \psi _R-\psi _L \right) -\left( \phi _R-\phi _L \right) \frac{B_{x,R}+B_{x,L}}{2}. 
\end{equation}
\end{definition}

In this paper, we employ the entropy conservative flux proposed by Chandrashekar and Klingengberg \cite{chandrashekar2016entropy}, along with the HLL numerical flux with suitable wave speed estimate \cite{bouchut2007multiwave, bouchut2010multiwave} as the entropy stable flux. For more details, see Appendix \ref{app}.

It was proved in \cite{chen2017entropy} that the nodal DG framework achieves the discrete entropy stability for a class of hyperbolic conservation laws, including the compressible Euler equation. Moreover, with the $k+2$ Gauss-Lobatto quadrature points, the truncation error of the scheme is of $(k+1)$-th order at each collocation point.
However, since the scheme \eqref{eq:nodalES} approximates the original conservative form of the MHD equations without preserving the divergence-free property, it cannot achieve the entropy stability. 
To address this, a globally divergence-free mechanism must be incorporated.

In the DG  framework \cite{balsara2021globally}, the globally divergence-free property can be achieved through two conditions: the magnetic field must be locally divergence-free within each cell, and it must be continuous in the normal direction across cell interfaces.
In particular, for the nodal DG formulation, the definition of a globally divergence-free magnetic field is given as follows. 

\begin{definition}
A magnetic field represented by nodal values on Gauss-Lobatto points $\{\mathbf B^{i,j}_{i_1,j_1}\}_{i_1,j_1=0,0}^{k+1,k+1}$ is globally divergence-free, if
\begin{equation}\label{eq:divfree1}
    \sum_{l=0}^{k+1} \left( \frac{2}{\Delta x}D_{i_1,l}B_{x,l,j_1}^{i,j}
    +\frac{2}{\Delta y}D_{j_1,l}B_{y,i_1,l}^{i,j} \right) =0,
\end{equation}
\begin{equation}\label{eq:divfree2}	
    B_{x,k+1,j_1}^{i,j}=B_{x,0,j_1}^{i+1,j} =:B^{i+1/2,j}_{x,j_1},\quad
    B_{y,i_1,k+1}^{i,j}=B_{y,i_1,0}^{i,j+1}=:B_{y,i_1}^{i,j+1/2}
		\end{equation}
  for any $i,j,i_1,j_1$.
\end{definition}

Note that condition \eqref{eq:divfree1} indicates that the divergence of the magnetic field is zero at each quadrature point. Meanwhile, condition \eqref{eq:divfree2} implies that $B_{x,h}$ and $B_{y,h}$ are continuous along the $x$-direction and $y$-direction, respectively. With a globally divergence-free magnetic field satisfying \eqref{eq:divfree1}-\eqref{eq:divfree2}, we are able to establish the following theorem for the proposed nodal DG scheme \eqref{eq:nodalES}. The proof is given in Appendix \ref{app3}.

\begin{theorem}\label{thm:ES}

    Assume that the boundary conditions are periodic, compactly supported, or reflective. 
    If $\mathbf U_h$ has a globally divergence-free magnetic field satisfying \eqref{eq:divfree1} and \eqref{eq:divfree2}, then the scheme \eqref{eq:nodalES} is entropy conservative within a single element
    \begin{equation} \label{eq:entropy_conser}
        \begin{aligned}
		&\frac{\mathrm{d}}{\mathrm{d}t}\left( \frac{\Delta x\Delta y}{4}\sum_{i_1,j_1=0}^{k+1}{\omega _{i_1}\omega _{j_1}\mathcal U_{i_1,j_1}^{i,j}} \right) \\&\quad 
        =-\sum_{j_1=0}^{k+1}{\frac{\Delta y}{2}\omega _{j_1}\left( \mathcal{F} _{k+1,j_1}^{*,i,j}-\mathcal{F} _{0,j_1}^{*,i,j} \right)}
        -\sum_{i_1=0}^{k+1}{\frac{\Delta x}{2}\omega _{i_1}\left( \mathcal{G} _{i_1,k+1}^{*,i,j}-\mathcal{G} _{i_1,0}^{*,i,j} \right)},
        \end{aligned}
    \end{equation}
    and entropy stable on the entire computational domain $\Omega$ in the sense of
    \begin{equation}\label{eq:entropy_stable}
		\frac{\mathrm{d}}{\mathrm{d}t}\int_{\Omega}{\mathcal U\left( \mathbf{U}_h \right) \mathrm{d}x\mathrm{d}y}\approx \frac{\mathrm{d}}{\mathrm{d}t}\left( \frac{\Delta x\Delta y}{4}\sum_{i,j=1}^{N_x,N_y}{\sum_{i_1,j_1=0}^{k+1}{\omega _{i_1}\omega _{j_1}\mathcal U_{i_1,j_1}^{i,j}}} \right) \le 0,
    \end{equation}
    where
    \begin{align*}
		\mathcal{F} _{k+1,j_1}^{*,i,j}&= \left( \mathbf{V}_{k+1,j_1}^{i,j} \right) ^T\mathbf{F}_{k+1,j_1}^{*,i,j}  -  \psi _{F,k+1,j_1}^{i,j} +\phi _{k+1,j_1}^{i,j}B_{x,j_1}^{i+1/2,j},
		\\
		\mathcal{F} _{0,j_1}^{*,i,j}&= \left( \mathbf{V}_{0,j_1}^{i,j} \right) ^T\mathbf{F}_{0,j_1}^{*,i,j}  -\psi _{F,0,j_1}^{i,j} +\phi _{0,j_1}^{i,j}B_{x,j_1}^{i-1/2,j},
		\\
		\mathcal{G} _{i_1,k+1}^{*,i,j}&= \left( \mathbf{V}_{i_1,k+1}^{i,j} \right) ^T\mathbf{G}_{i_1,k+1}^{*,i,j}  -  \psi _{G,i_1,k+1}^{i,j} +\phi _{i_1,k+1}^{i,j}B_{y,i_1}^{i,j+1/2},
		\\
		\mathcal{G} _{i_1,0}^{*,i,j}& = \left( \mathbf{V}_{i_1,0}^{i,j} \right) ^T \mathbf{G}_{i_1,0}^{*,i,j} -\psi _{G,i_1,0}^{i,j} +\phi _{i_1,0}^{i,j}B_{y,i_1}^{i,j-1/2}.
    \end{align*}
\end{theorem}

\subsection{Globally divergence-free method for $\mathbf B$}\label{sect:GDF}
To achieve the globally diver-\\gence-free property under the DG framework, in \cite{balsara2021globally} the authors proposed approximating the magnetic field components on the cell interface as univariate polynomials and utilizing them to reconstruct the magnetic field within each cell, ensuring it is divergence-free internally and continuous in the normal direction. 
The reconstructed magnetic field then replaces the one solved directly from the conservation laws \eqref{eq:MHD}.
Moreover, to guarantee the existence of such a reconstruction, the degree of polynomials at the interfaces should be one degree lower than that of the polynomials within the cell. In particular, we seek  the univariate polynomials of degree $k$ 
$$b_x^{i+1/2,j}(y)\in P^{k} \left(I_{i+1/2,j}^y\right), \quad 
b_y^{i,j+1/2}(x)\in P^{k} \left(I_{i,j+1/2}^x\right),$$
which approximate $B_x$ and $B_y$ on cell interface $I_{i+1/2,j}^y$ and $I_{i,j+1/2}^x$, respectively.
Moreover, denote by $E_z=u_yB_x-u_xB_y$, then $B_x$ and $B_y$ satisfy
\begin{equation}\label{eq:B_interface}
    \frac{\partial B_x}{\partial t}+\frac{\partial E_z}{\partial y}=0, \quad 
    \frac{\partial B_y}{\partial t}-\frac{\partial E_z}{\partial x}=0,
\end{equation}
which can be treated as one-dimensional PDEs, and the following DG weak formulations are  used to update $b_x^{i+1/2,j}$ and $b_y^{i,j+1/2}$ on interfaces:
\begin{equation}\label{eq:scheme_bx}
\begin{aligned}
    \int_{I^y_{i+1/2,j}} \frac{\partial b_x^{i+1/2,j}}{\partial t} w \mathrm{d}y 
    =& \int_{I^y_{i+1/2,j}} \hat{E}_z^{i+1/2,j} \frac{\partial w}{\partial y} \mathrm{d}y 
    - \tilde{E}_z^{i+1/2,j+1/2} w(y^-_{j+1/2})\\
    & 
    + \tilde{E}_z^{i+1/2,j-1/2} w(y^+_{j-1/2}) , 
    \qquad \forall w\in P^k\left( I^y_{i+1/2,j} \right),
\end{aligned}
\end{equation}
\begin{equation}\label{eq:scheme_by}
\begin{aligned}
    \int_{I^x_{i,j+1/2}} \frac{\partial b_y^{i,j+1/2}}{\partial t} w \mathrm{d}x =& -\int_{I^x_{i,j+1/2}} \hat{E}_z^{i,j+1/2} \frac{\partial w}{\partial x} \mathrm{d}x 
    + \tilde{E}_z^{i+1/2,j+1/2} w(x^-_{i+1/2})\\
    & 
    - \tilde{E}_z^{i-1/2,j+1/2} w(x^+_{i-1/2}),
    \qquad \forall w\in P^k\left( I^x_{i,j+1/2} \right).
\end{aligned}
\end{equation}
Note that $\hat E_{z}^{i+1/2,j}$ is the 7th component of $-\hat{\mathbf F}^{i+1/2,j}$, and $\hat E^{i,j+1/2}_{z}$ is the 6th component of $\hat{\mathbf G}^{i,j+1/2}$, respectively. $w^{\pm}$ are the one-side limiting values of $w$ at the interface.
$\tilde E_{z}^{i+1/2,j+1/2}$ is the two-dimensional Riemann solver defined at the vertex $(x_{i+1/2},y_{j+1/2})$, which is single-valued and depends on the limiting states $\mathbf{U}_{i+1/2,j+1/2}^{\pm,\pm}$ from four surrounding cells
 $$
\tilde{E}_z^{i+1/2, j+1/2}=\tilde{E}_z\left( 
\mathbf{U}_{i+1/2,j+1/2}^{+,+},\mathbf{U}_{i+1/2,j+1/2}^{-,+},\mathbf{U}_{i+1/2,j+1/2}^{+,-},\mathbf{U}_{i+1/2,j+1/2}^{-,-} 
\right) .
$$
Here, we employ the two-dimensional HLL flux  $\tilde{E}_z$, see \cite{chandrashekar2020constraint}. Details can be found in Appendix \ref{app2}.

Then, we are able to reconstruct a globally divergence-free magnetic field $\mathbf B_h=\left(B_{x,h}, B_{y,h}\right)\in [V_h]^2$  from $b_x$ and $b_y$ defined along the interfaces, satisfying 
\begin{align}
&1.\quad 
B_{x,h}\left( x_{i\pm 1/2}^\mp,y \right) =b_x^{i\pm 1/2,j}\left( y \right)\quad \mathrm{on}\ \ I_{i\pm1/2,j}^y,\label{eq:cons1}
\\&2.\quad  B_{y,h}\left( x,y_{j\pm 1/2}^\mp \right) =b_y^{i,j\pm 1/2}\left( x \right) \quad \mathrm{on}\ \ I_{i,j\pm 1/2}^x,\label{eq:cons2}
\\&3.\quad  \nabla\cdot\mathbf B_h=0\quad \mathrm{in}\ \ K_{ij}, \label{eq:cons3}&
\end{align}
for any $i,\,j$. Importantly, to ensure the existence of such reconstruction, the 1D polynomials $b_x$ and $b_y$ must satisfy a specific property.  Integrating the divergence-free
condition \eqref{eq:divfree0} on cell $K_{ij}$ and applying Gauss theorem, we get
\begin{align*}
0=& \int_{K_{ij}}{\nabla \cdot \mathbf{B}\,\mathrm{d}x\mathrm{d}y}
=\int_{\partial K_{ij}}{\mathbf{B}\cdot \mathbf{n}\,\mathrm{d}s} \\
=& \int_{I_{i+1/2,j}^y}{B_x\,\mathrm{d}y} -\int_{I_{i-1/2,j}^y}{B_x\,\mathrm{d}y} +\int_{I_{i,j+1/2}^x}{B_y\,\mathrm{d}x} -\int_{I_{i,j-1/2}^x}{B_y\,\mathrm{d}x}. 
\end{align*}
This leads to a necessary condition for $b_x$ and $b_y$ \cite{balsara2021globally}
\begin{equation}\label{eq:cellcons}
\Delta y \left( \bar{b}_x^{i+1/2, j} -\bar{b}_x^{i-1/2,j} \right) 
+ \Delta x \left( \bar{b}_y^{i,j+1/2} - \bar{b}_y^{i,j-1/2} \right) =0,
\end{equation}
where $\bar{(\cdot)}$ denotes the cell average. 
We call \eqref{eq:cellcons}  the \textit{cell-average constraint}. In \cite{balsara2021globally}, it is shown that if the cell-average constraint \eqref{eq:cellcons} holds initially, then the numerical solution of the semi-discrete scheme \eqref{eq:scheme_bx}-\eqref{eq:scheme_by} preserves \eqref{eq:cellcons} at any times.

The detail of the reconstruction procedure will be discussed later.

\begin{remark}
  Due to the reconstruction of $\mathbf B_h$, Theorem \ref{thm:ES} no longer holds. In order to restore the entropy stability, it is important to note that the entropy \eqref{eq:entropy_pair} depends only on density and pressure. Hence, a correction step for the total energy $\mathcal{E}_h$ is needed to ensure that the pressure remains unchanged at each quadrature point. The detail will be provided for fully discrete formulation in Section \ref{sec:framework}. 
\end{remark}

\subsection{Least-square magnetic field reconstruction}\label{sect:rec}

\begin{figure}[htb!]
    \centering
    \includegraphics[width=0.35\linewidth]{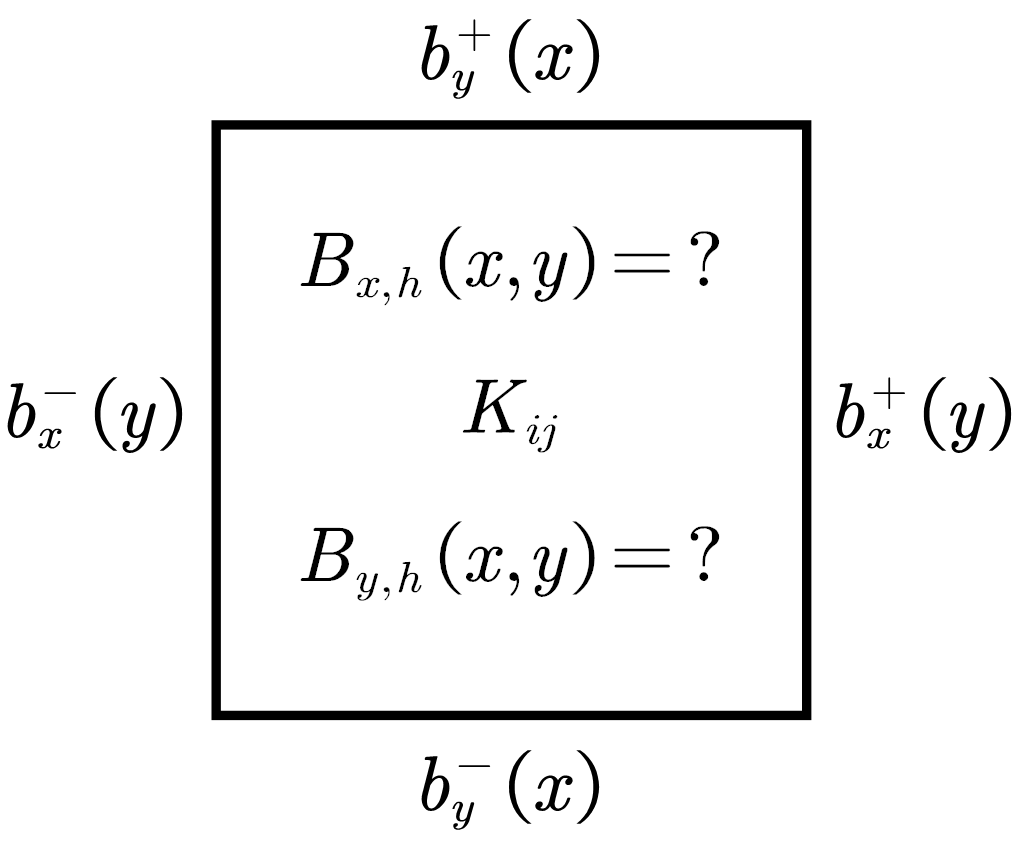}
    \caption{The notations of the reconstruction problem.}
    \label{rec}
\end{figure}

Now, we discuss the globally divergence-free reconstruction. For simplicity of notation, we will omit the subscript $i,\,j$, and denote the magnetic fields along the interfaces of cell $K$ as $b_x^\pm, b_y^\pm$, and the magnetic field inside cell $K$ as $\mathbf B_h=(B_{x,h}, B_{y,h})$, which are illustrated in Fig \ref{rec}. 
In \cite{fu2018globally, chandrashekar2020constraint, balsara2021globally}, $\mathbf B_h\in [P^{k+1}(K)]^2$ is reconstructed from the  polynomials in $P^{k}$ at interfaces.
Since we consider $\mathbf B_h\in [Q^{k+1}(K)]^2$ in the nodal DG framework,  a new technique based on the least-squares (LS) method is developed to directly compute the nodal values of $\mathbf B_h$.

Note that the degrees of approximate polynomials for $b_x$ and $b_y$ differ from that for $\mathbf B_h$. For convenience, we use Legendre basis functions $\{\phi_{l}\}_{l=0}^k$ rather than the Lagrange basis at the interface to represent $b_x^{\pm}$ and $b_y^{\pm}$
$$
b_{x}^{\pm}(y)=\sum_{l=0}^k{b_{x}^{\left( l \right) ,\pm}\phi _{l}\left( Y_j(y) \right)},\quad 
b_y^\pm(x) = \sum\limits_{l=0}^kb_y^{(l),\pm}\phi_l(X_i(x)),
$$
with $\phi_l$ is the Legendre polynomial of degree $l$ defined on the reference interval $[-1,1]$.
We organize the coefficients as vectors
$$ \mathbf b_x^\pm =(b_x^{(0),\pm},\cdots,b_x^{(k),\pm})^T,\quad 
\mathbf b_y^\pm =(b_y^{(0),\pm},\cdots,b_y^{(k),\pm})^T. $$
Meanwhile, we also organzine the nodal values of $B_{x,h}$ and $B_{y,h}$ as vectors
$$\mathbf B_x =\left[B_{x,0,0}, B_{x,1,0}, \cdots, B_{x,k+1,k+1} \right]^T,\quad  
\mathbf B_y=\left[B_{y,0,0}, B_{y,1,0}, \cdots, B_{y,k+1,k+1}\right]^T,$$
and denote $\mathbf B^{rec}=\left[ \mathbf B_x^T,\mathbf B_y^T \right]^T$. 
Then, the reconstruction conditions \eqref{eq:cons1}-\eqref{eq:cons3} are equivalent to the following linear system:
\begin{equation}\label{eq:sys}
		 \left[ \begin{matrix}
		 	\left( \Delta y/\Delta x \right) I_{k+2}\otimes D&		D\otimes I_{k+2}\\
		 	I_{k+2}\otimes r_0&		O\\
		 	I_{k+2}\otimes l_0&		O\\
		 	O&		r_0\otimes I_{k+2}\\
		 	O&		l_0\otimes I_{k+2}\\
		 \end{matrix} \right] \mathbf B^{rec} =\left[ \begin{array}{l}
		 	\mathbf{0}\\
		 	V\mathbf b_{x}^{+}\\
		 	V\mathbf b_{x}^{-}\\
		 	V\mathbf b_{y}^{+}\\
		 	V\mathbf b_{y}^{-}\\
		 \end{array} \right], 
\end{equation}
where $l_0=[1,0,\cdots,0]\in\mathbb R^{k+2},\ r_0=[0,\cdots,0,1]\in\mathbb R^{k+2}$, and the Vandermonde matrix $V\in\mathbb R^{(k+2)\times (k+1)}$ is defined by $V_{ij}=\phi_j(X_i)$. 
It is obvious that there are $N_d=2(k+2)^2$ unknowns for $\mathbf B^{rec}$. However, there are only $(k+2)^2 - 1 + 4(k+1)$ linearly independent conditions in \eqref{eq:sys}. Furthermore, the cell-average constraint \eqref{eq:cellcons} reduces the number of conditions by one \cite{balsara2021globally}. Thus, the total number of linearly independent conditions is given by
$$ N_c=4(k+1)+(k+2)^2-2. $$
Since $N_d\ge N_c$ for any $k\in\mathbb N$, the reconstruction exists but may lack uniqueness.  Note that if we require $b_x, b_y\in P^{k+1}$, it is possible that $N_c>N_d$, which would result in the reconstruction not existing. 
This explains why the degree of the approximate polynomials for $b_x$ and $b_y$ is chosen to be $k$ rather than $k+1$. 

To construct the reconstruction uniquely, we employ the least-squares method. Denote \eqref{eq:sys} as 
\begin{equation}\label{eq:sys2}
   A\mathbf{B}^{rec}=\mathbf b,
\end{equation}
and the magnetic field inside cell $K$ before reconstruction as $\tilde {\mathbf B}$.  
Then $\mathbf{B}^{rec}$ is obtained by solving the following constrained optimization problem 
\begin{equation}\label{eq:opt}
    \min\limits_{\mathbf B } \left\| \mathbf {B}-\tilde {\mathbf B} \right\| _{M_2}^{2},\quad \mathrm{s.t.}\ \ A\mathbf {B}=\mathbf b,
\end{equation}
where
$$
    \left\| \mathbf{X} \right\| _{M_2}^{2}=\mathbf{X}^T M_2\mathbf{X},\quad 
    M_2=\left[ \begin{matrix}
		M\otimes M&		O\\
		O&		M\otimes M\\
	\end{matrix} \right]. 
$$
Assume the singular value decomposition of $A$ is $A=USV^T$, then \eqref{eq:sys2} can be written as $ SV^T \mathbf{B}^{rec}=U^T \mathbf{b}$.  By removing the zero rows from the above system, we obtain an equivalent system with full row rank
$$
	A_1 \mathbf{B}^{rec}=\mathbf{b}_1.
	$$
Lastly, $\mathbf B^{rec}$ is computed by solving the following system
 \begin{equation}\label{eq:LS}
	\left[ \begin{matrix}
		M_2&		A_{1}^{T}\,\\
		A_1&		O\\
	\end{matrix} \right] \left[ \begin{array}{c}
		\mathbf B^{rec}\\
		\boldsymbol{\lambda}\\
	\end{array} \right] =\left[ \begin{array}{c}
		M_2 \tilde {\mathbf B}\\
		\mathbf b_1\\
	\end{array} \right], 
\end{equation}
where $\boldsymbol{\lambda}$ is the Lagrangian multiplier. Since $A_1$ is a full row-rank matrix, the above system \eqref{eq:LS} has a unique solution. 
A more detailed discussion is given in Appendix \ref{app4}.

\section{Limiting strategy and fully discrete formulation}
\label{sec4}
\setcounter{equation}{0}

\subsection{Limiting strategy}\label{sect:limiter}

As mentioned in \cite{chen2017entropy}, entropy stable schemes offer enhanced robustness. However, for problems involving strong shocks, a nonlinear limiter is still required to suppress spurious oscillations. The main challenge lies in ensuring that entropy does not increase after applying the limiter. Many existing limiters such as the popular TVB limiter with characteristic decomposition cannot guarantee this property \cite{liu2018entropy}. To address the difficulty, we recall the following lemma in \cite{chen2017entropy}.

\begin{lemma}\label{lem:limiter}
Suppose $\omega_j>0$ for $1\le j\le N$ with $\displaystyle\sum\limits_{j=1}^{N}\omega_j=1$. Define the average $\bar{\mathbf U} = \displaystyle\sum\limits_{j=1}^{N}\omega_j\mathbf U_j$ and the rescaled values towards the average
$$\tilde{\mathbf U}_j=\bar{\mathbf U}+\theta(\mathbf U_j-\bar{\mathbf U}),\quad j=1,\ldots,N,\quad 0\le\theta\le 1.$$
Then, for any convex function $\mathcal U$, we have
$$
\sum_{j=1}^{N}{\omega _j\mathcal U\left( \tilde{\mathbf U}_j \right)}\le \sum_{j=1}^{N}{\omega _j\mathcal U\left( \mathbf{U}_j \right)}.
$$
\end{lemma}

The theorem ensures that rescaling 
$$\tilde{\mathbf U}^{i,j}_{i_1,j_1}= \bar{\mathbf U}^{i,j} +\theta_{ij} (\mathbf U_{i_1,j_1}^{i,j}-\bar{\mathbf U}^{i,j}), 
\quad 0\le i_1,j_1\le k+1$$
by the same parameter $0\leq \theta_{ij}\le 1$
will not increase the discrete entropy of a cell. Inspired by the oscillation-free DG method \cite{lu2021oscillation} and its subsequent works \cite{peng2024oedg, wei2024jump}, we employ a similar strategy to determine the parameter $\theta_{ij}$ by measuring the jump intensity   of the numerical solution as well as its derivatives at interfaces. Specifically,  $\theta_{ij}$ is given by
\begin{equation}
\theta _{ij}=\exp \left( -\sigma _{ij}\Delta t \right),\quad \sigma _{ij}=\max_{1\le m\le 8} \sigma _{ij}^{\left( m \right)},
\end{equation}
where $\Delta t$ denotes the time step and 
$$
\sigma _{ij}^{\left( m \right)}=c_{f,ij}\sqrt{\Delta x\Delta y}\sum_{l=0}^{k+1}{\sum_{l_1+l_2=l}l(l+1){\Delta x^{l_1}\Delta y^{l_2}\left\| \left[ \mspace{-2.5mu}\left[  D_x^{l_1}D_y^{l_2} u^{\,(m),i,j} \right]\mspace{-2.5mu}\right] \right\| _{\partial K_{ij}}}},
$$
with $D_x = I_{k+2}\otimes D,\ D_y=D\otimes I_{k+2}$, ${u}^{(m)}$ being the $m$-th component of $\mathbf U$
and $\left\| [\mspace{-2.5mu}[\cdot]\mspace{-2.5mu}] \right\|_{\partial K_{ij}}$ representing the jump across the cell boundary 
$$\begin{aligned}
\left\| [\mspace{-2.5mu}[ u ]\mspace{-2.5mu}] \right\| _{\partial K_{ij}} =\sum_{i_1=0}^{k+1} \omega _{i_1}&\left( \left| u^{i,j}_{0,i_1}-u^{i-1,j}_{k+1,i_1} \right|
+ \left| u^{i+1,j}_{0,i_1}-u^{i,j}_{k+1,i_1} \right|\right.
\\&\left.+\left| u^{i,j}_{i_1,0}-u^{i,j-1}_{i_1,k+1} \right|
+\left| u^{i,j+1}_{i_1,0}-u^{i,j}_{i_1,k+1} \right| \right).
\end{aligned}$$
For the parameter $c_f$, we set 
$$
c_{f,ij}=\frac{c_0}{4k(k+1)}\max_{i_1,j_1}\left\{\frac{1}{H_{i_1,j_1}^{i,j}}\right\},\quad  H=\frac{\mathcal{E} +p^{\star}}{\rho},
$$
where $c_0>0$ is a free-parameter. 

In addition, we need to apply a limiter to the magnetic field on interfaces to ensure the globally divergence-free property. For example, for the magnetic field $b_x$ at interface $I_{i+1/2,j}^y$, we propose limiting its high-order moments using $\theta_{ij}$ in its adjacent cells, while keeping the average unchanged.
Specifically, we set 
$$\tilde{b}_{x}^{\left( l \right)} =\theta_{i+1/2,j}b_{x}^{\left( l \right)}, \quad\text{for} \quad l=1,\ldots,k,$$
with $\theta_{i+1/2,j}=\min\left\{\theta_{ij},\theta_{i+1,j}\right\}$. 
$b_y$ is limited similarly. Then, the interior magnetic field is reconstructed from $\tilde b_x$ and $\tilde b_y$.
Note that the limiting of $b_x$ and $b_y$ will preserve the discrete entropy in each cell through an energy correction, which will be introduced later. 
Furthermore, since the limiter also maintains the cell averages, the cell-average constraint \eqref{eq:cellcons}  remains satisfied. Hence, we can reconstruct a globally divergence-free magnetic field without spurious oscillations.

\subsection{The fully discrete ES-GDF framework} \label{sec:framework}
For time discretization, we employ the 10-stage 4th-order strong stability-preserving Runge–Kutta (SSP-RK4) scheme described in \cite{ketcheson2008highly}. For illustrative purposes, we present the algorithm flowchart for the proposed ES-GDF scheme coupled with the forward Euler method to update the solution from time $t^n$ to $t^{n+1}=t^n+\Delta t$. This approach can be easily extended to the SSP-RK4 time discretization.

\textbf{Step 1}\quad Update the numerical solution in cells based on the nodal DG scheme \eqref{eq:nodalES} and obtain
$\tilde{\mathbf U}^{n+1} = (\tilde\rho, \widetilde{\rho \mathbf u}, \tilde{\mathcal E}, \tilde{\mathbf B})^T$.

\textbf{Step 2}\quad Update the magnetic field $ \tilde {b}_{x}^{n+1}$ and $\tilde{b}_{y}^{n+1}$ on the interfaces based on \eqref{eq:scheme_bx}-\eqref{eq:scheme_by}.

\textbf{Step 3}\quad Apply the limiter to $\tilde b_x^{n+1}, \tilde b_y^{n+1}$ on the interfaces and $\tilde{\mathbf U}^{n+1}$ in cells to obtain $(b_x^{n+1}, b_y^{n+1})$ and $\mathbf U^{*,n+1} = (\rho_h^{n+1}, (\rho \mathbf u)_h^{n+1}, \mathcal E_h^*, \mathbf B_h^*)^T$.

\textbf{Step 4}\quad Reconstruct the internal magnetic field $\mathbf B^{n+1}$ using $b_x^{n+1}, b_y^{n+1}$, and correct the total energy $\mathcal E^{n+1}$. 
Moreover, in order to maintain the discrete entropy, we also adjust the energy to ensure that the pressure remains unchanged, 
\begin{equation}\label{eq:modene} 
\mathcal E_h^{n+1} = \mathcal E_h^* + \frac{1}{2} \left( \left\|\mathbf B_h^{n+1}\right\|^2 - \left\| \mathbf B_h^*\right\|^2 \right), 
\end{equation}
at each Gauss-Lobatto point.  
Finally, we can obtain $$\mathbf U_h^{n+1} = (\rho_h^{n+1}, (\rho \mathbf u)_h^{n+1}, \mathcal E_h^{n+1}, \mathbf B_h^{n+1})^T.$$

\begin{remark}
The key step the ES-GDF DG framework is the correction procedure \eqref{eq:modene} in \textbf{Step 4}. Note that if we simply set $\mathcal E_h^{n+1}=\mathcal E_h^*$, the change in the magnetic field would lead to a corresponding change in pressure, and there would be no guarantee that 
$$\mathcal U^{tot}((\rho_h^{n+1},(\rho\mathbf u)_h^{n+1},\mathcal E_h^*,\mathbf B_h^{n+1})^T)\le \mathcal U^{tot}(\mathbf U^{*,n+1}).$$
Here, $\mathcal U^{tot}(\mathbf U_h)$ denotes the discrete total entropy of $\mathbf U_h$ on $\Omega$. 
Note that this treatment will result in the loss of local energy conservation. Meanwhile, the corresponding conservation error in each cell $K_{ij}$ is associated with changes in the magnetic energy density due to the reconstruction:
$$
\bar{\mathcal{E}}^{n+1,i,j}-\bar{\mathcal{E}}^{*,i,j}
=\sum_{i_1,j_1=0}^{k+1} \frac{1}{8} \, {\omega _{i_1}\omega _{j_1} \left( \left\| \mathbf{B}_{i_1,j_1}^{n+1,i,j} \right\| ^2-\left\| \mathbf{B}_{i_1,j_1}^{*,i,j} \right\| ^2 \right)}. $$
 Thus, as long as the proposed reconstruction procedure maintains high-order accuracy, the conservative error of $\mathcal E$ is also high-order accurate \cite{balsara1999staggered}.  Furthermore, we may add  an additional constraint to the optimization problem \eqref{eq:opt} to ensure 
 the local magnetic energy is conserved, 
 hence preserving the local total energy conservation. However, it is nontrivial to rigorously show that the admissible set is nonempty. We will leave the investigation for future work.

\end{remark}

\section{Numerical experiments}
\label{sec5}
\setcounter{equation}{0}
\setcounter{table}{0}
\setcounter{figure}{0}

In this section, we simulate the ideal MHD equation to verify the accuracy and efficiency of our scheme. 
We use the time step
$$
\Delta t=\frac{\mathrm{CFL}}{\alpha _x/\Delta x+\alpha _y/\Delta y},\quad \mathrm{CFL}=0.45
$$
for all tests. Here, $\alpha_x=\max|\lambda(\partial \mathbf F/\partial \mathbf U)|$, $\alpha_y=\max|\lambda(\partial \mathbf G/\partial \mathbf U)|$, and $\lambda(\cdot)$ is the eigenvalue of the matrix. We use uniform meshes for convenience. 
To compare the performance of the schemes, we denote ``ES-GDF" as our proposed scheme, and ``ES" as the original entropy stable nodal DG scheme designed in \cite{liu2018entropy}. 
Without special declaration, we only present the results of $k=2$. For the example containing strong shock, the limiter introduced in Section \ref{sect:limiter} will be used.\\

\begin{figure}[htbp!]
    \centering
    \subfigure[ES.]{
        \includegraphics[width=0.4\linewidth]{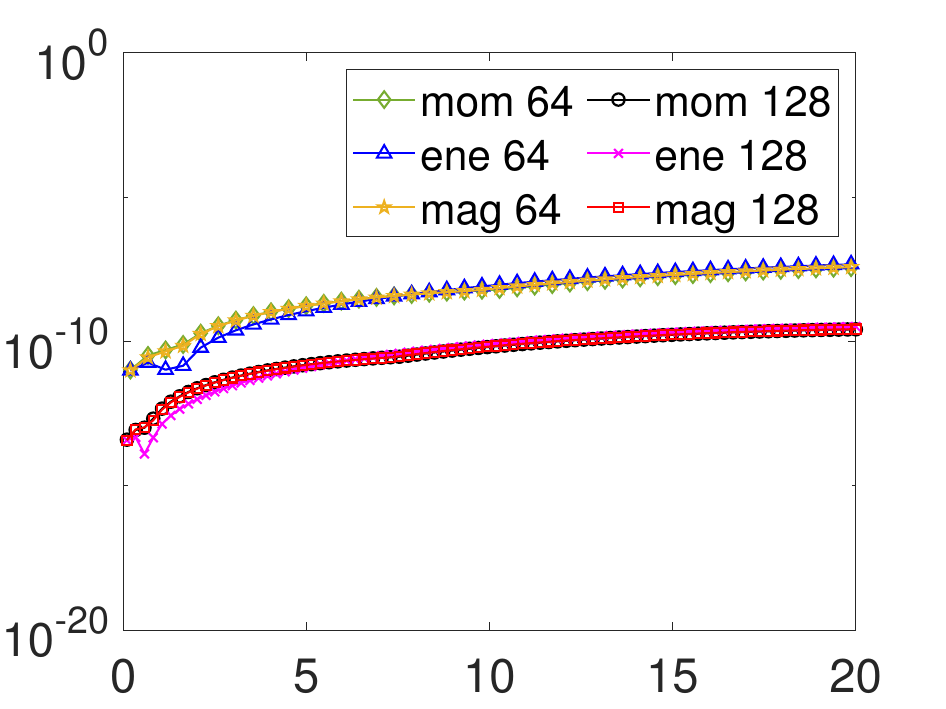}}
    \subfigure[ES-GDF.]{
        \includegraphics[width=0.4\linewidth]{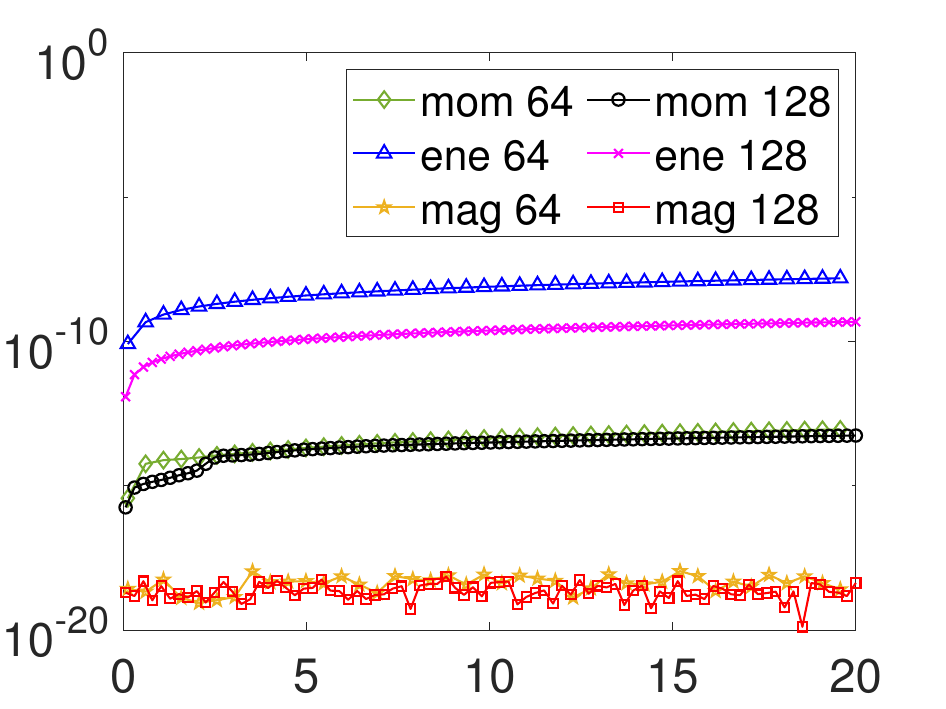}}  
    \caption{Example 5.1:  Smooth MHD vortex. 
    Time evolution of the absolute values of relative deviation in the conservative variables on different meshes for $k=2$.}
    \label{figcons}
\end{figure}

\textit{Example} 5.1 (Smooth MHD vortex).  First, we consider the MHD vortex problem, which is originally introduced by Shu \cite{shu2020essentially} in the hydrodynamical system, and was extended to MHD equations by Balsara \cite{balsara2004second}. Here we use the same setup in \cite{li2005locally}. We use a uniform mesh with $N=N_x=N_y$.
In Table \ref{tab2}, we present the $L^2$ errors and orders of the numerical solution with $k=1,2,3$ at \(T=20\), when the exact solution coincides with the initial data. It can be seen that our schemes also maintain the designed accuracy. 
In Fig \ref{figcons}, we plot the time evolution of the absolute values of relative deviation in the conservative variables for ES-GDF and ES schemes on $N_x\times N_y=64\times 64$ and $N_x\times N_y=128\times 128$ meshes with $k=2$. 
We recall that both two schemes are conservative in density. Hence, we only plot momentum, total energy, and magnetic field here (denoted by ``mom", ``ene", and ``mag"). 
The error of momentum and magnetic field is calculated in $L^1$ vector norm. It is observed that both schemes are not conservative in total energy, and the errors are at the same level. 
However, for momentum and magnetic field, the proposed ES-GDF scheme maintains the conservation errors at machine accuracy, while the ES scheme shows non-conservation. 

 \begin{table}[htb!]
		\centering
        \caption{Example 5.1: Smooth MHD vortex.  The $L^2$ errors and orders at final time $T = 20$.}
        \renewcommand{\arraystretch}{1.25}
		\setlength{\tabcolsep}{2.2mm}{
			\begin{tabular}{|c|c|c|c|c|c|c|c|c|}
\hline $N$ & $\rho$ & order & $\rho u_x$ & order & $B_x$ & order & $\mathcal E$ & order \\ 			
				\hline \multicolumn{9}{|c|}{$k=1$} \\				
				\hline  
                     32 &3.14e-04 & -- &1.51e-03 & --&4.87e-03 &-- &2.47e-03 &-- \\
                \hline  
                     64 &9.62e-05 & 1.71 &2.23e-04 & 2.76 &9.93e-04 &2.29 &4.40e-04 &2.49 \\
                \hline  
                    128 &1.69e-05 & 2.51 &3.69e-05 & 2.60 &1.50e-04 &2.72 &7.39e-05 &2.57 \\
                \hline  
                    256 &3.03e-06 & 2.48 &7.55e-06 & 2.29 &2.54e-05 &2.57 &1.45e-05 &2.35 \\								\hline \multicolumn{9}{|c|}{$k=2$} \\				
				\hline  
				$     32  $ &3.02e-05 & --&5.82e-05 &-- &1.36e-04 & --&9.80e-05 &-- \\
				\hline  
				$     64  $ &3.49e-06 & 3.11 &3.83e-06 &3.92 &1.39e-05 &3.29 &9.05e-06 &3.44 \\
				\hline  
				$    128  $ &3.28e-07 & 3.41 &3.47e-07 &3.46 &2.13e-06 &2.71 &9.17e-07 &3.30 \\
				\hline  
				$    256 $ &3.85e-08 & 3.09 &4.24e-08 &3.03 &2.84e-07 &2.91 &1.11e-07 &3.05 \\					\hline \multicolumn{9}{|c|}{$k=3$} \\				
				\hline  
                     32 &1.96e-06 & -- &2.59e-06 & -- &7.29e-06 & --&5.58e-06 &-- \\
                \hline  
                     64 &1.49e-07 & 3.72 &1.55e-07 & 4.06 &4.29e-07 &4.09 &3.68e-07 &3.92 \\
                \hline  
                    128 &5.49e-09 & 4.76 &6.22e-09 & 4.64 &2.65e-08 &4.02 &1.44e-08 &4.68 \\
                \hline  
                    256 &2.48e-10 & 4.47 &3.35e-10 & 4.21 &1.65e-09 &4.00 &7.14e-10 &4.33 \\			
				\hline  
		\end{tabular}} 
		\label{tab2}
	\end{table}

\begin{figure}[htbp!]
		\centering
		\subfigure[$\rho$.]{
			\includegraphics[width=0.31\linewidth]{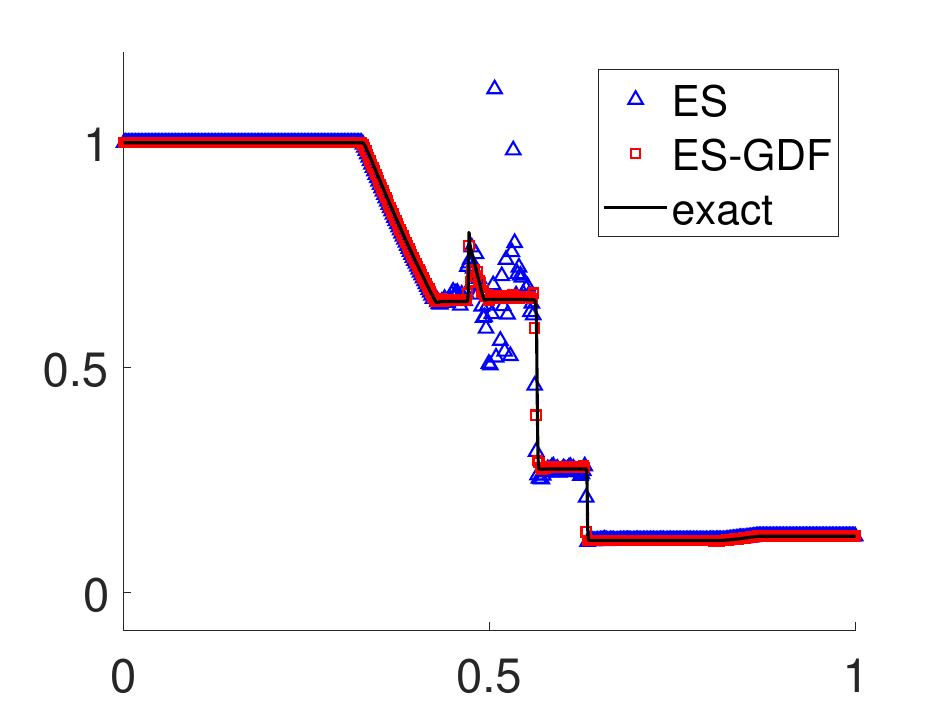}}
		\subfigure[$p$.]{
			\includegraphics[width=0.31\linewidth]{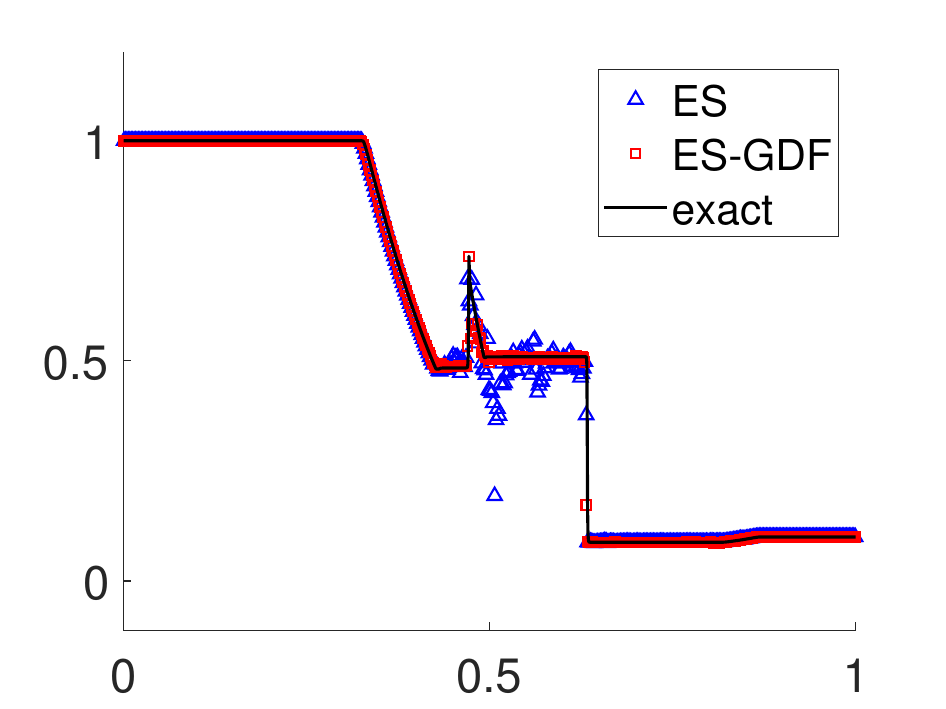}}
		\subfigure[$u_{||}=(2u_x+u_y)/\sqrt{5}$.]{
			\includegraphics[width=0.31\linewidth]{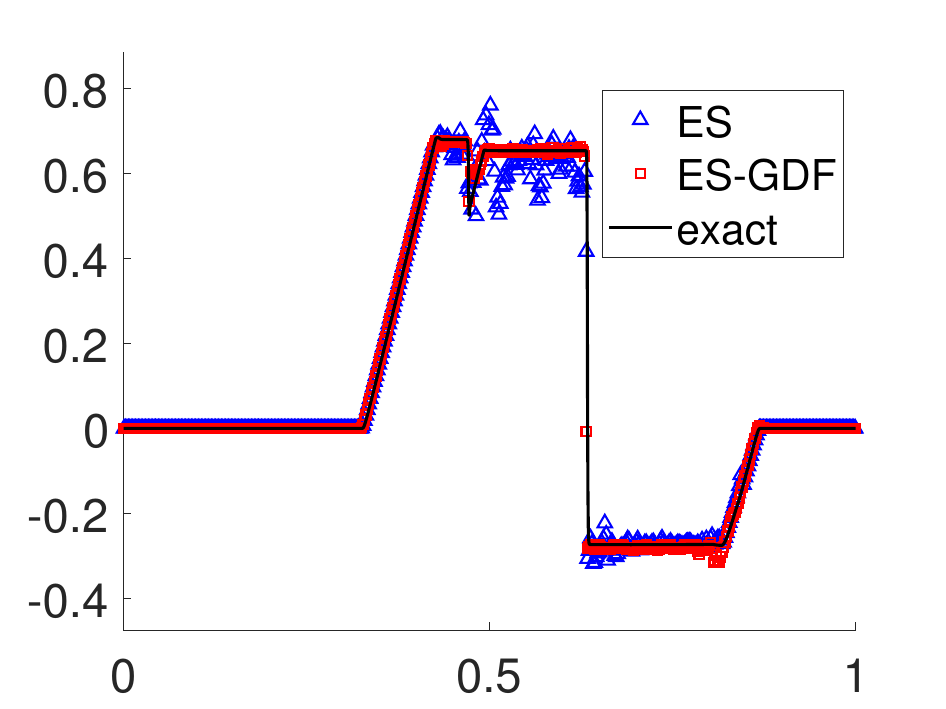}}
		\subfigure[$u_{\bot}=(-u_x+2u_y)/\sqrt{5}$.]{
			\includegraphics[width=0.31\linewidth]{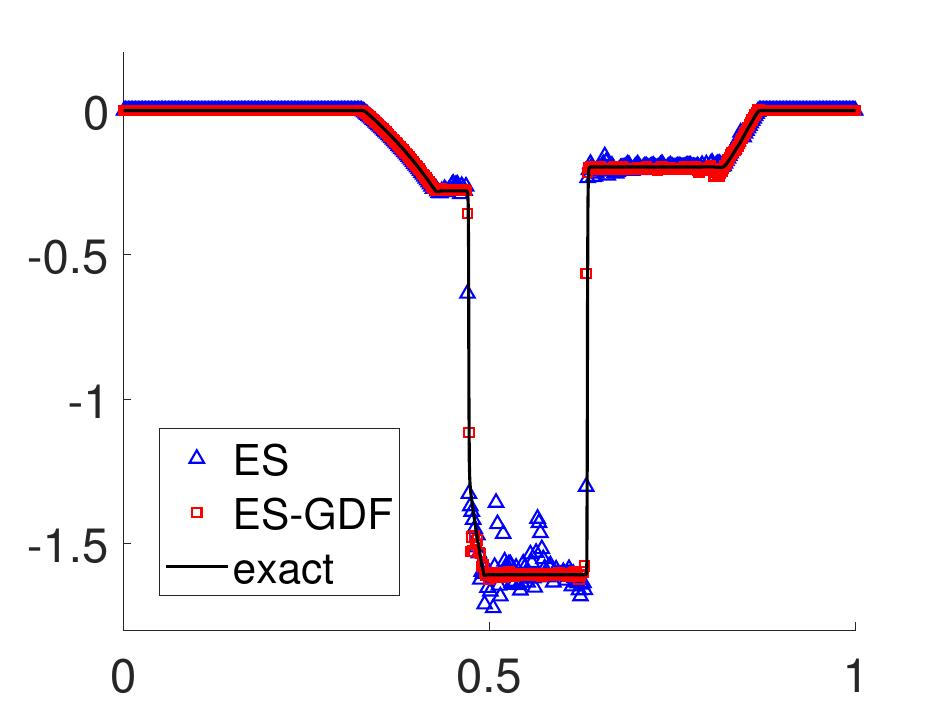}}
		\subfigure[$B_{||}=(2B_x+B_y)/\sqrt{5}$.]{
			\includegraphics[width=0.31\linewidth]{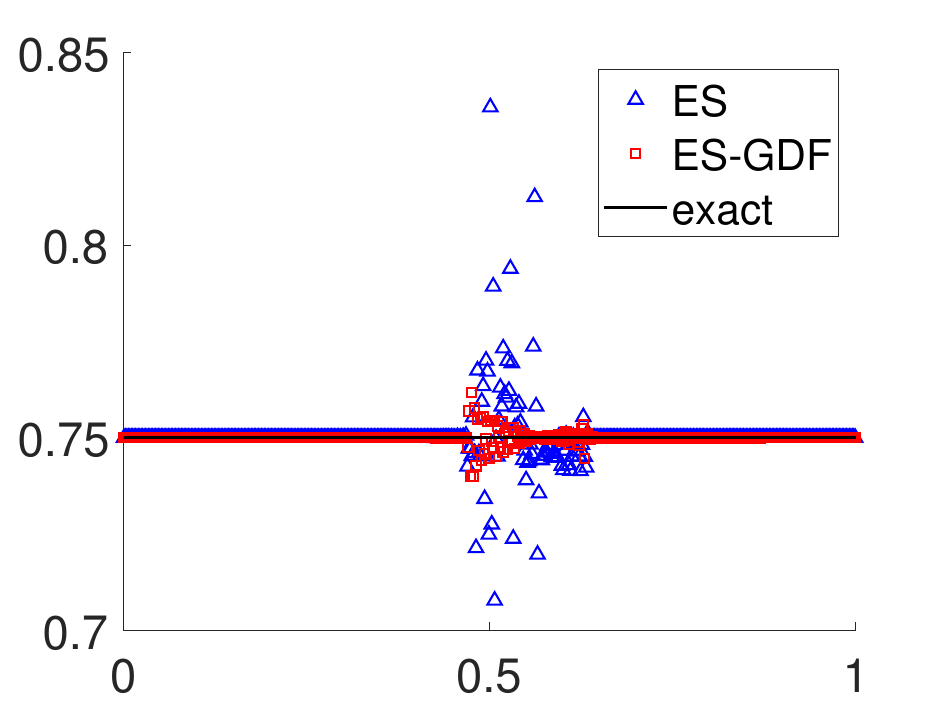}}
		\subfigure[$B_{\bot}=(-B_x+2B_y)/\sqrt{5}$.]{
			\includegraphics[width=0.31\linewidth]{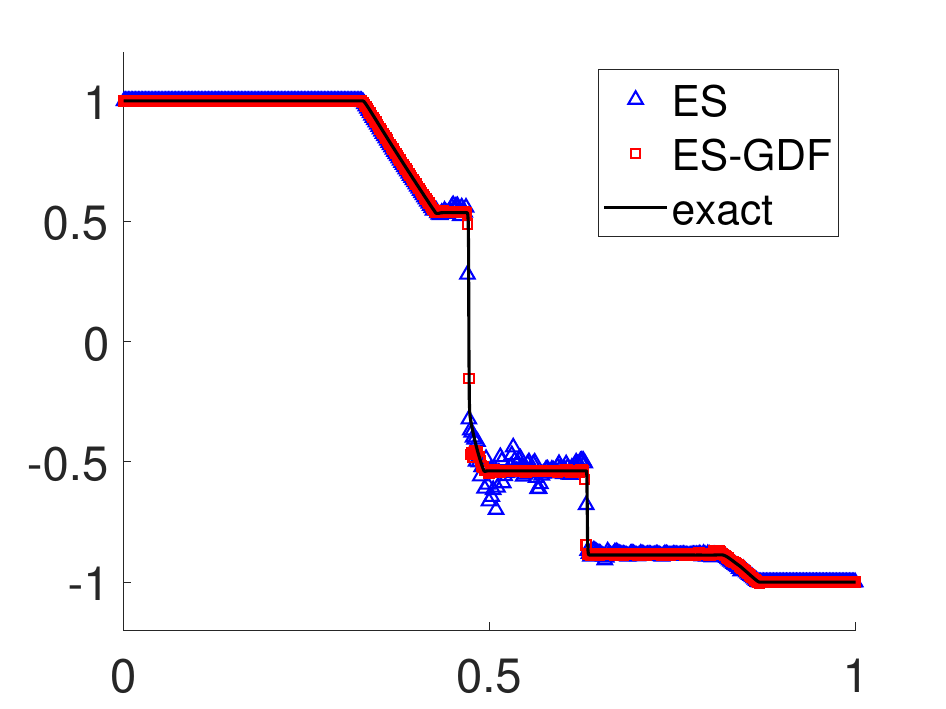}}		
		\caption{Example 5.2: The rotated shock tube. The numerical solution at \(T = 0.2/\sqrt{5}\) on \(N_x\times N_y = 512\times 2\) meshes. No limiter is applied.}
		\label{figBrio}
	\end{figure}

\begin{figure}[htbp!]
		\centering
		\subfigure[$\rho$.]{
			\includegraphics[width=0.31\linewidth]{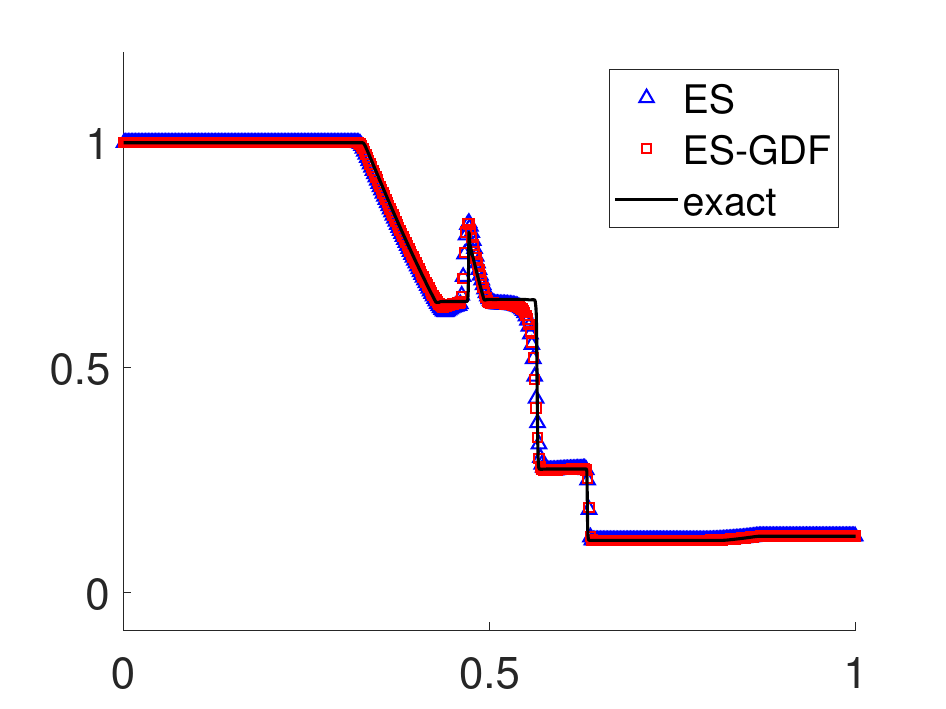}}
		\subfigure[$p$.]{
			\includegraphics[width=0.31\linewidth]{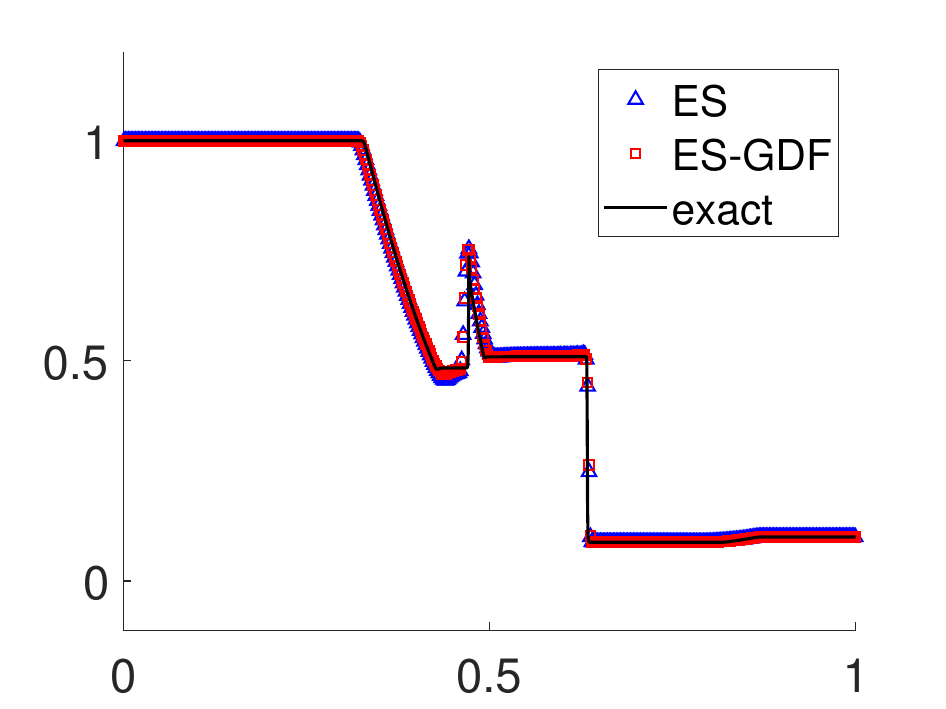}}
		\subfigure[$u_{||}=(2u_x+u_y)/\sqrt{5}$.]{
			\includegraphics[width=0.31\linewidth]{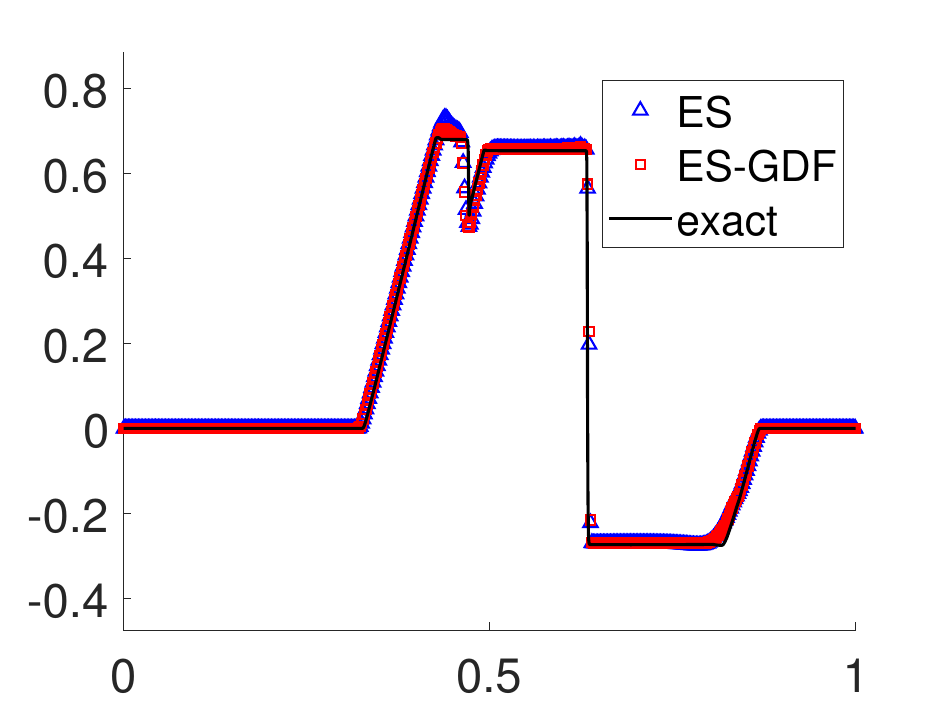}}
		\subfigure[$u_{\bot}=(-u_x+2u_y)/\sqrt{5}$.]{
			\includegraphics[width=0.31\linewidth]{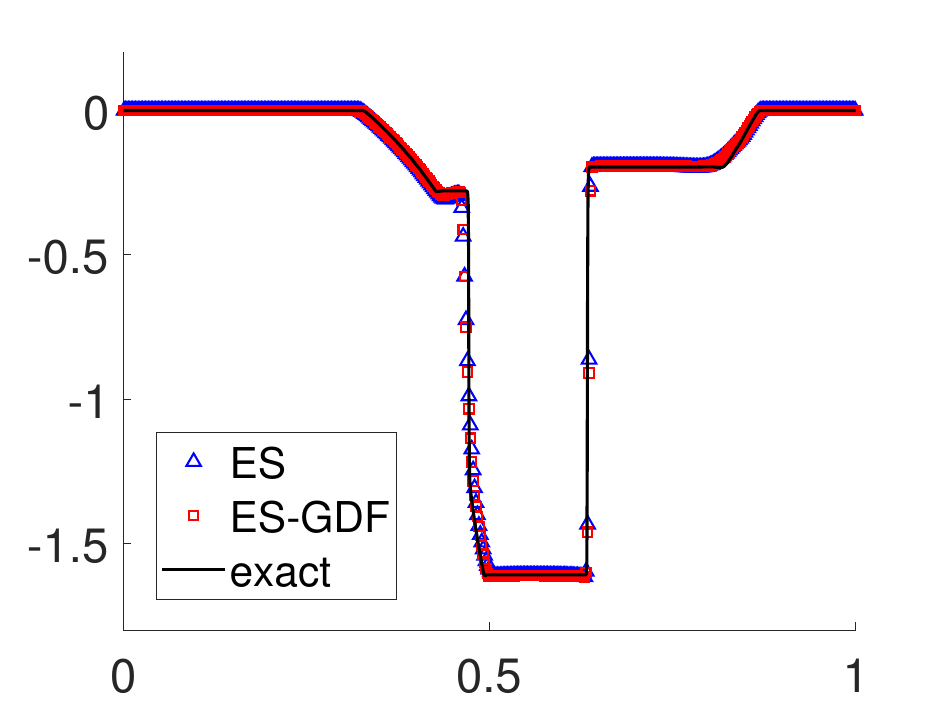}}
		\subfigure[$B_{||}=(2B_x+B_y)/\sqrt{5}$.]{
			\includegraphics[width=0.31\linewidth]{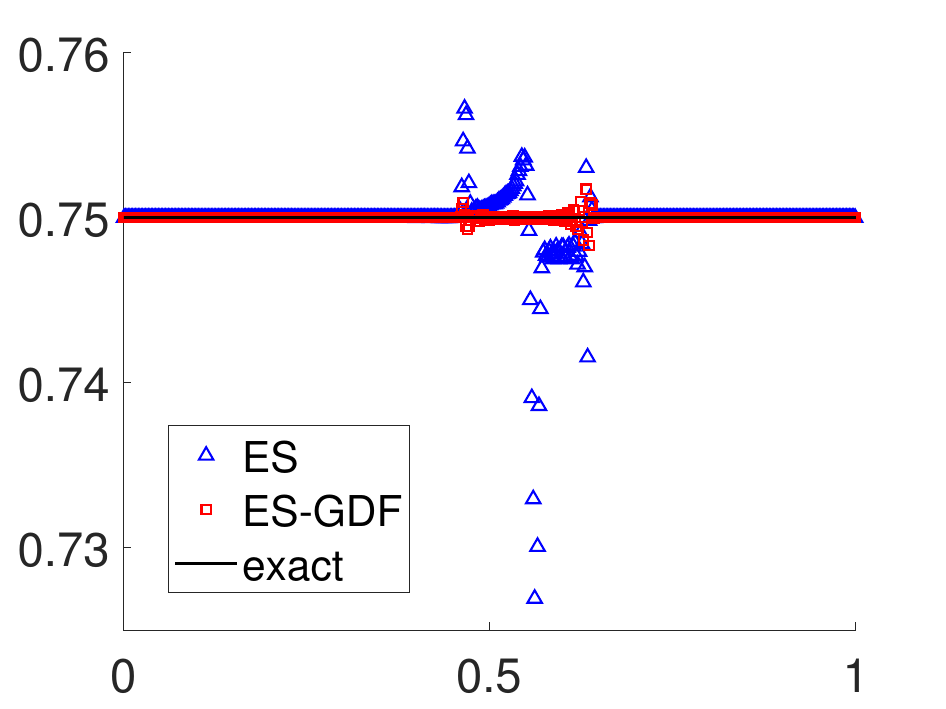}}
		\subfigure[$B_{\bot}=(-B_x+2B_y)/\sqrt{5}$.]{
			\includegraphics[width=0.31\linewidth]{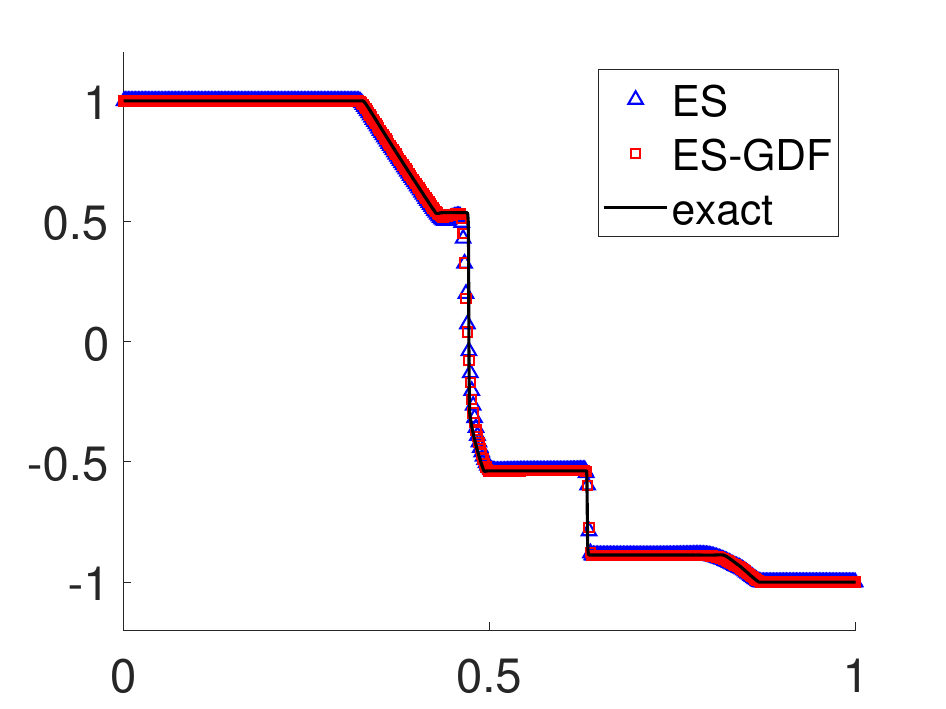}}	
		\caption{Example 5.2: The rotated shock tube. The numerical solution at \(T = 0.2/\sqrt{5}\) on \(N_x\times N_y = 512\times 2\) meshes with the limiter.}
		\label{figBriolim}
	\end{figure}

\textit{Example} 5.2 (Rotated shock tube).	We consider a 2D Riemann problem obtained by rotating the one-dimensional Brio-Wu shock tube \cite{brio1988upwind} with an angle \(\alpha =\arctan(0.5) \). We want to remark that in this problem, the constancy of the parallel component $B_{||}:=B_x\cos\alpha + B_y\sin\alpha$ is difficult to obtain in numerical schemes which do not satisfy the discrete divergence-free property. 
    The computational domain is \(\Omega=[0,1]\times[(2/N_x)\cot\alpha]\). The Dirichlet boundary conditions are applied to the left and right boundaries according to the initial condition. The top and bottom boundaries are imposed with the shifted periodic type according to the translational symmetry, as detailed in \cite{toth2000b}. 
    The ``exact" solution is computed by a classical third-order one-dimensional DG scheme using TVB limiter with $M=1$ on 10,000 cells for the non-rotated version, which can be treated as a one-dimensional problem. 
    The initial data is given by
	$$
	\left( \rho,u_x,u_y,u_z,B_x,B_y,B_z,p \right) =\begin{cases}
		\left( 1,0,0,0,\frac{0.5}{\sqrt{5}},\frac{2.75}{\sqrt{5}},0,1 \right) , \quad\quad\quad\,\, 2x+y<1,\\
		\left( 0.125,0,0,0,\frac{2.5}{\sqrt{5}},\frac{-1.25}{\sqrt{5}},0,0.1 \right) , \,\, 2x+y\ge 1.\\
	\end{cases}
	$$
	Notably, without the entropy stable treatment, both the standard nodal DG scheme and GDF nodal DG scheme will break down in the first several steps. In Fig \ref{figBrio}, we present the result at $T=0.1\cos\alpha=0.2/\sqrt 5$ on \(N_x\times N_y=512\times 2\) meshes without limiter. 
    It is clear that the ES-GDF DG scheme produces a more stable result, while the ES DG scheme appears unstable in the range \(0.45<x<0.62\). This may be caused by the non-constancy of $B_{||}$. 
    The results with the limiter is shown in Fig \ref{figBriolim}, we can see the ES-GDF DG scheme gives a more compatible simulation, especially for $B_{||}$. \\

\begin{figure}[htbp!]
		\centering
		\subfigure[ES.]{
        \includegraphics[width=0.35\linewidth]{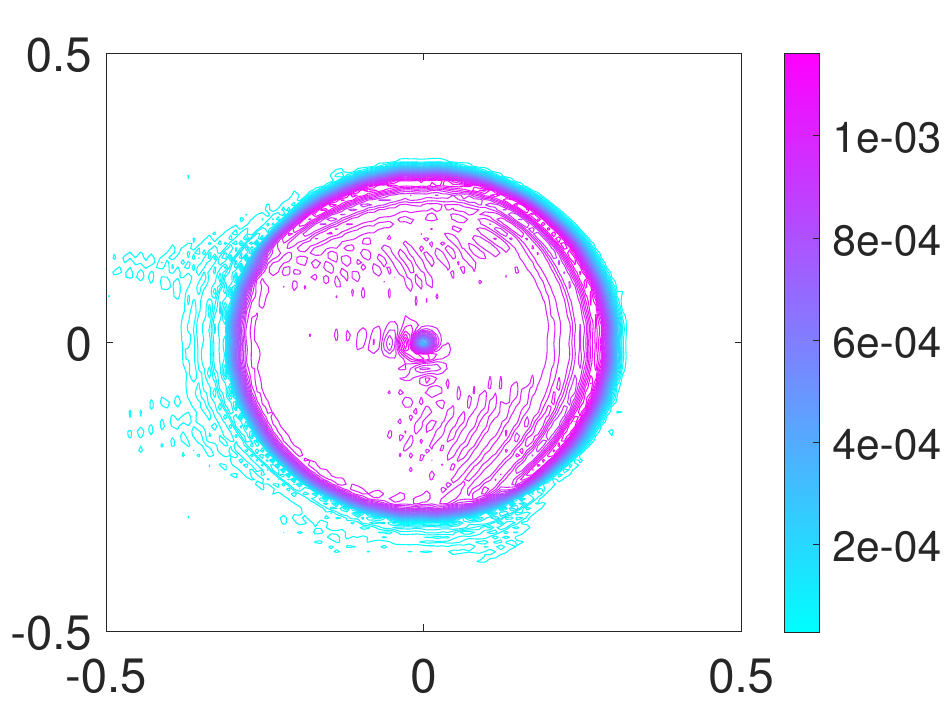}}
		\subfigure[ES-GDF.]{
        \includegraphics[width=0.35\linewidth]{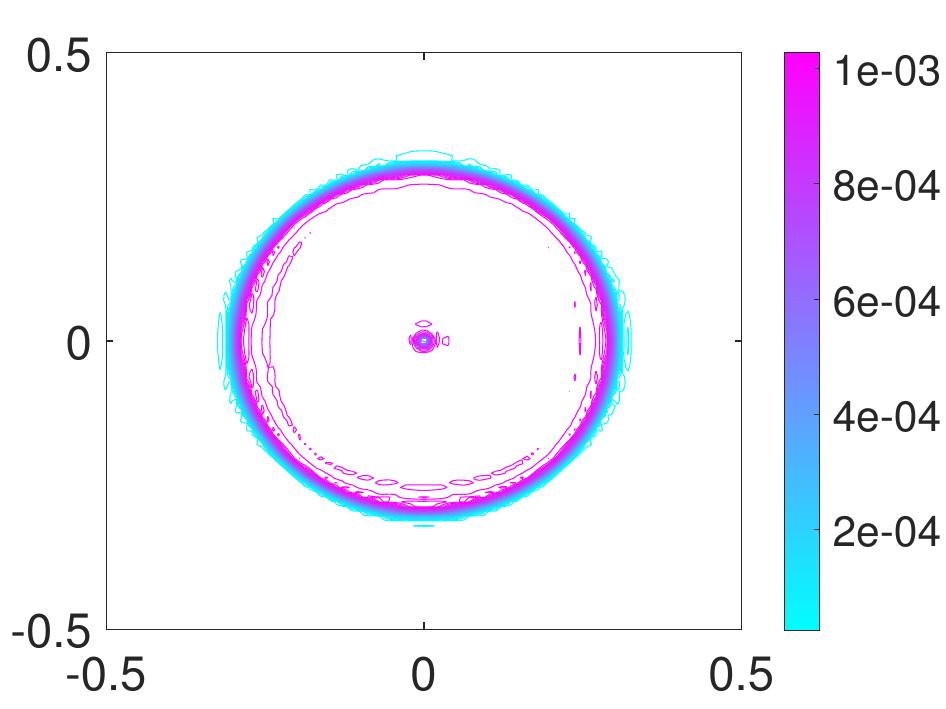}}		
		\caption{Example 5.3: Magnetic field loop. The numerical solution of $\sqrt{B_x^2+B_y^2}$ at $T = 2$ on $N_x\times N_y = 240\times 120$ meshes with 40 contour lines. No limiter is applied.}
		\label{figLoop1}
	\end{figure}

\begin{figure}[htbp!]
		\centering
		\subfigure[Total entropy.]{
        \includegraphics[width=0.35\linewidth]{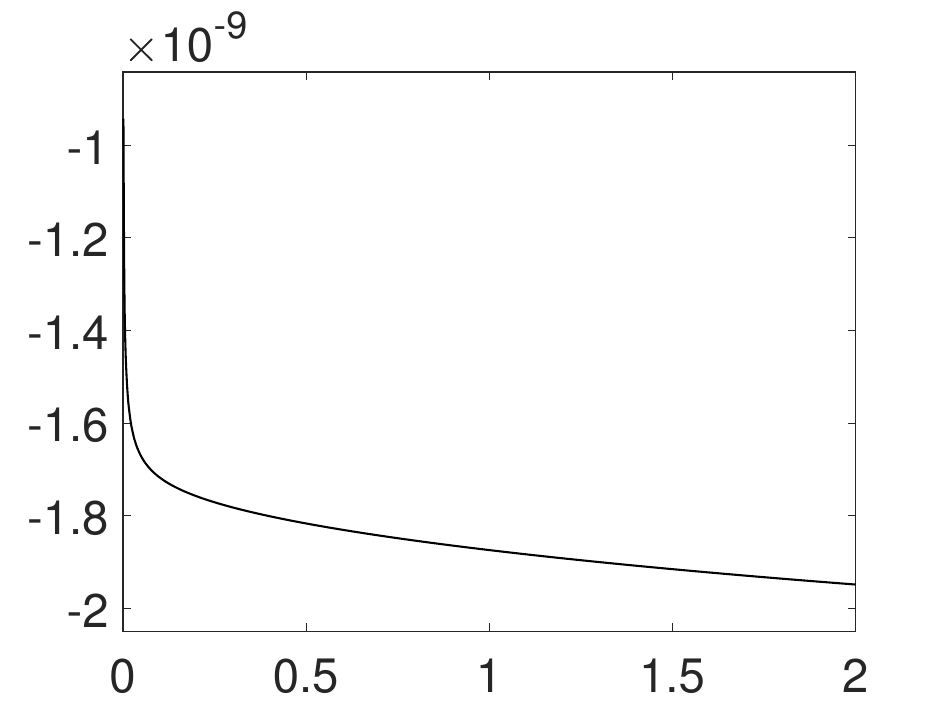}}
		\subfigure[Divergence.]{
        \includegraphics[width=0.35\linewidth]{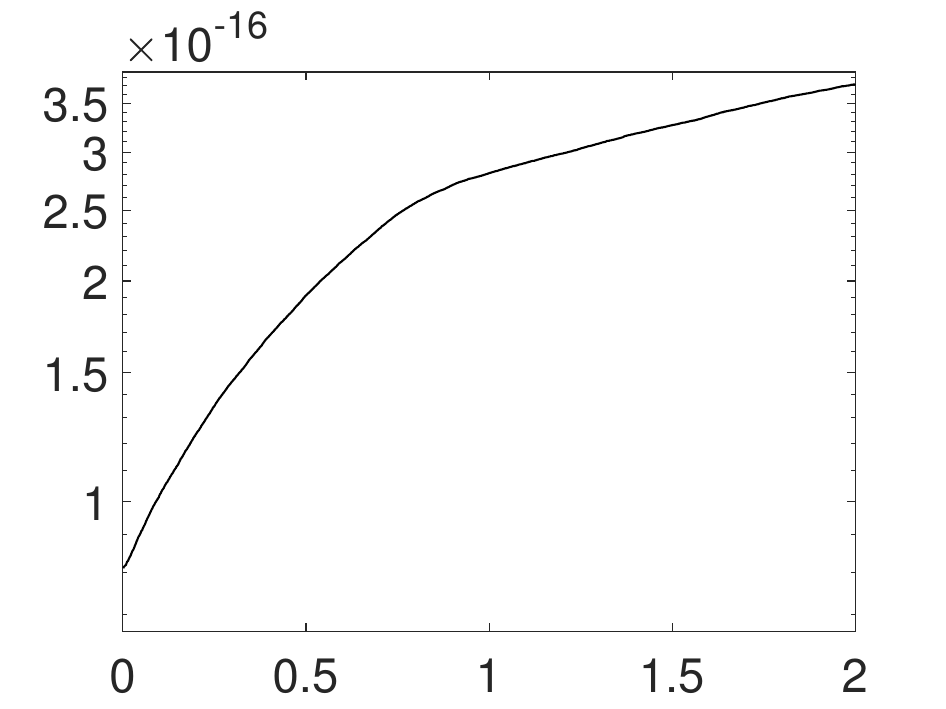}}		
		\caption{Example 5.3: Magnetic field loop. The time evolution of total entropy and divergence norm on $N_x\times N_y=240\times 120$ meshes. No limiter is applied.}
		\label{figLoop2}
	\end{figure}

\textit{Example} 5.3 (Magnetic field loop).	We consider the magnetic field loop problem, which is first introduced in \cite{gardiner2005unsplit}, and our setup is the same as in \cite{li2011central}. This test involves the advection of magnetic field loop over a periodic domain. 
    The magnetic field will transport over the domain and return to its initial position at positive integer times. This problem can be computed without any limiter with ES treatment. 
    In Fig \ref{figLoop1}, we present the result of $\sqrt{B_x^2+B_y^2}$ at $T=2$ on $N_x\times N_y=240\times 120$ meshes. Here, we also compare it with the ES DG scheme \cite{liu2018entropy}, showing that the ES-GDF scheme better captures the features of the solution. In Fig \ref{figLoop2}, we present the evolution of total entropy and divergence over time, where the divergence is defined as \cite{cockburn2004locally}
	$$ \left\| \mathrm{div}\mathbf B \right\|=\sum\limits_{K\in\mathcal K}\left(\int_{K}\left|\nabla\cdot\mathbf B\right|\mathrm dx\mathrm dy+\int_{\partial K}|[\mspace{-2.5mu}[\mathbf B\cdot\mathbf n]\mspace{-2.5mu}]|\mathrm ds\right). $$
	It can be seen that the entropy is non-increasing over time, and the total divergence of the magnetic field closes to machine error.\\

\textit{Example} 5.4 (Kelvin-Helmholtz instability). 
We consider the Kelvin-Helmholtz instability problem and follow the setup in \cite{mignone2010high}. This instability has been widely studied in the literature, and it clearly shows the advantages of maintaining a divergence-free magnetic field. 
The ES-GDF DG scheme can simulate this problem without any limiters, while the non-ES DG schemes and the ES DG scheme without GDF treatment will blow up if the limiter is not applied.

\begin{figure}[htbp!]
    \centering
    \includegraphics[width=0.8\linewidth]{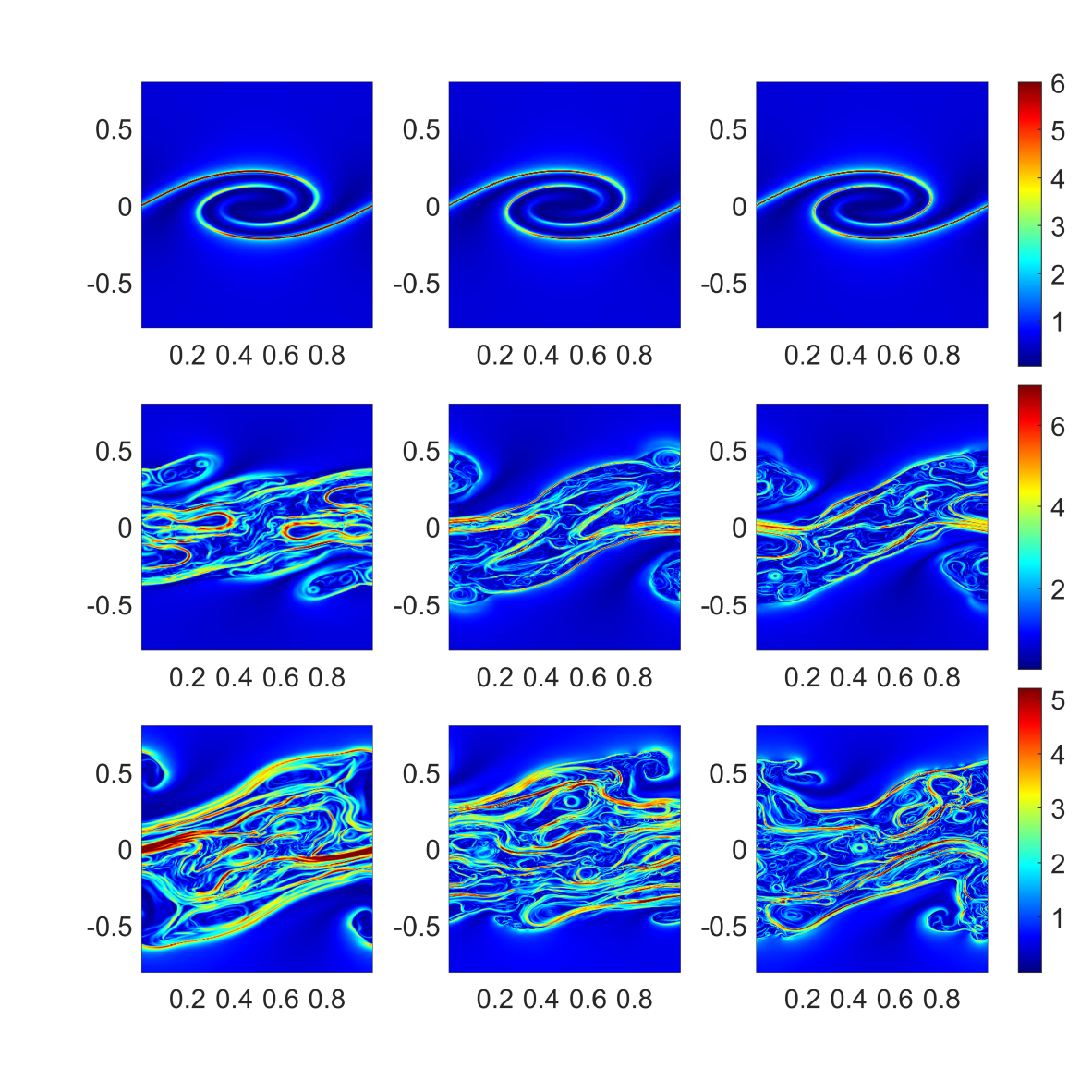}
    \caption{Example 5.4: Kelvin-Helmholtz instability. The numerical solution of $B_p/B_t$ on $N_x\times N_y = 256\times 512$ meshes. From top to bottom: $T=5,12,20$; from left to right: $k=1,2,3$. No limiter is applied.}
    \label{figKH}
\end{figure}

\begin{figure}[htbp!]
    \centering
    \subfigure[$N_x\times N_y=256\times 512$.]{
        \includegraphics[width=0.35\linewidth]{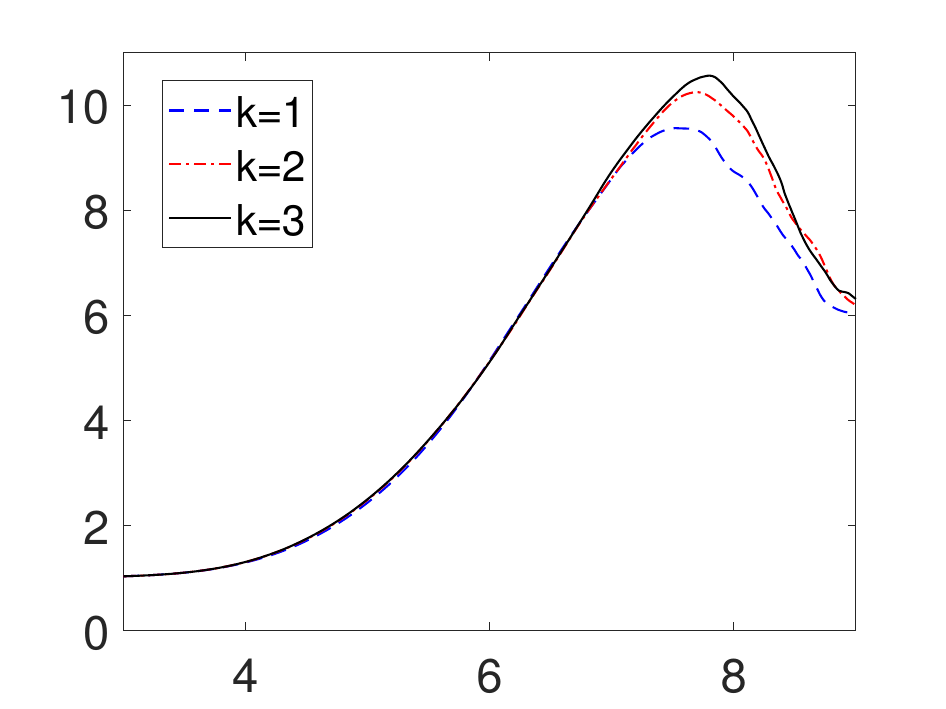}}
    \subfigure[$k=2$.]{
        \includegraphics[width=0.35\linewidth]{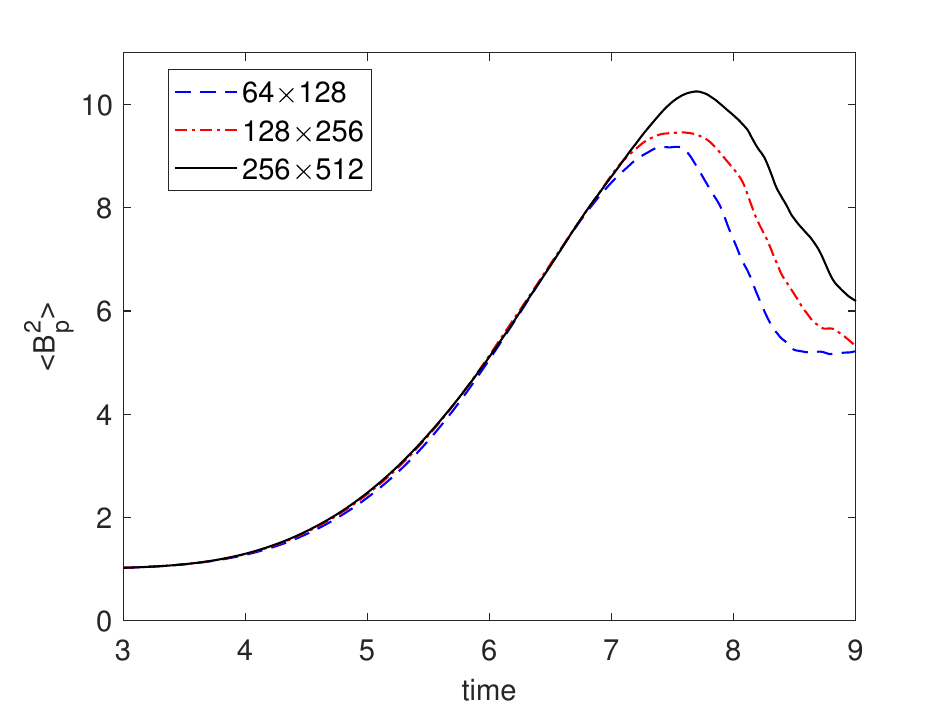}}    
    \caption{Example 5.4:  Kelvin-Helmholtz instability. The time evolution of $<B_p^2>$ with different $k$ and meshes. No limiter is applied.} 
    \label{figKH1D}
\end{figure}

\begin{figure}[htbp!]
    \centering
    \subfigure[Total entropy.]{
        \includegraphics[width=0.35\linewidth]{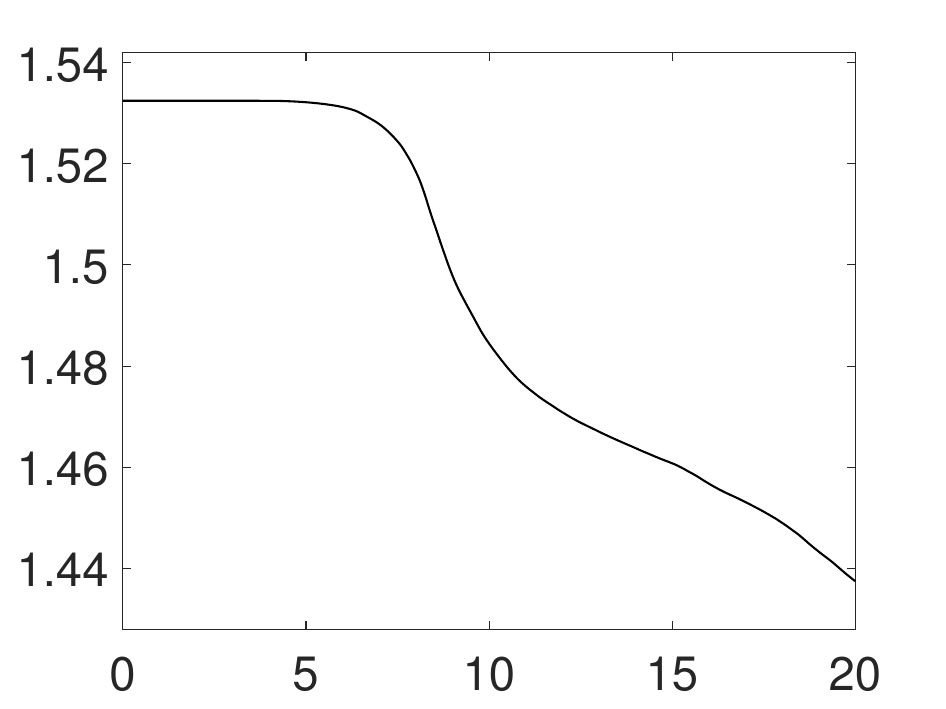}}
    \subfigure[Divergence.]{
        \includegraphics[width=0.35\linewidth]{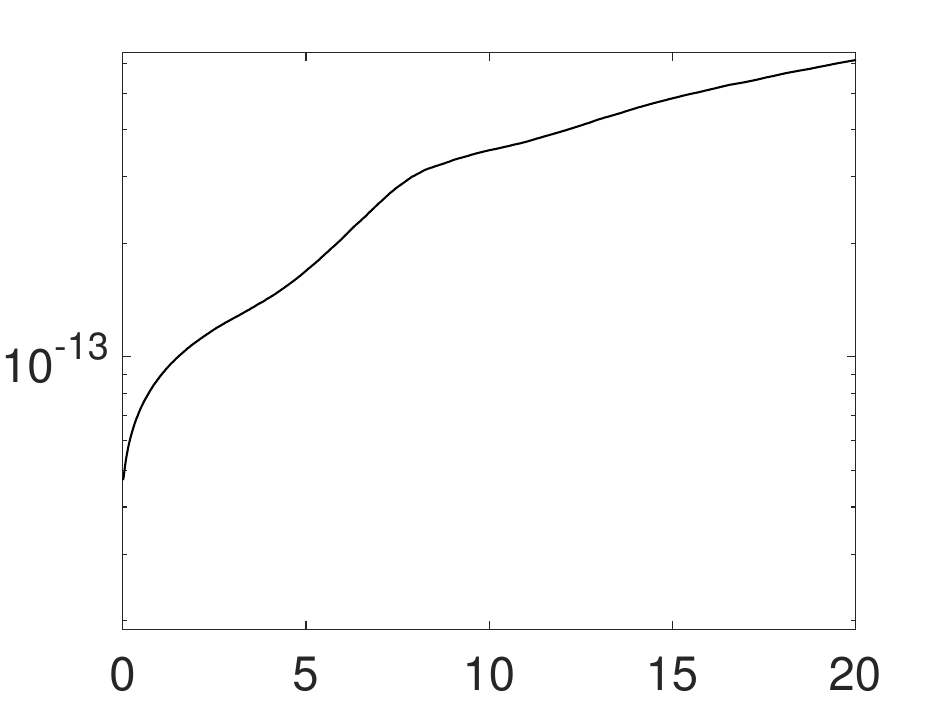}}
    \caption{Example 5.4: Kelvin-Helmholtz instability. The time evolution of total entropy and divergence norm with time on mesh with $N_x\times N_y=256\times 512$ for $k=2$. No limiter is applied.}
    \label{figKH2}
\end{figure}

We simulate this problem on mesh with $N_x\times N_y=256\times 512$ and the results of $B_p/B_t$ at $T=5, 12, 20$ with $k=1,2,3$ are shown in Fig \ref{figKH}, where $B_p=\sqrt{B_x^2+B_y^2},\  B_t=B_z$. 
In Fig \ref{figKH1D}, we plot the evolution of the poloidal magnetic energy $<B_p^2>$  with different $k$ and meshes, where 
$$ <B_p^2(t)>= \left( \int_\Omega B_p^2(t)\,\mathrm dx\mathrm dy\right)  / 
\left(\int_{\Omega} B_p^2(0)\,\mathrm dx\mathrm dy\right) . $$
For $t\le 5$, the perturbation follows a linear growth phase during which magnetic field lines wound up through the formation of a typical cat’s eye vortex structure \cite{malagoli1996nonlinear,jones1997mhd}. 
When $t\ge 8$, field amplification is eventually prevented by tearing mode instabilities, leading to reconnection events capable of expelling magnetic flux from the vortex \cite{malagoli1996nonlinear}. 
One can clearly see that small-scale structures are captured well with the proposed ES-GDF scheme. 
Further, we can see the poloidal magnetic energy grows faster in the transition to turbulence when a higher-order scheme is selected. This indicates that the overall numerical dissipation is smaller for higher-order scheme. The evolution of total entropy and divergence are shown in Fig \ref{figKH2}, indicating that the entropy is non-increase for reflective walls.\\

\textit{Example} 5.5 (Rotor). This benchmark test problem is first introduced in \cite{balsara1999staggered}, and we follow the setup in \cite{toth2000b}. It describes the spinning of a dense rotating disc of fluid in the center while ambient fluid is at rest. The magnetic field wraps around the rotating dense fluid turning it into an oblate shape. If the numerical scheme does not sufficiently control the divergence error in the magnetic field, distortion can be observed in Mach number \cite{li2005locally}.  For this problem, the limiter is applied. 

\begin{figure}[htbp!]
    \centering
    \subfigure[$200\times 200$.]{
        \includegraphics[width=0.31\linewidth]{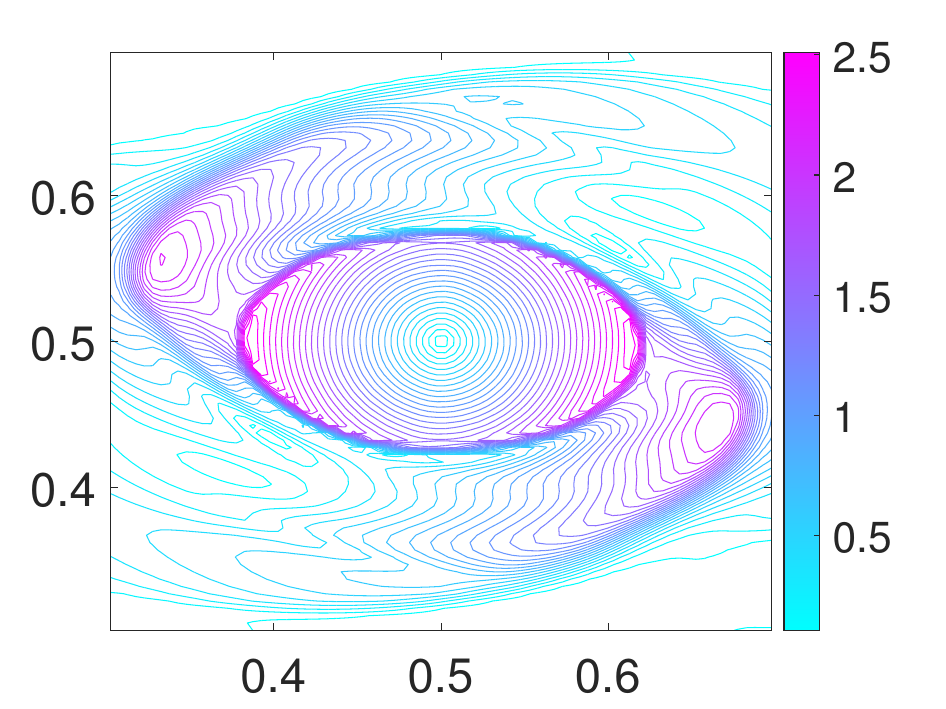}}
    \subfigure[$400\times 400$.]{
        \includegraphics[width=0.31\linewidth]{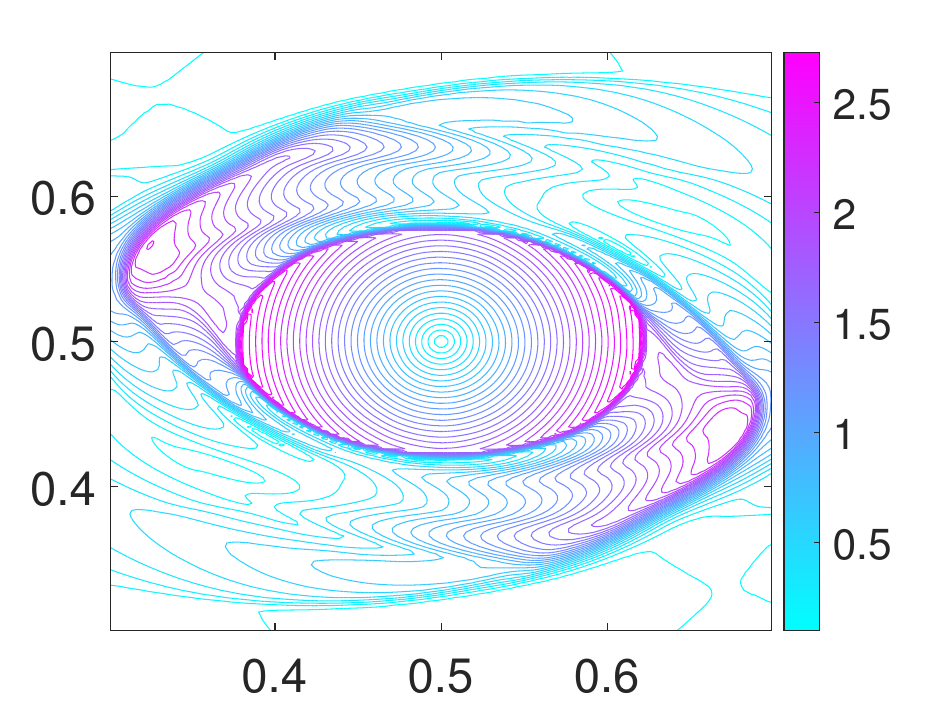}}
    \subfigure[$800\times 800$.]{
        \includegraphics[width=0.31\linewidth]{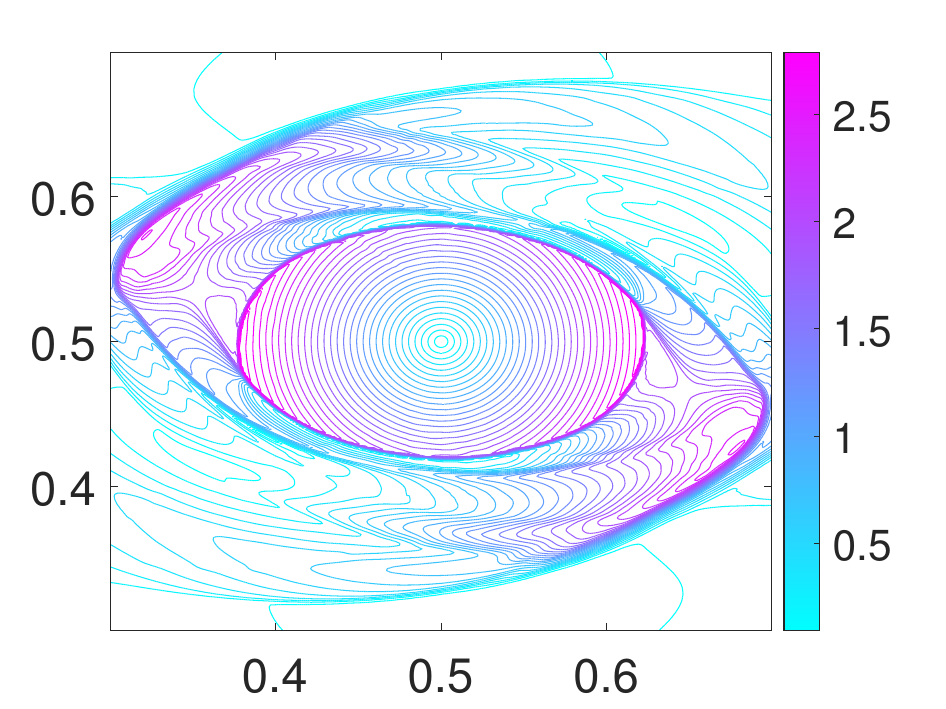}}   
    \caption{Example 5.5: Rotor. The Mach number of the central area at $T=0.295$ with different meshes. 30 contour lines are used.}
    \label{figRotor2}
\end{figure}

\begin{figure}[htbp!]
    \centering
    \subfigure[ES.]{
        \includegraphics[width=0.35\linewidth]{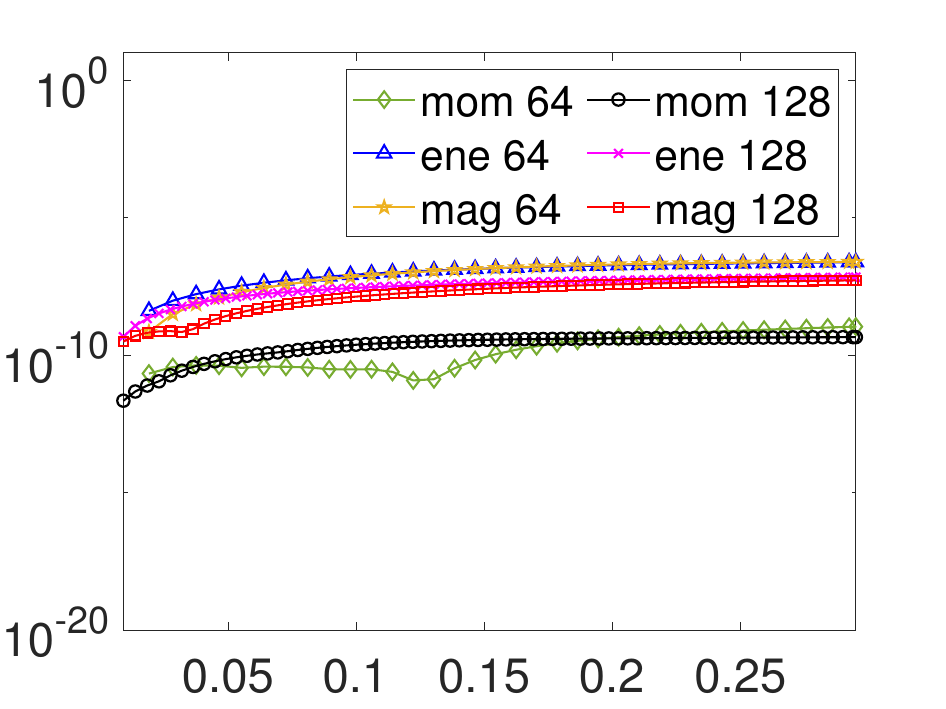}}
    \subfigure[ES-GDF.]{
        \includegraphics[width=0.35\linewidth]{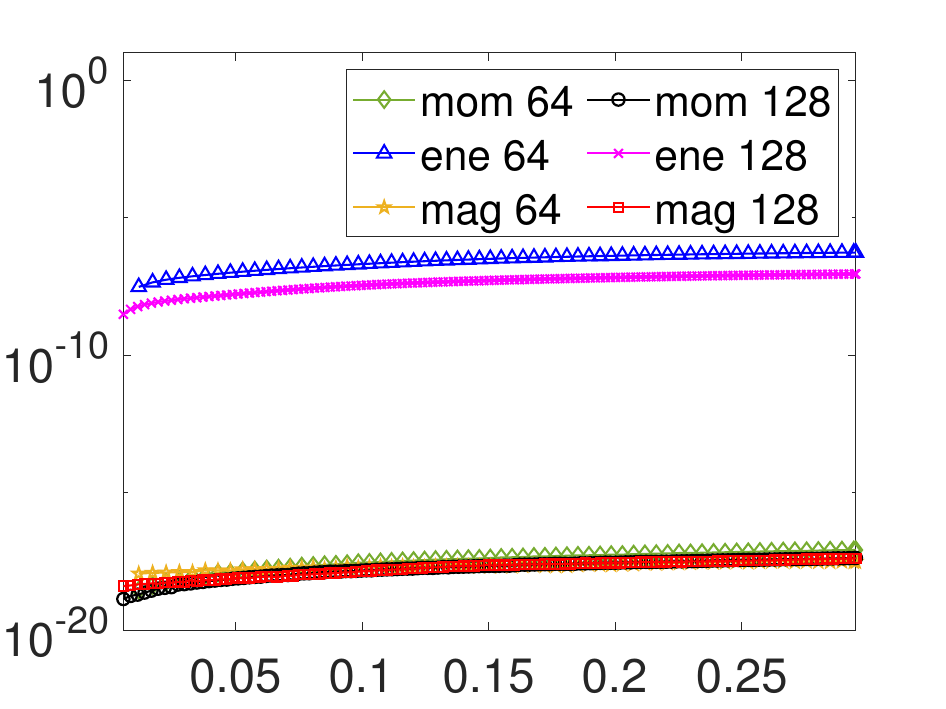}}
    \caption{Example 5.5: Rotor. The conservative errors for momentum, total energy, and magnetic field on different meshes.}
    \label{figRotorcons}
\end{figure}

\begin{figure}[htbp!]
    \centering
    \subfigure[Total entropy.]{
        \includegraphics[width=0.35\linewidth]{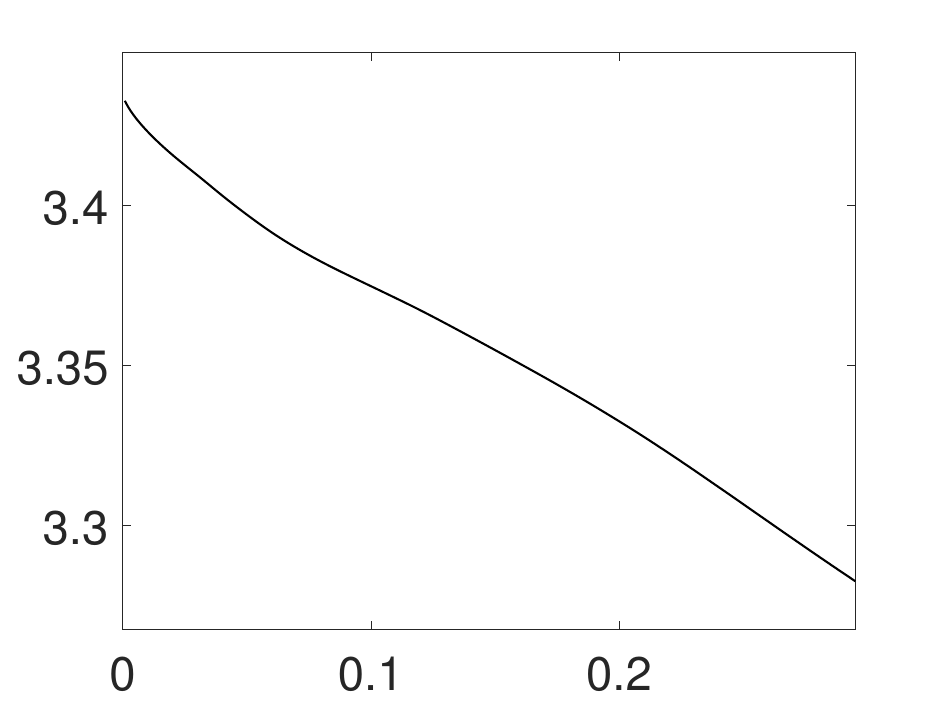}}
    \subfigure[Divergence.]{
        \includegraphics[width=0.35\linewidth]{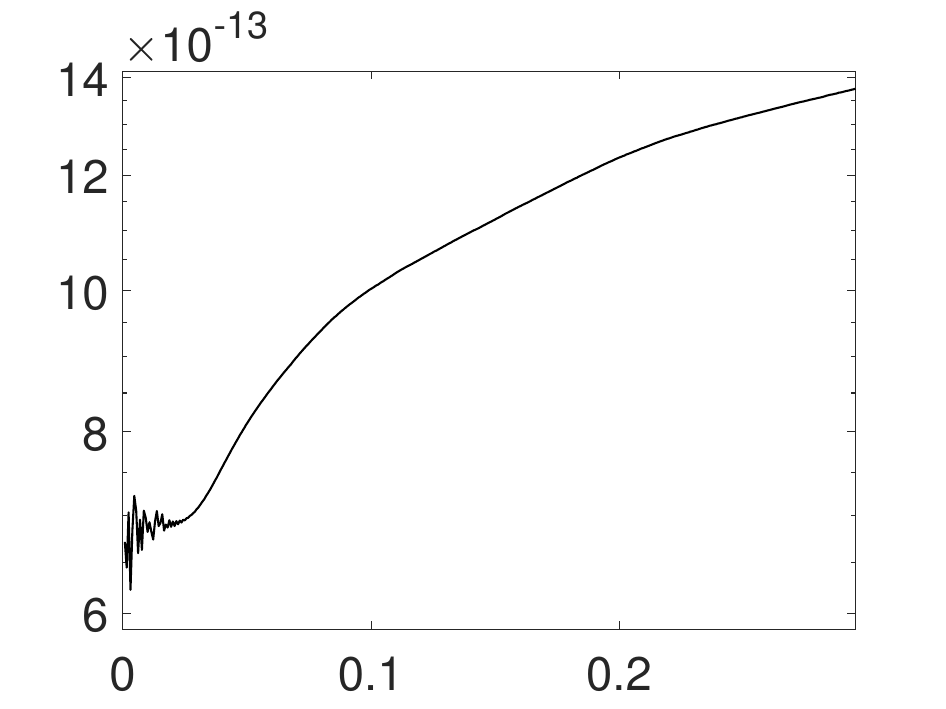}}
    \caption{Example 5.5: Rotor. The development of total entropy and divergence norm with time on $N_x\times N_y=400\times 400$ meshes.}
    \label{figRotor4}
\end{figure}

In Fig \ref{figRotor2}, we present the results of Mach number in the central area at $T=0.295$ on different meshes, and no distortion is obtained. In Fig \ref{figRotorcons}, we plot the conservative errors of ES-GDF and ES scheme. In Fig \ref{figRotor4}, we present the evolution of total entropy and divergence over time with $N_x\times N_y=400\times 400$. It can be seen that our scheme conserves the momentum and magnetic field, dissipates the entropy, and maintains the globally divergence-free property for discontinuous problems.\\

\begin{figure}[htbp!]
    \centering
    \subfigure[$\rho$.]{
        \includegraphics[width=0.31\linewidth]{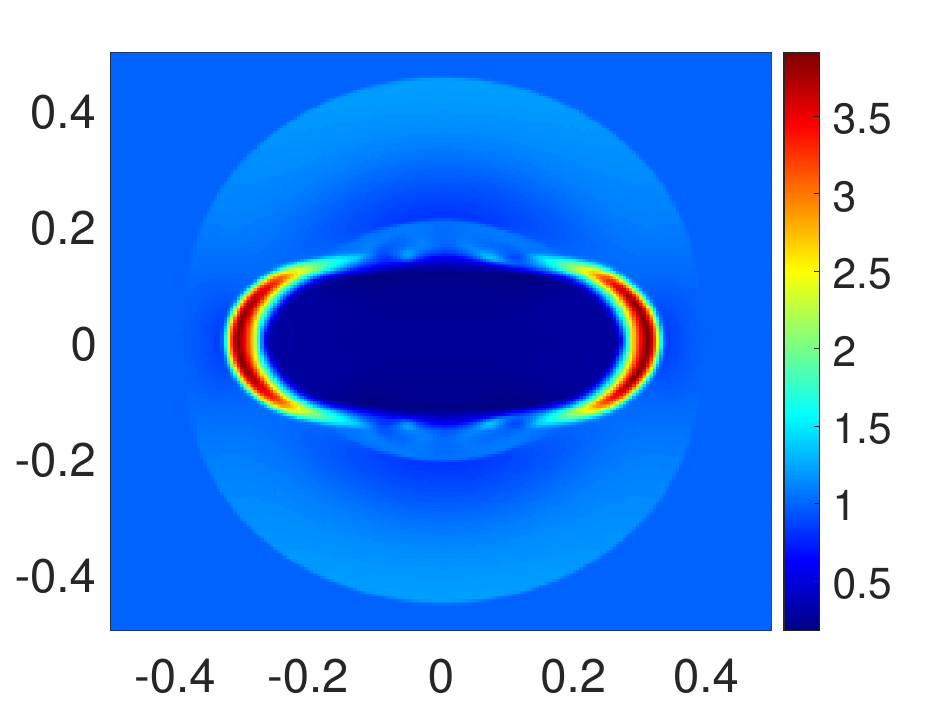}}
    \subfigure[$u_x^2+u_y^2$.]{
        \includegraphics[width=0.31\linewidth]{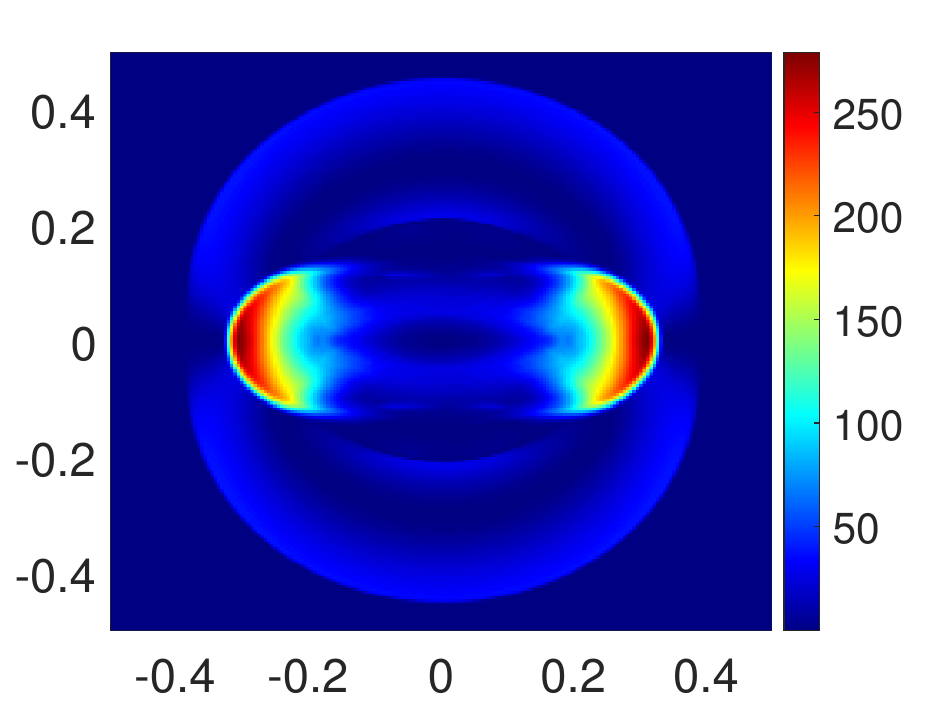}}
    \subfigure[$B_x^2+B_y^2$.]{
        \includegraphics[width=0.31\linewidth]{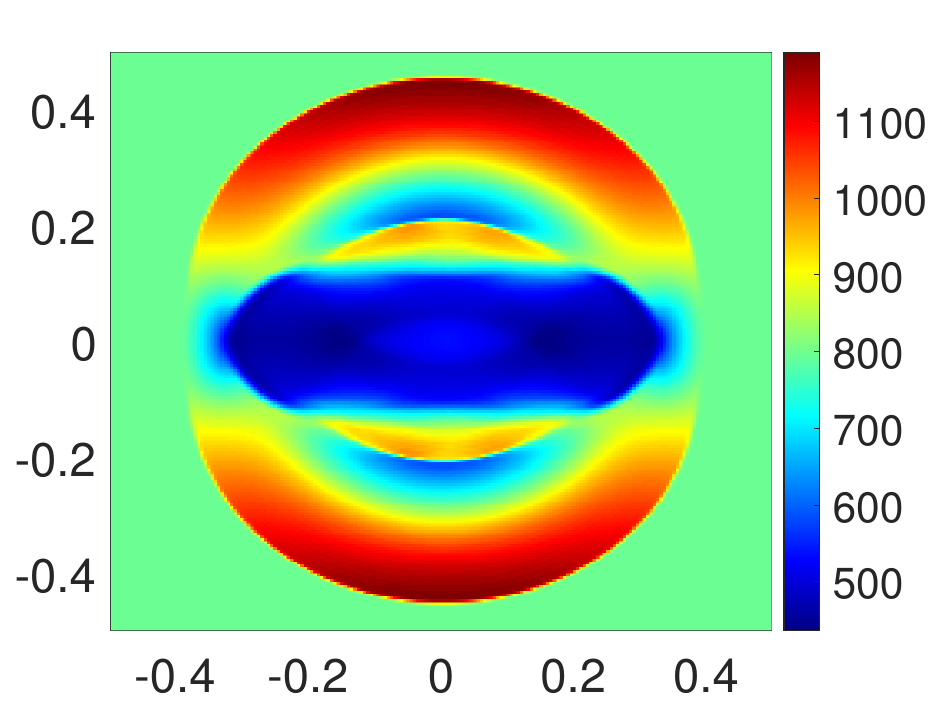}}
    \caption{Example 5.6: MHD blast. The numerical solution at $T = 0.01$ with $N_x\times N_y = 200\times 200$.}
    \label{figBlast1}
\end{figure}

\textit{Example} 5.6 (MHD blast). 
We consider the MHD blast problem \cite{balsara1999staggered}. This is a challenging test problem and is often used as a benchmark for testing the robustness of the numerical algorithms in terms of maintaining positivity of solutions. 
The fluid pulse has very small plasma beta about $2.51\times 10^{-4}$.  For this problem, the limiter is applied. 
The result at $T=0.01$ on mesh with  $N_x\times N_y=200\times 200$ is shown in Fig \ref{figBlast1}. Compared to the results in \cite{fu2018globally}, even without adding positivity limiters, we obtain non-negative pressure at all time layers, where the minimum pressure at final time is about $0.0952$. This may indicate that entropy stability can enhance computational stability for this problem. \\

\begin{figure}[htbp!]
    \centering
    \subfigure[$\ln\rho$.]{
        \includegraphics[width=0.31\linewidth]{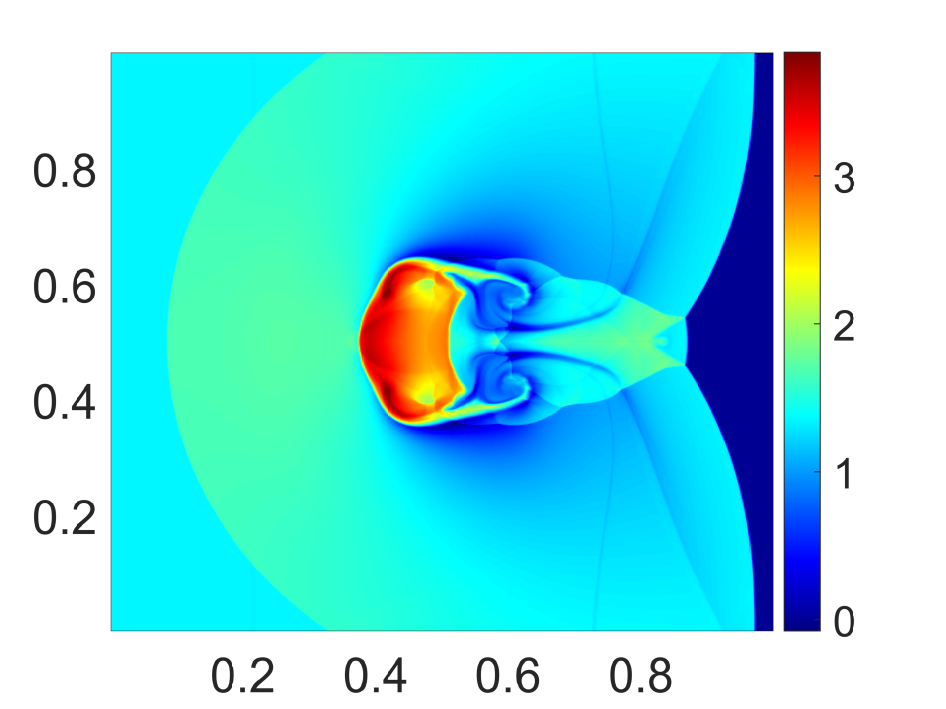}}
    \subfigure[$p$.]{
        \includegraphics[width=0.31\linewidth]{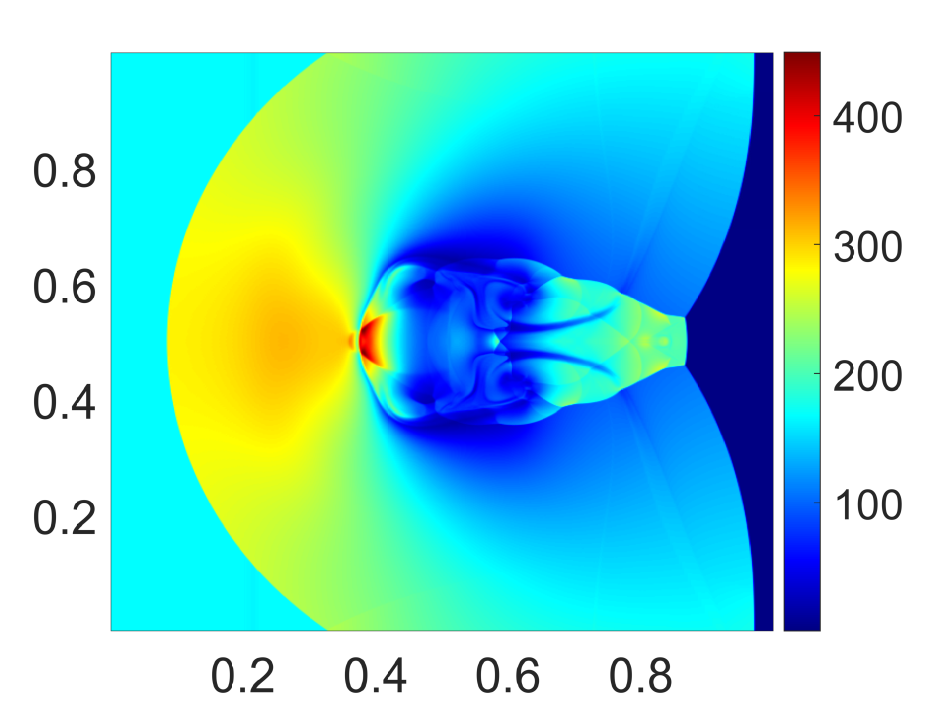}}
    \subfigure[$\left\|\mathbf B\right\|$.]{
        \includegraphics[width=0.31\linewidth]{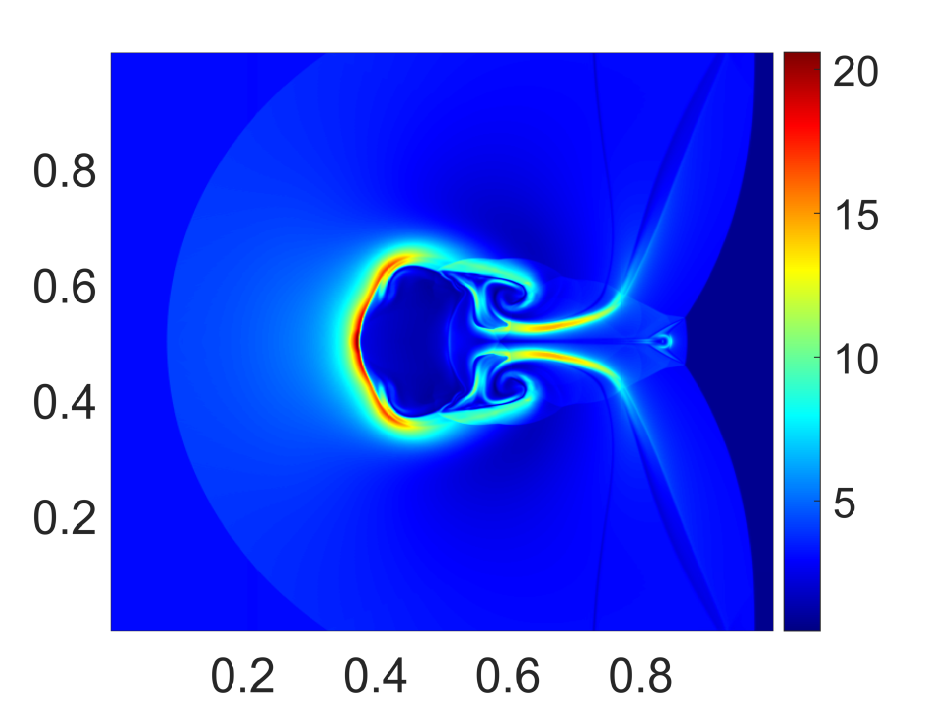}}  
    \caption{Example 5.7: Cloud shock interaction. The numerical solution at $T = 0.06$ with $N_x\times N_y = 600\times 600$.}
    \label{figCloud}
\end{figure}

\textit{Example} 5.7 (Cloud shock interaction). 
Finally, we consider the cloud shock problem \cite{rossmanith2006unstaggered}, and use the same setup in \cite{christlieb2014finite}. This problem simulates an MHD shock propagating toward a stationary bubble. 
After the shock passes through the bubble, very complex structures will appear in the computational domain. Those structures around the bubble are susceptible to numerical dissipations and low dissipative schemes are advantageous to obtain less smeared structures. In this problem, the limiter is applied. In Fig \ref{figCloud}, we present the results at $T=0.06$ with $600\times 600$. It can be seen that the ES-GDF DG scheme resolves shocks and other complex features very well without introducing negative pressure, indicates the low-dissipation and robustness of our scheme. The solution matches well with the
results in the literature, such as those in \cite{christlieb2014finite, wu2018provably, christlieb2018high}.

\section{Concluding remarks}\label{sec6}

In this paper, we present a globally divergence-free entropy stable nodal DG method that directly solves the conservative form of ideal MHD equation. By using the constraint-preserving schemes \cite{balsara2021globally} to update the magnetic field at interfaces, we can obtain the globally divergence-free internal magnetic field by a novel least-square reconstruction, where $Q^{k+1}$ basis is used. Utilizing the SBP property of stiffness matrix, we are able to prove the entropy stability of semi-discrete scheme. For problems containing strong shock, we propose a new limiting strategy, which can control the numerical oscillation and preserve the divergence-free property. 
And employing a energy correction step, we can guarantee the total entropy will not increase.
Numerical results verify that our scheme can achieve the designed accuracy, dissipate the entropy, and control the divergence norm at machine error level. 
Future works may include the extension to unstructured meshes, and the application of other equations.

\appendix

\section{Entropy conservative flux and entropy stable flux}\label{app}

\subsection{Entropy conservative flux}

For the MHD equations, Chandrashekar and Klingengberg \cite{chandrashekar2016entropy} suggested the following entropy conservative flux:
\begin{align*}
\mathbf{F}_{1}^{S}&=\hat{\rho}\bar{u}_x,
\\
\mathbf{F}_{2}^{S}&=\frac{\bar{\rho}}{2\bar{\beta}}+\bar{u}_x\mathbf{F}_{1}^{S}+\frac{1}{2}\overline{\left\| \mathbf{B} \right\| ^2}-\bar{B}_x\bar{B}_x,
\\
\mathbf{F}_{3}^{S}&=\bar{u}_y\mathbf{F}_{1}^{S}-\bar{B}_x\bar{B}_y,
\\
\mathbf{F}_{4}^{S}&=\bar{u}_z\mathbf{F}_{1}^{S}-\bar{B}_x\bar{B}_z,
\\
\mathbf{F}_{6}^{S}&=0,
\\
\mathbf{F}_{7}^{S}&=\frac{1}{\bar{\beta}}\left( \overline{\beta u_x}\bar{B}_y-\overline{\beta u_y}\bar{B}_x \right), 
\\
\mathbf{F}_{8}^{S}&=\frac{1}{\bar{\beta}}\left( \overline{\beta u_x}\bar{B}_z-\overline{\beta u_z}\bar{B}_x \right),
\\
\mathbf{F}_{5}^{S}&=\frac{1}{2}\left[ \frac{1}{\left( \gamma -1 \right) \hat{\beta}}-\overline{\left\| \mathbf{u} \right\| ^2} \right] \mathbf{F}_{1}^{S}+\bar{u}_x\mathbf{F}_{2}^{S}+\bar{u}_y\mathbf{F}_{3}^{S}+\bar{u}_z\mathbf{F}_{4}^{S}
\\
&+\bar{B}_x\mathbf{F}_{6}^{S}+\bar{B}_y\mathbf{F}_{7}^{S}+\bar{B}_z\mathbf{F}_{8}^{S}-\frac{1}{2}\bar{u}_x\overline{\left\| \mathbf{B} \right\| ^2}+\left( \bar{u}_x\bar{B}_x+\bar{u}_y\bar{B}_y+\bar{u}_z\bar{B}_z \right) \bar{B}_x, 
\end{align*}
where $\overline{(a,b)}=(a+b)/2$ is the arithmetic average, and $\hat{(\cdot)}$ is the logarithmic average of two strictly positive quantities as
$$ \hat\alpha = \frac{\alpha_r-\alpha_l}{\ln\alpha_r-\ln\alpha_l}. $$
The formula of $\mathbf G^S$ is similar.

\subsection{Entropy stable flux}

For entropy stable fluxes, it is well known that the Godunov
flux based on exact Riemann solvers is by definition entropy stable. Meanwhile, many approximate Riemann solvers, such as the HLL
or Lax–Friedrichs fluxes are also entropy stable if the estimates of left and right local wave speed $S_L,\ S_R$ satisfy $S_L\le S_L^{real},\ S_R\ge S_R^{real}$, where $S_{L,R}^{real}$ are the real minimal and maximal wave speed of the Riemann problem at the interface \cite{chen2017entropy}. In particular, the HLL flux is given by
\begin{equation*}
\hat{\mathbf{F}}^{HLL}\left( \mathbf{U}_L,\mathbf{U}_R \right) =\left\{\begin{aligned}
	&\mathbf{F}(\mathbf{U}_L),\quad S_L\ge 0,\\
	&\frac{S_R\mathbf{F}_L-S_L\mathbf{F}_R+S_RS_L\left( \mathbf{U}_R-\mathbf{U}_L \right)}{S_R-S_L},\quad S_L<0<S_R,\\
	&\mathbf{F}(\mathbf{U}_R),\quad S_R\le 0,\\
\end{aligned}\right.
\end{equation*}
If we denote the non-positive part of $S_L$ by $\mathcal S_L=\min\{S_L,0\}$, and the non-negative part of $S_R$ by $\mathcal S_R=\max\{S_R,0\}$, then the HLL flux can be written by
\begin{equation}\label{eq:HLL} \hat{\mathbf F}^{HLL}(\mathbf U_L,\mathbf U_R) = \frac{\mathcal S_R\mathbf F_L-\mathcal S_L\mathbf F_R+\mathcal S_R\mathcal S_L(\mathbf U_R-\mathbf U_L)}{\mathcal S_R-\mathcal S_L}.
\end{equation}

However, the computation of $S_L$ and $S_R$ is non-trivial.  To address this, Toro recommends the two-rarefaction approximate in \cite{toro2013riemann}, then Guermond and Popov \cite{guermond2016fast} prove that the two-rarefaction approximated wave speeds indeed provide the correct bounds for Euler equations with $1<\gamma\le 5/3$. 
However, for MHD system, the Riemann problem becomes more complex and is difficult to analyze. In \cite{bouchut2007multiwave, bouchut2010multiwave}, another ``relaxation" approach is used, and the 3-wave and 7-wave approximate Riemann solver are designed, which are proved to be entropy stable. Further, if we simply use the HLL solver \eqref{eq:HLL}, but the wavespeeds $S_L,\ S_R$ are chosen by the above 3-wave solver, then it is also entropy stable \cite{bouchut2010multiwave}. 

Specifically, for two given state $\mathbf U_L,\mathbf U_R$, the explicit formula of $S_L(\mathbf U_L,\mathbf U_R)$ and $S_R(\mathbf U_L,\mathbf U_R)$ along $x$-direction are given by 
$$
S_L=u_L-\mathcal C_L,\quad   S_R=u_R+\mathcal C_R,
$$
where
$$
\mathcal{C} _L=c_{f,L}^{0}+\alpha \left( \left( u_L-u_R \right) _++\frac{\left( p_R-p_L \right) _+}{\rho _Lc_{f,L}+\rho _Rc_{f,R}} \right), 
$$
$$
\mathcal{C} _R=c_{f,R}^{0}+\alpha \left( \left( u_L-u_R \right) _++\frac{\left( p_R-p_L \right) _+}{\rho _Lc_{f,L}+\rho _Rc_{f,R}} \right), 
$$
$$
c_{f,L}^{0}=\left\{ \frac{1}{2}\left( a_L+\frac{\left\| \mathbf{B}_L \right\| ^2}{\rho _Lx_L}+\sqrt{\left( a_L+\frac{\left\| \mathbf{B}_L \right\| ^2}{\rho _Lx_L} \right) ^2-4a_L\frac{B_{x,L}^{2}}{\rho _Lx_L}} \right) \right\} ^{1/2},
$$
$$
c_{f,R}^{0}=\left\{ \frac{1}{2}\left( a_R+\frac{\left\| \mathbf{B}_R \right\| ^2}{\rho _Rx_R}+\sqrt{\left( a_R+\frac{\left\| \mathbf{B}_R \right\| ^2}{\rho _Rx_R} \right) ^2-4a_R\frac{B_{x,R}^{2}}{\rho _Rx_R}} \right) \right\} ^{1/2},
$$and$$
\begin{aligned}
X_L=&\frac{1}{c_{f,L}}\left[ \left( u_L-u_R \right) _++\frac{\left( p_R-p_L \right) _+}{\rho _Lc_{f,L}+\rho _Rc_{f,R}} \right],
\\
X_R=&\frac{1}{c_{f,R}}\left[ \left( u_L-u_R \right) _++\frac{\left( p_R-p_L \right) _+}{\rho _Lc_{f,L}+\rho _Rc_{f,R}} \right], 
\\
x_L=&1-\frac{X_L}{1+\alpha X_L},\quad x_R=1-\frac{X_R}{1+\alpha X_R}
\end{aligned}$$
with $a=\sqrt{\gamma p/\rho},\ \alpha=(\gamma + 1)/2,\ (\cdot)_+=\max\{\cdot,0\}$, and $c_{f,L},c_{f,R}$ denotes the largest eigenvalue of $\partial\mathbf F(\mathbf U)/\partial \mathbf U$ at $\mathbf U=\mathbf U_L,\mathbf U_R$, respectively.

\section{Two-dimensional HLL Riemann solver}\label{app2} 

\begin{figure}[htb!]
    \centering
    \includegraphics[width=0.35\linewidth]{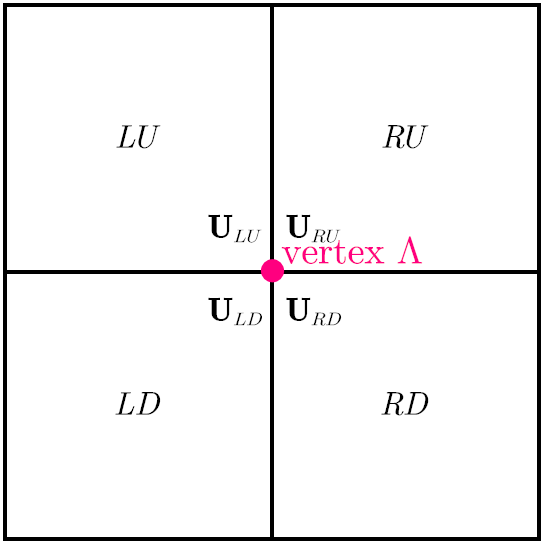}
    \caption{The notations around a vertex $\Lambda$.}
    \label{figflux}
\end{figure}

For four given stages 
at a vertex $\Lambda$, the two-dimensional HLL flux 
$\tilde{E}_z|_{\Lambda}=\tilde{E}_z\left( 
\mathbf{U}_{LD},\mathbf{U}_{LU},\mathbf{U}_{RD},\mathbf{U}_{RU} 
\right) 
$ is calculated by following \cite{chandrashekar2020constraint}. To help with the presentation, we illustrate the notations of states around a vertex $\Lambda$ in Fig \ref{figflux}.

First, the wave speeds are chosen by
$$\begin{aligned}
S_R|_\Lambda&=\max \left\{ S_R\left( \mathbf{U}_{LU},\mathbf{U}_{RU} \right) ,S_R\left( \mathbf{U}_{LD},\mathbf{U}_{RD} \right) \right\} ,
\\S_L|_\Lambda&=\min \left\{ S_L\left( \mathbf{U}_{LU},\mathbf{U}_{RU} \right) ,S_L\left( \mathbf{U}_{LD},\mathbf{U}_{RD}\right)  \right\} ,
\\
S_U|_\Lambda&=\max \left\{ S_R\left( \mathbf{U}_{RD},\mathbf{U}_{RU} \right) ,S_R\left( \mathbf{U}_{LD},\mathbf{U}_{LU} \right) \right\} ,
\\S_D|_\Lambda&=\min \left\{ S_L\left( \mathbf{U}_{RD},\mathbf{U}_{RU} \right) ,S_R\left( \mathbf{U}_{LD},\mathbf{U}_{LU} \right) \right\}\end{aligned}
$$
to match the wave speeds of the nodal DG scheme  
in the cell. Denote $$\mathcal S_R=\max\{S_R,0\},\quad \mathcal S_U=\max\{S_U,0\}, \quad \mathcal S_L=\min\{S_L,0\},\quad  \mathcal S_D=\min\{S_D,0\},$$
then the value of $\tilde E_z$ is given by
$$\begin{aligned}
\tilde E_z^{HLL}&=\frac{1}{4}\left( E_{z}^{*,R}+E_{z}^{*,L}+E_{z}^{*,U}+E_{z}^{*,D} \right) 
\\
&\quad -\frac{1}{4}\mathcal{S} _U\left( B_{x}^{U,*}-B_{x}^{**} \right) -\frac{1}{4}\mathcal{S} _D\left( B_{x}^{D,*}-B_{x}^{**} \right) 
\\
&\quad +\frac{1}{4}\mathcal{S} _R\left( B_{y}^{R,*}-B_{y}^{**} \right) +\frac{1}{4}\mathcal{S} _L\left( B_{y}^{L,*}-B_{y}^{**} \right), 
\end{aligned}
$$
where
$$
E_{z}^{*,R}=\mathbf{\hat{G}}_{6}^{HLL}\left( \mathbf{U}_{RD},\mathbf{U}_{RU} \right) , \quad 
E_{z}^{*,L}=\mathbf{\hat{G}}_{6}^{HLL}\left( \mathbf{U}_{LD},\mathbf{U}_{LU} \right), 
$$
$$
E_{z}^{*,U}=-\mathbf{\hat{F}}_{7}^{HLL}\left( \mathbf{U}_{LU},\mathbf{U}_{RU} \right) , \quad 
E_{z}^{*,D}=-\mathbf{\hat{F}}_{7}^{HLL}\left( \mathbf{U}_{LD},\mathbf{U}_{RD} \right), 
$$
and
$$\begin{aligned}
B_{x}^{**}=\frac{1}{2\Delta \mathcal S}&\left[ 
	2\mathcal{S} _R\mathcal{S} _UB_{x}^{RU}-2\mathcal{S} _L\mathcal{S} _UB_{x}^{LU}-2\mathcal{S} _R\mathcal{S} _DB_{x}^{RD}+2\mathcal{S} _L\mathcal{S} _DB_{x}^{LD}\right.
 \\ &-\mathcal{S} _R\left( E_{z}^{RU}-E_{z}^{RD} \right) +\mathcal{S} _L\left( E_{z}^{LU}-E_{z}^{LD} \right) 
 \\ &\left.-\left( \mathcal{S} _R-\mathcal{S} _L \right) \left( E_{z}^{*,U}-E_{z}^{*,D} \right) \right], 
 \end{aligned}
$$
$$\begin{aligned}
B_{y}^{**}=\frac{1}{2\Delta\mathcal S}&\left[ 
	2\mathcal{S} _R\mathcal{S} _UB_{y}^{RU}-2\mathcal{S} _L\mathcal{S} _UB_{y}^{LU}-2\mathcal{S} _R\mathcal{S} _DB_{y}^{RD}+2\mathcal{S} _L\mathcal{S} _DB_{y}^{LD}\right.\\
	&+\mathcal{S} _U\left( E_{z}^{RU}-E_{z}^{LU} \right) -\mathcal{S} _D\left( E_{z}^{RD}-E_{z}^{LD} \right)\\
	&\left.+\left( \mathcal{S} _U-\mathcal{S} _D \right) \left( E_{z}^{*,R}-E_{z}^{*,L} \right) \right], \end{aligned}
$$
where $\Delta \mathcal S=(\mathcal S_R-\mathcal S_L)(\mathcal S_U-\mathcal S_D)$.

\section{Proof of Theorem \ref{thm:ES}} \label{app3}
\begin{proof}
Firstly, we prove the local entropy conservative \eqref{eq:entropy_conser}. 
In each cell $K_{i,j}$, with the help of the semi-discrete nodal DG scheme \eqref{eq:nodalES} we have that
\begin{align*}
	\frac{\mathrm{d}}{\mathrm{d}t}&\left( \frac{\Delta x\Delta y}{4}\sum_{i_1,j_1=0}^{k+1}{\omega _{i_1}\omega _{j_1}\mathcal U_{i_1,j_1}^{i,j}} \right) =\frac{\Delta x\Delta y}{4}\sum_{i_1,j_1=0}^{k+1}{\omega _{i_1}\omega _{j_1}\left( \mathbf{V}_{i_1,j_1}^{i,j} \right) ^T\frac{\mathrm{d}\mathbf{U}_{i_1,j_1}^{i,j}}{\mathrm{d}t}}
	\\
	=&-\frac{\Delta y}{2}\sum_{j_1=0}^{k+1}\omega _{j_1}{\sum_{i_1,l=0}^{k+1}{2S_{i_1,l}}\left( \mathbf{V}_{i_1,j_1}^{i,j} \right) ^T\mathbf{F}^S\left( \mathbf{U}_{i_1,j_1}^{i,j},\mathbf{U}_{l,j_1}^{i,j} \right)}
	\\
	&+\frac{\Delta y}{2}\sum_{j_1=0}^{k+1}\omega _{j_1}{\sum_{i_1=0}^{k+1}\tau _{i_1}{\left( \left( \mathbf{V}_{i_1,j_1}^{i,j} \right) ^T\mathbf{F}_{i_1,j_1}^{i,j}-\left( \mathbf{V}_{i_1,j_1}^{i,j} \right) ^T\mathbf{F}_{i_1,j_1}^{*,i,j} \right)}}
	\\
	&-\frac{\Delta x}{2}\sum_{i_1=0}^{k+1}\omega _{i_1}{\sum_{j_1,l=0}^{k+1}{2S_{j_1,l}}\left( \mathbf{V}_{i_1,j_1}^{i,j} \right) ^T\mathbf{G}^S\left( \mathbf{U}_{i_1,j_1}^{i,j},\mathbf{U}_{i_1,l}^{i,j} \right)}
	\\
	&+\frac{\Delta x}{2}\sum_{i_1=0}^{k+1}\omega _{i_1}{\sum_{j_1=0}^{k+1}\tau _{j_1}{\left( \left( \mathbf{V}_{i_1,j_1}^{i,j} \right) ^T\mathbf{G}_{i_1,j_1}^{i,j}-\left( \mathbf{V}_{i_1,j_1}^{i,j} \right) ^T\mathbf{G}_{i_1,j_1}^{*,i,j} \right)}}
	\\
	=:&\  T_1+T_2+T_3+T_4. 
\end{align*}
 For $T_1$, utilizing the SBP property yields
 \begin{align*}
	-\frac{2}{\Delta y}T_1&=\sum_{j_1=0}^{k+1}\omega_{j_1}{\sum_{i_1,l=0}^{k+1}{\left( S_{i_1,l}+B_{i_1,l}-S_{l,i_1} \right)}\left( \mathbf{V}_{i_1,j_1}^{i,j} \right) ^T\mathbf{F}^S\left( \mathbf{U}_{i_1,j_1}^{i,j},\mathbf{U}_{l,j_1}^{i,j} \right)}
	\\
    &=\sum_{j_1=0}^{k+1}\omega_{j_1}{\sum_{i_1=0}^{k+1}{\tau _{i_1}}}\left( \mathbf{V}_{i_1,j_1}^{i,j} \right) ^T  \mathbf{F}_{i_1,j_1}^{i,j}
	\\
	&\quad +\sum_{j_1=0}^{k+1}\omega_{j_1}{\sum_{i_1,l=0}^{k+1}{\left( S_{i_1,l}-S_{l,i_1} \right)}\left( \mathbf{V}_{i_1,j_1}^{i,j} \right) ^T\mathbf{F}^S\left( \mathbf{U}_{i_1,j_1}^{i,j},\mathbf{U}_{l,j_1}^{i,j} \right)}.
\end{align*}
According to the definition of entropy conservative flux \eqref{eq:ESflux}, we get
\begin{align*}
	\sum_{j_1=0}^{k+1}\ &\omega_{j_1}{\sum_{i_1,l=0}^{k+1}{\left( S_{i_1,l}-S_{l,i_1} \right)}\left( \mathbf{V}_{i_1,j_1}^{i,j} \right) ^T\mathbf{F}^S\left( \mathbf{U}_{i_1,j_1}^{i,j},\mathbf{U}_{l,j_1}^{i,j} \right)}
	\\&=\sum_{j_1=0}^{k+1}\omega_{j_1}{\sum_{i_1,l=0}^{k+1}{S_{i_1,l}}\left( \mathbf{V}_{i_1,j_1}^{i,j}-\mathbf{V}_{l,j_1}^{i,j} \right) ^T\mathbf{F}^S\left( \mathbf{U}_{i_1,j_1}^{i,j},\mathbf{U}_{l,j_1}^{i,j} \right)}
	\\
	&=\sum_{j_1=0}^{k+1}\omega_{j_1}{\sum_{i_1,l=0}^{k+1}{S_{i_1,l}}\left( \psi _{F,i_1,j_1}^{i,j}-\psi _{F,l,j_1}^{i,j} \right)}
	\\
	&\quad -\sum_{j_1=0}^{k+1}\omega_{j_1}{\sum_{i_1,l=0}^{k+1}{S_{i_1,l}}\left( \phi _{i_1,j_1}^{i,j}-\phi _{l,j_1}^{i,j} \right) \frac{B_{x,i_1,j_1}^{i,j}+B_{x,l,j_1}^{i,j}}{2}}.
\end{align*}
 For the first term, we have
$$
	\sum_{j_1=0}^{k+1}\omega_{j_1}{\sum_{i_1,l=0}^{k+1}{S_{i_1,l}}\left( \psi _{F,i_1,j_1}^{i,j}-\psi _{F,l,j_1}^{i,j} \right)}
    =-\sum_{j_1=0}^{k+1}\omega_{j_1}{\sum_{i_1=0}^{k+1}{\tau _{i_1}\psi _{F,i_1,j_1}^{i,j}}}.
$$
 For the second term, we have
\begin{align*}
	&\sum_{j_1=0}^{k+1}\ \omega_{j_1}{\sum_{i_1,l=0}^{k+1}{S_{i_1,l}}\left( \phi _{i_1,j_1}^{i,j}-\phi _{l,j_1}^{i,j} \right) \frac{B_{x,i_1,j_1}^{i,j}+B_{x,l,j_1}^{i,j}}{2}}
	\\
	&=\frac{1}{2}\sum_{j_1=0}^{k+1}\omega_{j_1}{\sum_{i_1,l=0}^{k+1}{S_{i_1,l}}\left[ \left( \phi _{i_1,j_1}^{i,j}B_{x,i_1,j_1}^{i,j}-\phi _{l,j_1}^{i,j}B_{x,l,j_1}^{i,j} \right) +\left( \phi _{i_1,j_1}^{i,j}B_{x,l,j_1}^{i,j}-\phi _{l,j_1}^{i,j}B_{x,i_1,j_1}^{i,j} \right) \right]}
	\\
	&=\frac{1}{2}\sum_{j_1=0}^{k+1}\omega_{j_1}{\sum_{i_1=0}^{k+1}{\phi _{i_1,j_1}^{i,j}B_{x,i_1,j_1}^{i,j}\sum_{l=0}^{k+1}{S_{i_1,l}}}}
	 -\frac{1}{2}\sum_{j_1=0}^{k+1}\omega_{j_1}{\sum_{l=0}^{k+1}{\phi _{l,j_1}^{i,j}B_{x,l,j_1}^{i,j}\sum_{i_1=0}^{k+1}{S_{i_1,l}}}}
	\\
	&\quad +\frac{1}{2}\sum_{j_1=0}^{k+1}\omega_{j_1}{\sum_{i_1,l=0}^{k+1}{\phi _{i_1,j_1}^{i,j}B_{x,l,j_1}^{i,j}\left( S_{i_1,l}-S_{l,i_1} \right)}}
	\\
	&=\sum_{j_1=0}^{k+1}\omega_{j_1}{\sum_{i_1,l=0}^{k+1}{\phi _{i_1,j_1}^{i,j}B_{x,l,j_1}^{i,j}S_{i_1,l}}}
    -\sum_{j_1=0}^{k+1}\omega_{j_1}{\sum_{i_1=0}^{k+1}{\tau _{i_1}\phi _{i_1,j_1}^{i,j}B_{x,i_1,j_1}^{i,j}}}.
\end{align*}
 By the divergence-free property \eqref{eq:divfree2},
\begin{align*}
	-&\frac{2}{\Delta y}T_1
    =\sum_{j_1=0}^{k+1}\omega_{j_1}{\sum_{i_1=0}^{k+1}{\tau _{i_1}}}\left( \mathbf{V}_{i_1,j_1}^{i,j} \right) ^T\mathbf{F}_{i_1,j_1}^{i,j}
    -\sum_{j_1=0}^{k+1}\omega_{j_1}{\sum_{i_1=0}^{k+1}{\tau _{i_1}\psi _{F,i_1,j_1}^{i,j}}}
    \\
    &=\sum_{j_1=0}^{k+1}\omega_{j_1} \left( \left( \mathbf{V}_{k+1,j_1}^{i,j} \right) ^T\mathbf{F}_{k+1,j_1}^{i,j} - \left( \mathbf{V}_{0,j_1}^{i,j} \right) ^T\mathbf{F}_{0,j_1}^{i,j} \right) 
    -\sum_{j_1=0}^{k+1}\omega_{j_1}(\psi_{F,k+1,j_1}^{i,j}-\psi_{F,0,j_1}^{i,j} ) 
	\\
	&-\sum_{j_1=0}^{k+1}\omega_{j_1}{\sum_{i_1,l=0}^{k+1}{\phi _{i_1,j_1}^{i,j}B_{x,l,j_1}^{i,j}S_{i_1,l}}}
    +\sum_{j_1=0}^{k+1}\omega_{j_1}(\phi _{k+1,j_1}^{i,j}B_{x,j_1}^{i+1/2,j} -\phi _{0,j_1}^{i,j}B_{x,j_1}^{i-1/2,j} ).
\end{align*}
 Hence,
$$
    T_1+T_2=-\frac{\Delta y}{2}\sum_{j_1=0}^{k+1}{\omega _{j_1}\left( \mathcal{F} _{k+1,j_1}^{*,i,j}-\mathcal{F} _{0,j_1}^{*,i,j} \right)}
    +\frac{\Delta y}{2}\sum_{j_1=0}^{k+1}{\omega _{j_1}\sum_{i_1,l=0}^{k+1}{\phi _{i_1,j_1}^{i,j}B_{x,l,j_1}^{i,j}S_{i_1,l}}}.
$$
Similarly,
$$
    T_3+T_4=-\frac{\Delta x}{2}\sum_{i_1=0}^{k+1}{\omega _{i_1}\left( \mathcal{G} _{i_1,k+1}^{*,i,j}-\mathcal{G} _{i_1,0}^{*,i,j} \right)}
    +\frac{\Delta x}{2}\sum_{i_1=0}^{k+1}{\omega _{i_1}\sum_{j_1,l=0}^{k+1}{\phi _{i_1,j_1}^{i,j}B_{y,i_1,l}^{i,j}S_{j_1,l}}}.
$$
 Again, using the divergence-free property \eqref{eq:divfree1} yields
 $$\begin{aligned}
	\frac{\Delta y}{2}&\sum_{j_1=0}^{k+1}\,{\omega _{j_1}\sum_{i_1,l=0}^{k+1}{\phi _{i_1,j_1}^{i,j}B_{x,l,j_1}^{i,j}S_{i_1,l}}}+\frac{\Delta x}{2}\sum_{i_1=0}^{k+1}{\omega _{i_1}\sum_{j_1,l=0}^{k+1}{\phi _{i_1,j_1}^{i,j}B_{y,i_1,l}^{i,j}S_{j_1,l}}}
	\\
	&=\frac{\Delta x\Delta y}{4}\sum_{i_1,j_1=0}^{k+1}{\omega _{i_1}\omega _{j_1}\phi _{i_1,j_1}^{i,j}\sum_{l=0}^{k+1}{\left( \frac{2}{\Delta x}D_{i_1,l}B_{x,l,j_1}^{i,j}+\frac{2}{\Delta y}D_{j_1,l}B_{y,i_1,l}^{i,j} \right)}}=0.
	\end{aligned}
	$$
Therefore,
 \begin{equation}\begin{aligned}\label{eq:EC}
	\frac{\mathrm{d}}{\mathrm{d}t}&\left( \frac{\Delta x\Delta y}{4}\sum_{i_1,j_1=0}^{k+1}{\omega _{i_1}\omega _{j_1}\mathcal U_{i_1,j_1}^{i,j}} \right) 
 =T_1+T_2+T_3+T_4
	\\
	&=-\frac{\Delta y}{2}\sum_{j_1=0}^{k+1}{\omega _{j_1}\left( \mathcal{F} _{k+1,j_1}^{*,i,j}-\mathcal{F} _{0,j_1}^{*,i,j} \right)}-\frac{\Delta x}{2}\sum_{i_1=0}^{k+1}{\omega _{i_1}\left( \mathcal{G} _{i_{1},k+1}^{*,i,j}-\mathcal{G} _{i_1,0}^{*,i,j} \right)}.\end{aligned}
	\end{equation}
 Next, we prove the global entropy stable. At interface $I_{i+1/2,j}^y$ inside the computational domain $\Omega$, by using the definition of entropy stable flux,  we can obtain
 \begin{equation}\label{eq:FhatES}\begin{aligned}
	-\left( \mathcal{F} _{k+1,j_1}^{*,i,j}\right.&\left.-\ \mathcal{F} _{0,j_1}^{*,i+1,j}\right)
    \\&=\left( \mathbf{V}_{0,j_1}^{i+1,j}-\mathbf{V}_{k+1,j_1}^{i,j} \right) ^T\mathbf{\hat{F}}^{i+1/2,j}_{j_1} -\left( \psi _{F,0,j_1}^{i+1,j}-\psi _{F,k+1,j_1}^{i,j} \right) 
    \\&\quad +\frac{B_{x,0,j_1}^{i+1,j}+B_{x,k+1,j_1}^{i,j}}{2}\left( \phi _{0,j_1}^{i+1,j}-\phi _{k+1,j_1}^{i,j} \right) 
	\\
	&\le 0.\end{aligned}
\end{equation}
If the boundaries are periodic or compactly supported, then \eqref{eq:FhatES} also hold for the boundaries of $\Omega$. 
Similarly, $-\left(\mathcal{G} _{i_1,k+1}^{*,i,j}-\mathcal{G} _{i_1,0}^{*,i,j+1}\right)\le 0.$
Furthermore, summing \eqref{eq:EC} over $i,j$ yields
 $$
	 \frac{\mathrm{d}}{\mathrm{d}t}\left( \frac{\Delta x\Delta y}{4}\sum_{i,j=1}^{N_x,N_y}{\sum_{i_1,j_1=0}^{k+1}{\omega _{i_1}\omega _{j_1}\mathcal U_{i_1,j_1}^{i,j}}} \right) \le 0.
	$$

     For reflective walls, without loss of generality, we assume the left boundary $x_{1/2}=a$ is a reflective wall, i.e.
 $$
\mathbf U_{int}=(\rho,u_x,u_y,u_z,\mathcal E,B_x,B_y,B_z)^T, 
 \quad \mathbf U_{out}=(\rho,-u_x,u_y,u_z,\mathcal E,-B_x,B_y,B_z)^T. $$
We can verify that $S_R=-S_L=:S$, thus the HLL flux is equivalent to the Lax-Friedrichs flux. Notice that the condition \eqref{eq:divfree2} leads to $ B_x=0 $, hence
$$\begin{aligned}
\mathbf{\hat{F}}(\mathbf U_{out},\mathbf U_{int})&=\frac{1}{2}\left[ \mathbf{F}\left( \mathbf{U}_{out} \right) +\mathbf{F}\left( \mathbf{U}_{int} \right) \right] -\frac{S}{2}\left( \mathbf{U}_{int}-\mathbf{U}_{out} \right) 
\\&=\left[ 0,p^{\star}-\rho u_{x}S,0,0,0,0,0,0 \right] ^T,
\end{aligned}$$
and 
$$
\mathbf{V}_{int}-\mathbf{V}_{out}=4\beta u_x,\quad \psi_{int}-\psi _{out}=2\rho u_{x}+2\beta u_{x}\left\| \mathbf{B} \right\| ^2.
$$
Then, it can be calculated that
$$\begin{aligned}-\left(\mathcal F^{*,int}-\mathcal F^{*,out}\right)&=
\left( \mathbf{V}_{int}-\mathbf{V}_{out} \right) ^T\mathbf{\hat{F}}\left( \mathbf{U}_{out},\mathbf{U}_{int} \right) -\left( \psi _{int}-\psi _{out} \right) 
\\
&=4\frac{\rho}{2p}u_x\left( p+\frac{1}{2}\left\| \mathbf{B} \right\| ^2-\rho u_xS \right) -2\rho u_x+2\frac{\rho}{2p}u_x\left\| \mathbf{B} \right\| ^2
\\
&=-4\beta \rho u_x^2S\le 0,
\end{aligned}
$$
i.e. the entropy dissipate property \eqref{eq:FhatES} also holds for reflective walls.
And the entropy dissipate on the whole domain holds as well. 
\end{proof}

\section{The reconstruction for small $k$}\label{app4} 

We give a brief discussion about the LS reconstruction problem in Section 3.3 for small $k$. We rewrite the linear system (3.18) as
$$
\mathbf{G}\left[ \begin{array}{c}
	\mathbf{B}^{rec}\\
	\boldsymbol{\lambda }\\
\end{array} \right] =\left[ \begin{array}{c}
	M_2\mathbf{\tilde{\mathbf B}}\\
	\mathbf{b}_1\\
\end{array} \right], 
$$
then $\mathbf G^{-1}$ can be given by
$$
\mathbf{G}^{-1}=\left[ \begin{matrix}
	M^{-1}-M^{-1}A_1^T\mathbf{S}A_1M^{-1}&		M^{-1}A_1^T\mathbf{S}\\
	-\mathbf{S}A_1M^{-1}&		\mathbf{S}\\
\end{matrix} \right], 
$$
where $\mathbf S=(A_1MA_1^T)^{-1}$ is the Schur complement of $M$. In the following, we will assume $\Delta x=\Delta y$.

To help understand, we first consider the simplest case $k=0$, and the 2 points Gauss-Lobatto quadrature is used. Now we have
$$
M=\left[ \begin{matrix}
	1&		0\\
	0&		1\\
\end{matrix} \right] =I,\quad D=\left[ \begin{matrix}
	0.5&		-0.5\\
	0.5&		-0.5\\
\end{matrix} \right]. 
$$
Assume $$ b_x^+ = a_0^+,\quad b_x^- = a_0^-,\quad b_y^+ = b_0^+,\quad b_y^- = b_0^-. $$
Then, the original linear system of reconstruction (3.16) is
\begin{equation}\label{eq:P0-1}
\left[ \begin{matrix}
	-1&		1&		0&		0&		-1&		0&		1&		0\\
	-1&		1&		0&		0&		0&		-1&		0&		1\\
	0&		0&		-1&		1&		-1&		0&		1&		0\\
	0&		0&		-1&		1&		0&		-1&		0&		1\\
	0&		1&		0&		0&		0&		0&		0&		0\\
	0&		0&		0&		1&		0&		0&		0&		0\\
	1&		0&		0&		0&		0&		0&		0&		0\\
	0&		0&		1&		0&		0&		0&		0&		0\\
	0&		0&		0&		0&		0&		0&		1&		0\\
	0&		0&		0&		0&		0&		0&		0&		1\\
	0&		0&		0&		0&		1&		0&		0&		0\\
	0&		0&		0&		0&		0&		1&		0&		0\\
\end{matrix} \right] \left[ \begin{matrix}
	B_{x,00}\\
	B_{x,10}\\
	B_{x,01}\\
	B_{x,11}\\
	B_{y,00}\\
	B_{y,10}\\
	B_{y,01}\\
	B_{y,11}\\
\end{matrix} \right] =\left[ \begin{matrix}
	0\\
	0\\
	0\\
	0\\
	a_{0}^{+}\\
	a_{0}^{+}\\
	a_{0}^{-}\\
	a_{0}^{-}\\
	b_{0}^{+}\\
	b_{0}^{+}\\
	b_{0}^{-}\\
	b_{0}^{-}\\
\end{matrix} \right],
\end{equation}
which is a $12\times 8$ linear system. To simplify the form, we have multiplied a factor 2 for the first four rows. It is notable that the last 8 rows of $A$ already form a full rank $8\times 8$ matrix, hence the system \eqref{eq:P0-1} has at most one solution. Since the number of equations is larger than DOFs, it seems that there may have ``\,$0 = d$\," rows with $d\ne 0$ in the simplified stepped system, which will lead to the system non-solvable. However, applying some elementary row operations to $A$ yields

$$
\left[ \begin{matrix}
	0&		0&		0&		0&		0&		0&		0&		0\\
	0&		0&		0&		0&		0&		0&		0&		0\\
	0&		0&		0&		0&		0&		0&		0&		0\\
	0&		0&		0&		0&		0&		0&		0&		0\\
	0&		1&		0&		0&		0&		0&		0&		0\\
	0&		0&		0&		1&		0&		0&		0&		0\\
	1&		0&		0&		0&		0&		0&		0&		0\\
	0&		0&		1&		0&		0&		0&		0&		0\\
	0&		0&		0&		0&		0&		0&		1&		0\\
	0&		0&		0&		0&		0&		0&		0&		1\\
	0&		0&		0&		0&		1&		0&		0&		0\\
	0&		0&		0&		0&		0&		1&		0&		0\\
\end{matrix} \right] \left[ \begin{matrix}
	B_{x,00}\\
	B_{x,10}\\
	B_{x,01}\\
	B_{x,11}\\
	B_{y,00}\\
	B_{y,10}\\
	B_{y,01}\\
	B_{y,11}\\
\end{matrix} \right] =\left[ \begin{matrix}
	a_{0}^{-}-a_{0}^{+}+b_{0}^{-}-b_{0}^{+}\\
	a_{0}^{-}-a_{0}^{+}+b_{0}^{-}-b_{0}^{+}\\
	a_{0}^{-}-a_{0}^{+}+b_{0}^{-}-b_{0}^{+}\\
	a_{0}^{-}-a_{0}^{+}+b_{0}^{-}-b_{0}^{+}\\
	a_{0}^{+}\\
	a_{0}^{+}\\
	a_{0}^{-}\\
	a_{0}^{-}\\
	b_{0}^{+}\\
	b_{0}^{+}\\
	b_{0}^{-}\\
	b_{0}^{-}\\
\end{matrix} \right]. 
$$
Hence, the system is solvable if and only if the cell-average constraint is satisfied, which indicates that $$a_0^+-a_0^-+b_0^+-b_0^-=0.$$
Then, the unique solution of \eqref{eq:P0-1} is
\begin{equation}\begin{aligned}\label{eq:P0-solution}
B_{x,00}=B_{x,01}=a_{0}^{-},&\quad B_{x,10}=B_{x,11}=a_{0}^{+},
\\
B_{y,00}=B_{y,10}=b_{0}^{-},&\quad B_{y,01}=B_{y,11}=b_{0}^{+}.
\end{aligned}
\end{equation}
On the other hand, if we directly consider the system (3.18), then we can see $A_1$ is a $8\times 8$ full rank matrix. Since $M=I$, we can derive that $\mathbf S=(A_1^{T})^{-1}A_1^{-1}$, and
$$
\mathbf{G}^{-1}=\left[ \begin{matrix}
	O&		A_1^{-1}\\
	-\left( A_1^T \right) ^{-1}&		(A_1^T)^{-1}A_1^{-1}\\
\end{matrix} \right]. 
$$
Therefore, we have
$$
\left[ \begin{array}{c}
	\mathbf{B}^{rec}\\
	\mathbf{\lambda }\\
\end{array} \right] =\left[ \begin{matrix}
	O&		A_1^{-1}\\
	-\left( A_1^T \right) ^{-1}&		(A_1^T)^{-1}A_1^{-1}\\
\end{matrix} \right] \left[ \begin{array}{c}
	\mathbf{B}^{\mathrm{DG}}\\
	\mathbf{b}_1\\
\end{array} \right], 
$$
so $\mathbf B^{rec}= A_1^{-1}\mathbf b_1$, it is also the solution of $A\mathbf B=\mathbf b$, thus we recover \eqref{eq:P0-solution}. This is also consistent with the conclusion in \cite{balsara2021globally}.

\begin{figure}[htbp!]
    \centering
    \subfigure[$k=1$.]{
        \includegraphics[width=0.31\linewidth]{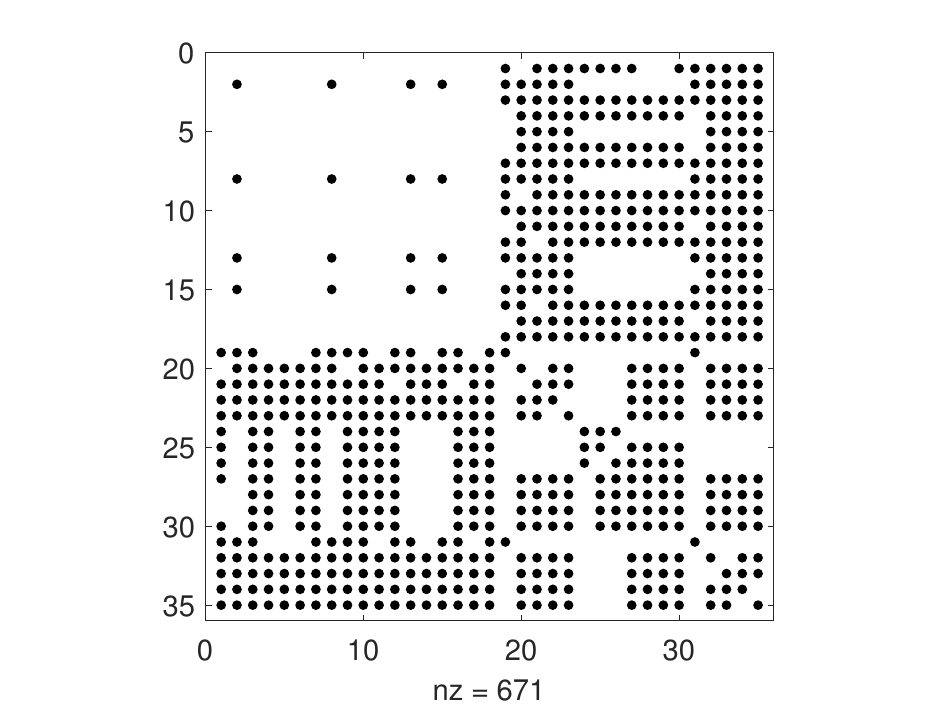}}
    \subfigure[$k=2$.]{
        \includegraphics[width=0.31\linewidth]{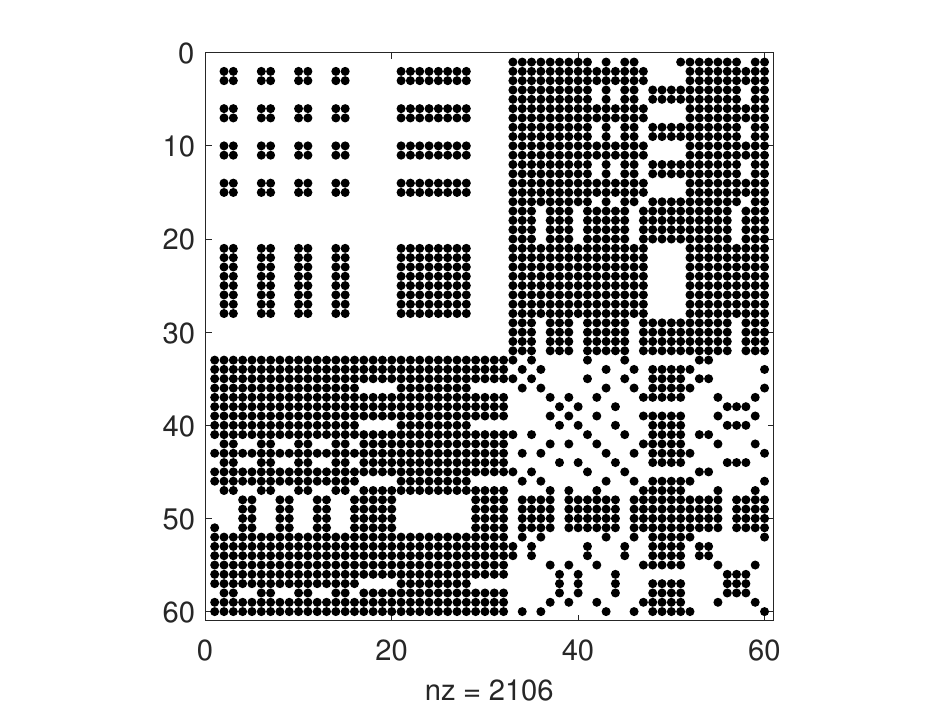}}
    \subfigure[$k=3$.]{
        \includegraphics[width=0.31\linewidth]{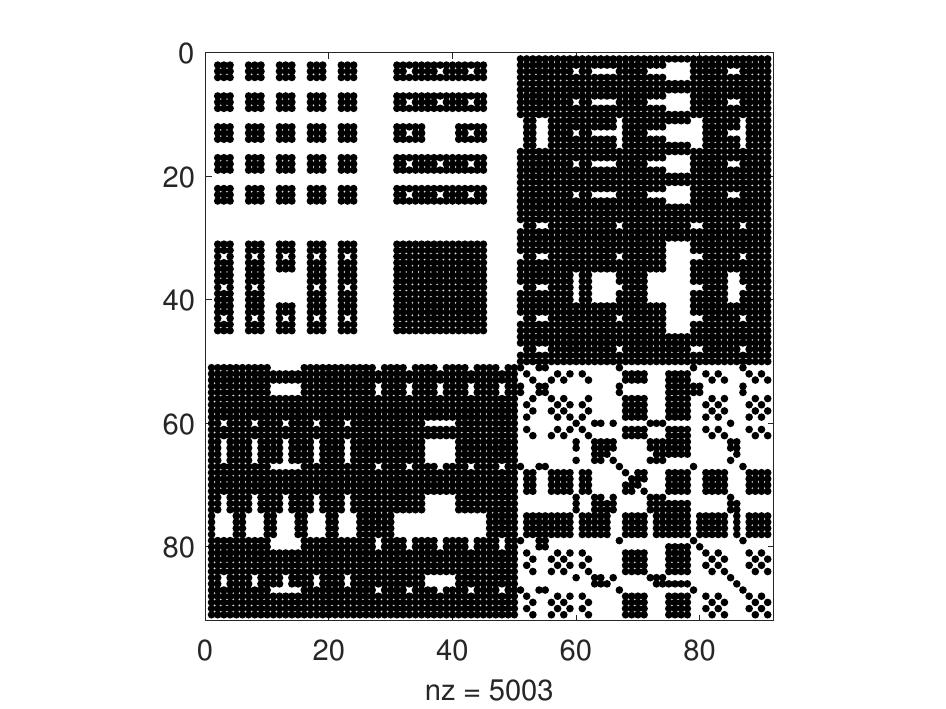}}
    \caption{The structure of nonzero elements in $\mathbf G^{-1}$. }
    \label{figmat}
\end{figure}

However, for $k> 0$, the system has become very complicated.  When $k=1$, the size of $A$ is a $21\times 18$, and the size of $\mathbf G$ is $35\times 35$. For $k=2$, the sizes are respectively $32\times 32$ and $60\times 60$. For $k=3$, the sizes are respectively $ 45\times 50 $ and $91\times 91$. Therefore, it is difficult to study the underlying mechanism or give the explicit formulation like \eqref{eq:P0-solution} for the reconstruction equation with general $k$. In Fig \ref{figmat}, we present the structure of nonzero elements in $\mathbf G^{-1}$ for $k=1,2,3$ as a reference. 

In addition, we give the full express of $\mathbf G^{-1}$ for $k=1$. The first 5 columns of $\mathbf G^{-1}$ are\newpage 
 $$\left[ \begin{matrix} 
 0 &0 &0 &0 &0 \\ 
0 &9/16 &0 &0 &0 \\ 
0 &0 &0 &0 &0 \\ 
0 &0 &0 &0 &0 \\ 
0 &0 &0 &0 &0 \\ 
0 &0 &0 &0 &0 \\ 
0 &0 &0 &0 &0 \\ 
0 &-9/16 &0 &0 &0 \\ 
0 &0 &0 &0 &0 \\ 
0 &0 &0 &0 &0 \\ 
0 &0 &0 &0 &0 \\ 
0 &0 &0 &0 &0 \\ 
0 &-9/16 &0 &0 &0 \\ 
0 &0 &0 &0 &0 \\ 
0 &9/16 &0 &0 &0 \\ 
0 &0 &0 &0 &0 \\ 
0 &0 &0 &0 &0 \\ 
0 &0 &0 &0 &0 \\ 
446/10487 &-361/4426 &446/10487 &0 &0 \\ 
0 &-407/5879 &1345/18033 &686/18395 &-407/5879 \\ 
194/2497 &-811/8786 &158/7263 &2699/54280 &-811/8786 \\ 
-293/5856 &482/5241 &-293/5856 &293/2928 &-403/2191 \\ 
31/13608 &-236/56361 &31/13608 &-31/6804 &225/26867 \\ 
-195/1351 &0 &195/1351 &-195/1351 &0 \\ 
1262/6927 &0 &-1262/6927 &-1217/3340 &0 \\ 
5143/432011 &0 &-5143/432011 &-10286/432011 &0 \\ 
-3599/5453 &0 &937/6487 &-115/589 &0 \\ 
0 &0 &-137/1263 &2539/11975 &0 \\ 
0 &0 &336/1285 &609/3007 &0 \\ 
-99/90076 &0 &-1471/2543 &513/7487 &0 \\ 
-338/963 &-265/724 &-338/963 &0 &0 \\ 
152/7003 &181/7664 &152/7003 &-304/7003 &-181/3832 \\ 
-883/3115 &-223/723 &-883/3115 &847/1494 &446/723 \\ 
89/1024 &884/2323 &213/499 &737/2869 &884/2323 \\ 
-413/985 &-469/1397 &-203/5977 &-526/2321 &-469/1397 \\ 
 \end{matrix} \right]. $$ 

\newpage
Column 6 to 10 of $\mathbf G^{-1}$ are
  $$\left[ \begin{matrix} 
 0 &0 &0 &0 &0 \\ 
0 &0 &-9/16 &0 &0 \\ 
0 &0 &0 &0 &0 \\ 
0 &0 &0 &0 &0 \\ 
0 &0 &0 &0 &0 \\ 
0 &0 &0 &0 &0 \\ 
0 &0 &0 &0 &0 \\ 
0 &0 &9/16 &0 &0 \\ 
0 &0 &0 &0 &0 \\ 
0 &0 &0 &0 &0 \\ 
0 &0 &0 &0 &0 \\ 
0 &0 &0 &0 &0 \\ 
0 &0 &9/16 &0 &0 \\ 
0 &0 &0 &0 &0 \\ 
0 &0 &-9/16 &0 &0 \\ 
0 &0 &0 &0 &0 \\ 
0 &0 &0 &0 &0 \\ 
0 &0 &0 &0 &0 \\ 
0 &-446/10487 &361/4426 &-446/10487 &446/10487 \\ 
686/18395 &1345/18033 &-407/5879 &0 &-158/7263 \\ 
2699/54280 &158/7263 &-811/8786 &194/2497 &1345/18033 \\ 
293/2928 &-293/5856 &482/5241 &-293/5856 &-31/13608 \\ 
-31/6804 &31/13608 &-236/56361 &31/13608 &-293/5856 \\ 
195/1351 &-195/1351 &0 &195/1351 &-195/1351 \\ 
1217/3340 &1262/6927 &0 &-1262/6927 &-5143/432011 \\ 
10286/432011 &5143/432011 &0 &-5143/432011 &1262/6927 \\ 
51/881 &739/2742 &0 &-96/3349 &491/831 \\ 
723/6241 &5078/11975 &0 &946/2781 &-1395/8672 \\ 
-162/1069 &625/1543 &0 &-2033/3601 &-341/1901 \\ 
-809/2711 &255/1846 &0 &-77/4190 &524/3107 \\ 
0 &338/963 &265/724 &338/963 &-338/963 \\ 
-304/7003 &152/7003 &181/7664 &152/7003 &-883/3115 \\ 
847/1494 &-883/3115 &-223/723 &-883/3115 &-152/7003 \\ 
737/2869 &213/499 &884/2323 &89/1024 &-203/5977 \\ 
-526/2321 &-203/5977 &-469/1397 &-413/985 &-213/499 \\ 
 \end{matrix} \right]. $$ 

\newpage

Column 11 to 15 of $\mathbf G^{-1}$ are
 $$\left[ \begin{matrix} 
 0 &0 &0 &0 &0 \\ 
0 &0 &-9/16 &0 &9/16 \\ 
0 &0 &0 &0 &0 \\ 
0 &0 &0 &0 &0 \\ 
0 &0 &0 &0 &0 \\ 
0 &0 &0 &0 &0 \\ 
0 &0 &0 &0 &0 \\ 
0 &0 &9/16 &0 &-9/16 \\ 
0 &0 &0 &0 &0 \\ 
0 &0 &0 &0 &0 \\ 
0 &0 &0 &0 &0 \\ 
0 &0 &0 &0 &0 \\ 
0 &0 &9/16 &0 &-9/16 \\ 
0 &0 &0 &0 &0 \\ 
0 &0 &-9/16 &0 &9/16 \\ 
0 &0 &0 &0 &0 \\ 
0 &0 &0 &0 &0 \\ 
0 &0 &0 &0 &0 \\ 
0 &-446/10487 &-361/4426 &0 &361/4426 \\ 
-2699/54280 &-194/2497 &811/8786 &811/8786 &811/8786 \\ 
686/18395 &0 &-407/5879 &-407/5879 &-407/5879 \\ 
31/6804 &-31/13608 &236/56361 &-225/26867 &236/56361 \\ 
293/2928 &-293/5856 &482/5241 &-403/2191 &482/5241 \\ 
-195/1351 &-195/1351 &0 &0 &0 \\ 
10286/432011 &-5143/432011 &0 &0 &0 \\ 
-1217/3340 &1262/6927 &0 &0 &0 \\ 
551/2084 &-169/2723 &0 &0 &0 \\ 
-823/16086 &301/5142 &0 &0 &0 \\ 
-3234/139709 &461/3464 &0 &0 &0 \\ 
-572/2423 &-1841/2873 &0 &0 &0 \\ 
0 &338/963 &-265/724 &0 &265/724 \\ 
847/1494 &-883/3115 &-223/723 &446/723 &-223/723 \\ 
304/7003 &-152/7003 &-181/7664 &181/3832 &-181/7664 \\ 
-526/2321 &-413/985 &-469/1397 &-469/1397 &-469/1397 \\ 
-737/2869 &-89/1024 &-884/2323 &-884/2323 &-884/2323 \\ 
 \end{matrix} \right]. $$ 

\newpage
Column 16 to 20 of $G^{-1}$ are
 $$\left[ \begin{matrix} 
 0 &0 &0 &446/10487 &0 \\ 
0 &0 &0 &-361/4426 &-407/5879 \\ 
0 &0 &0 &446/10487 &1345/18033 \\ 
0 &0 &0 &0 &686/18395 \\ 
0 &0 &0 &0 &-407/5879 \\ 
0 &0 &0 &0 &686/18395 \\ 
0 &0 &0 &-446/10487 &1345/18033 \\ 
0 &0 &0 &361/4426 &-407/5879 \\ 
0 &0 &0 &-446/10487 &0 \\ 
0 &0 &0 &446/10487 &-158/7263 \\ 
0 &0 &0 &0 &-2699/54280 \\ 
0 &0 &0 &-446/10487 &-194/2497 \\ 
0 &0 &0 &-361/4426 &811/8786 \\ 
0 &0 &0 &0 &811/8786 \\ 
0 &0 &0 &361/4426 &811/8786 \\ 
0 &0 &0 &446/10487 &-194/2497 \\ 
0 &0 &0 &0 &-2699/54280 \\ 
0 &0 &0 &-446/10487 &-158/7263 \\ 
446/10487 &0 &-446/10487 &-23/1712 &0 \\ 
-194/2497 &-2699/54280 &-158/7263 &0 &-287/6896 \\ 
0 &686/18395 &1345/18033 &0 &0 \\ 
-31/13608 &31/6804 &-31/13608 &0 &-94/5141 \\ 
-293/5856 &293/2928 &-293/5856 &0 &199/7414 \\ 
195/1351 &195/1351 &195/1351 &0 &0 \\ 
5143/432011 &-10286/432011 &5143/432011 &0 &0 \\ 
-1262/6927 &1217/3340 &-1262/6927 &0 &0 \\ 
1022/16221 &31/2753 &-251/6200 &0 &67/10675 \\ 
357/604 &401/8908 &-731/1459 &0 &-27/6463 \\ 
646/2335 &961/2904 &3219/8357 &0 &63/14102 \\ 
222/2929 &120/917 &3951/21250 &0 &56/50269 \\ 
-338/963 &0 &338/963 &-90/2261 &0 \\ 
-883/3115 &847/1494 &-883/3115 &0 &-65/1076 \\ 
-152/7003 &304/7003 &-152/7003 &0 &144/3743 \\ 
-413/985 &-526/2321 &-203/5977 &0 &199/1660 \\ 
-89/1024 &-737/2869 &-213/499 &0 &62/2495 \\ 
 \end{matrix} \right]. $$ 

 \newpage
 The 21 to 25 columns of $\mathbf G^{-1}$ are
  $$\left[ \begin{matrix} 
 0 &0 &0 &446/10487 &0 \\ 
0 &0 &0 &-361/4426 &-407/5879 \\ 
0 &0 &0 &446/10487 &1345/18033 \\ 
0 &0 &0 &0 &686/18395 \\ 
0 &0 &0 &0 &-407/5879 \\ 
0 &0 &0 &0 &686/18395 \\ 
0 &0 &0 &-446/10487 &1345/18033 \\ 
0 &0 &0 &361/4426 &-407/5879 \\ 
0 &0 &0 &-446/10487 &0 \\ 
0 &0 &0 &446/10487 &-158/7263 \\ 
0 &0 &0 &0 &-2699/54280 \\ 
0 &0 &0 &-446/10487 &-194/2497 \\ 
0 &0 &0 &-361/4426 &811/8786 \\ 
0 &0 &0 &0 &811/8786 \\ 
0 &0 &0 &361/4426 &811/8786 \\ 
0 &0 &0 &446/10487 &-194/2497 \\ 
0 &0 &0 &0 &-2699/54280 \\ 
0 &0 &0 &-446/10487 &-158/7263 \\ 
446/10487 &0 &-446/10487 &-23/1712 &0 \\ 
-194/2497 &-2699/54280 &-158/7263 &0 &-287/6896 \\ 
0 &686/18395 &1345/18033 &0 &0 \\ 
-31/13608 &31/6804 &-31/13608 &0 &-94/5141 \\ 
-293/5856 &293/2928 &-293/5856 &0 &199/7414 \\ 
195/1351 &195/1351 &195/1351 &0 &0 \\ 
5143/432011 &-10286/432011 &5143/432011 &0 &0 \\ 
-1262/6927 &1217/3340 &-1262/6927 &0 &0 \\ 
1022/16221 &31/2753 &-251/6200 &0 &67/10675 \\ 
357/604 &401/8908 &-731/1459 &0 &-27/6463 \\ 
646/2335 &961/2904 &3219/8357 &0 &63/14102 \\ 
222/2929 &120/917 &3951/21250 &0 &56/50269 \\ 
-338/963 &0 &338/963 &-90/2261 &0 \\ 
-883/3115 &847/1494 &-883/3115 &0 &-65/1076 \\ 
-152/7003 &304/7003 &-152/7003 &0 &144/3743 \\ 
-413/985 &-526/2321 &-203/5977 &0 &199/1660 \\ 
-89/1024 &-737/2869 &-213/499 &0 &62/2495 \\ 
 \end{matrix} \right]. $$ 

 \newpage
Column 26 to 30 of $\mathbf G^{-1}$ are
 $$\left[ \begin{matrix} 
 5143/432011 &-3599/5453 &0 &0 &-99/90076 \\ 
0 &0 &0 &0 &0 \\ 
-5143/432011 &937/6487 &-137/1263 &336/1285 &-1471/2543 \\ 
-10286/432011 &-115/589 &2539/11975 &609/3007 &513/7487 \\ 
0 &0 &0 &0 &0 \\ 
10286/432011 &51/881 &723/6241 &-162/1069 &-809/2711 \\ 
5143/432011 &739/2742 &5078/11975 &625/1543 &255/1846 \\ 
0 &0 &0 &0 &0 \\ 
-5143/432011 &-96/3349 &946/2781 &-2033/3601 &-77/4190 \\ 
1262/6927 &491/831 &-1395/8672 &-341/1901 &524/3107 \\ 
-1217/3340 &551/2084 &-823/16086 &-3234/139709 &-572/2423 \\ 
1262/6927 &-169/2723 &301/5142 &461/3464 &-1841/2873 \\ 
0 &0 &0 &0 &0 \\ 
0 &0 &0 &0 &0 \\ 
0 &0 &0 &0 &0 \\ 
-1262/6927 &1022/16221 &357/604 &646/2335 &222/2929 \\ 
1217/3340 &31/2753 &401/8908 &961/2904 &120/917 \\ 
-1262/6927 &-251/6200 &-731/1459 &3219/8357 &3951/21250 \\ 
0 &0 &0 &0 &0 \\ 
0 &67/10675 &-27/6463 &63/14102 &56/50269 \\ 
0 &-74/64347 &-111/20717 &-45/9634 &95/18561 \\ 
0 &53/12731 &-68/6223 &-45/20756 &152/19417 \\ 
0 &-90/9571 &21/29873 &-61/5987 &37/11708 \\ 
-149/3989 &0 &0 &0 &0 \\ 
0 &-49/1496 &110/8839 &81/1768 &145/3054 \\ 
-2/15 &429/14929 &-61/5587 &-91/2264 &-142/3409 \\ 
429/14929 &-797/5375 &107/6728 &56/3313 &182/5961 \\ 
-61/5587 &107/6728 &-297/2248 &-112/8103 &8/11573 \\ 
-91/2264 &56/3313 &-112/8103 &-425/2513 &-175/4867 \\ 
-142/3409 &182/5961 &8/11573 &-175/4867 &-1250/7501 \\ 
0 &0 &0 &0 &0 \\ 
0 &-332/6139 &59/9990 &-285/4963 &92/5557 \\ 
0 &131/5963 &-80/1293 &-155/11002 &277/6160 \\ 
0 &106/3253 &-178/6237 &87/4607 &179/15250 \\ 
0 &632/47779 &239/9859 &216/7151 &-429/16265 \\ 
 \end{matrix} \right]. $$

 \newpage
Column 31 to 35 of $\mathbf G^{-1}$ are
  $$\left[ \begin{matrix} 
 -338/963 &152/7003 &-883/3115 &89/1024 &-413/985 \\ 
-265/724 &181/7664 &-223/723 &884/2323 &-469/1397 \\ 
-338/963 &152/7003 &-883/3115 &213/499 &-203/5977 \\ 
0 &-304/7003 &847/1494 &737/2869 &-526/2321 \\ 
0 &-181/3832 &446/723 &884/2323 &-469/1397 \\ 
0 &-304/7003 &847/1494 &737/2869 &-526/2321 \\ 
338/963 &152/7003 &-883/3115 &213/499 &-203/5977 \\ 
265/724 &181/7664 &-223/723 &884/2323 &-469/1397 \\ 
338/963 &152/7003 &-883/3115 &89/1024 &-413/985 \\ 
-338/963 &-883/3115 &-152/7003 &-203/5977 &-213/499 \\ 
0 &847/1494 &304/7003 &-526/2321 &-737/2869 \\ 
338/963 &-883/3115 &-152/7003 &-413/985 &-89/1024 \\ 
-265/724 &-223/723 &-181/7664 &-469/1397 &-884/2323 \\ 
0 &446/723 &181/3832 &-469/1397 &-884/2323 \\ 
265/724 &-223/723 &-181/7664 &-469/1397 &-884/2323 \\ 
-338/963 &-883/3115 &-152/7003 &-413/985 &-89/1024 \\ 
0 &847/1494 &304/7003 &-526/2321 &-737/2869 \\ 
338/963 &-883/3115 &-152/7003 &-203/5977 &-213/499 \\ 
-90/2261 &0 &0 &0 &0 \\ 
0 &-65/1076 &144/3743 &199/1660 &62/2495 \\ 
0 &144/3743 &65/1076 &62/2495 &-199/1660 \\ 
0 &-17/3218 &573/3355 &298/4075 &-217/3070 \\ 
0 &573/3355 &17/3218 &-217/3070 &-298/4075 \\ 
0 &0 &0 &0 &0 \\ 
0 &0 &0 &0 &0 \\ 
0 &0 &0 &0 &0 \\ 
0 &-332/6139 &131/5963 &106/3253 &632/47779 \\ 
0 &59/9990 &-80/1293 &-178/6237 &239/9859 \\ 
0 &-285/4963 &-155/11002 &87/4607 &216/7151 \\ 
0 &92/5557 &277/6160 &179/15250 &-429/16265 \\ 
-202/581 &0 &0 &0 &0 \\ 
0 &-1189/1092 &0 &473/1203 &356/931 \\ 
0 &0 &-1189/1092 &-356/931 &473/1203 \\ 
0 &473/1203 &-356/931 &-1321/1514 &0 \\ 
0 &356/931 &473/1203 &0 &-1321/1514 \\ 
 \end{matrix} \right]. $$

\bibliographystyle{siamplain}
\bibliography{references}
\end{document}